\documentclass[a4paper,12pt]{article}
\usepackage{amssymb,amsmath,amscd,amsfonts,amsthm,mathtools}
\usepackage{latexsym,graphpap,xcolor}
\usepackage[utf8]{inputenc}
\usepackage[english]{babel}
\usepackage{color}

\usepackage[right=3cm,left=3cm,top=3cm,bottom=3cm]{geometry}

\newtheorem{theorem}{Theorem}[section]
\newtheorem{prop}[theorem]{Proposition}
\newtheorem{lemme}[theorem]{Lemma}
\newtheorem{coro}[theorem]{Corollary}

\newtheorem{defi}[theorem]{Definition}

\newenvironment{demo}{ \noindent 
\emph{\textbf{Proof:}}}{\hfill$\square$\\}

\newcommand{\RR}{\mathbb{R}}
\newcommand{\NN}{\mathbb{N}}
\newcommand{\CC}{\mathbb{C}}
\newcommand{\ZZ}{\mathbb{Z}}

\newcommand{\Cc}{\mathcal{C}}
\newcommand{\Ec}{\mathcal{E}}
\newcommand{\Fc}{\mathcal{F}}
\newcommand{\Kc}{\mathcal{K}}
\newcommand{\Lc}{\mathcal{L}}

\newcommand{\ul}{\mathrm{ul}}
\newcommand{\uls}{\mathrm{ul-s}}
\newcommand{\ulw}{\mathrm{ul-w}}
\newcommand{\ulx}{\mathrm{ul-*}}
\newcommand{\loc}{\mathrm{loc}}

\newcommand{\s}{{\text{\rm s}}}
\newcommand{\w}{{\text{\rm w}}}
\renewcommand{\d}{\,{\rm d}}
\newcommand{\id}{\,{\rm id}}

\newcommand{\grad}{\nabla}
\newcommand{\la}{\langle}
\newcommand{\ra}{\rangle}
\newcommand{\doublehookrightarrow}{\mathrel{\mathrlap{{\mspace{4mu}\lhook}}{\hookrightarrow}}}
\renewcommand{\div}{\operatorname{div}}
\newcommand{\Un}{1\hspace{-1.5mm}1}
\newcommand{\no}{n$^{\text{o}}$}

\newdimen\texpscorrection
\newdimen\figcenter
\def\figurewithtex #1 #2 #3 #4 #5\cr{\null
  {\goodbreak\figcenter=\hsize\relax
  \advance\figcenter by -#4truecm
  \divide\figcenter by 2
  \begin{figure}[hbt]
  \vskip #3truecm\noindent\hskip\figcenter
  \includegraphics{#1}{\hskip\texpscorrection\input #2 }
  \vskip 0.8truecm{\baselineskip=0.8\baselineskip
  \noindent \vbox{\noindent {\footnotesize #5}}\par}
  \end{figure}}}
\def\point#1 #2 #3 {\rlap{\kern #1 truecm
\raise #2 truecm \hbox{#3}}}

\numberwithin{equation}{section}

\begin{document}

\title{\bf Parabolic abstract evolution equations in cylindrical domains and uniformly local Sobolev spaces}

\author{Romain {\sc Joly}}

\maketitle

\bigskip

\begin{abstract}
In this article, we consider parabolic equations of the form 
$$\partial_t u(x,t)=\Delta u(x,t) - Bu(x,t) + F(u(x,t))$$
where $u$ is valued in a transverse Hilbert space $Y$ and $B$ is a positive self-adjoint operator on $Y$, which allows for a different diffusion mechanism in the transverse direction. Our aim is to consider solutions with infinite energy and we will study the Cauchy problem in the uniformly local spaces associated with the norm 
$$\|u\|_{L^2_\ul(\RR,Y)}= \sup_{a\in\RR} \|u\|_{L^2(B(a,1),Y)}.$$
For the classical parabolic equation, i.e. when $Y=\RR$, it is known that the Cauchy problem is ill-posed in the weak version of the uniformly local spaces but well-posed in a stronger version, where additional uniform continuity is  required. 
In this paper, we show that the linear operator $\partial^2_{xx} - B$ is not necessarily a sectorial operator in any version of the uniformly local Lebesgue space, due to the possible non-density of its domain. We then use the theory of parabolic abstract evolution equations to formulate a well-posed Cauchy problem, even in the weak version of the uniformly local space. 
In particular, we believe that this paper offers a new perspective on the comparison between the two versions of the uniformly local spaces and also provides a new natural example of differential operators with non-dense domain.\\[8mm]
\end{abstract}

\pagebreak

\section{Introduction and motivation}

The modelling of many physical, chemical or biological phenomena leads to semilinear parabolic PDEs. Several interesting patterns correspond to solutions with unbounded energy, such as non-zero constant solutions and travelling fronts. When studying different models, it becomes apparent that many arguments are the same, making it tempting to bring together all the studied models in one abstract framework. Our initial motivation was to properly set the Cauchy problem for one such abstract framework. This raises some interesting questions of fundamental analysis, making this study interesting in itself. 

In this introduction, we would like to motivate our study and to explain what we believe are its most interesting aspects. 
To keep the discussion as simple as possible, we have postponed the technical details, more general results and more precise discussions to the other sections. Nevertheless, this paper will not aim for the utmost generality.
Our aim is to provide a better understanding of the subject and introduce some tools and ideas for possible future applications. 

\subsection{The cylindrical framework}

Our first motivation was to introduce an abstract framework modeling several cylindrical geometries. To explain our formalism, we illustrate it with the following simplified version. Let $Y$ be a separable Hilbert space, embedded with the scalar product $\la\cdot|\cdot\ra_Y$ and its associated norm $|\cdot|_Y$, and let $B:D(B)\rightarrow Y$ be a positive self-adjoint operator of dense domain $D(B)\subset Y$. We consider an abstract parabolic equation of the type
\begin{equation}\label{eq_intro}
\partial_t u(x,t)=\partial_{xx}^2 u(x,t) - Bu(x,t) + L(u,\partial_x u) + F(u(x,t))~~t>0, x\in\RR
\end{equation}
where $F:Y\rightarrow Y$ is a nonlinearity and $L$ includes linear terms of lower order. Several interesting classes of models fit into this framework:
\begin{itemize}
\item {\it Advective reaction-diffusion equations in cylinders.} We consider a cylinder $\Omega=\RR\times\omega\subset\RR^3$, where $\omega\subset\RR^2$ is a bounded smooth domain. A point of $\Omega$ is denoted by $(x,y)\in\RR\times\omega$. We set $Y=L^2(\omega)$ and $B=- \Delta_y$ with Dirichlet boundary condition on $\partial\omega$. Let $v\in\Cc^1(\omega,\RR^2)$ be a smooth vector field and $L=v\cdot\grad_y$. The function $F:L^2(\omega)\rightarrow L^2(\omega)$ can simply be a Nemytskii operator associated to a real function $f\in\Cc^0(\omega\times\RR,\RR)$.  
Then our abstract equation \eqref{eq_intro} becomes a reaction-diffusion equation with convective term
$$\partial_t u(x,y,t)=\Delta u(x,y,t) + v(y)\cdot\nabla_y u(x,y,t) + f(y,u(x,y,t))$$ 
(actually, our framework will allow more freedom concerning the $x$ dependence but we keep the notations simple in this introduction). This kind of PDE is very common, see for example \cite{Berestycki_Nirenberg, Giletti, Muratov_Novaga_1, Muratov_Novaga_2, Panfilov, Pang_Wu, Xin}. These articles are particularly interested in front-like solutions that do not belong to $L^p(\RR\times\omega)$ for $p<+\infty$. The Cauchy problem can be studied in many frameworks. It seems that this paper is the first study in the uniformly local Sobolev spaces.

\item {\it Favorable corridor in unfavorable environment.} Consider $Y=L^2(\RR)$ and $B=-\Delta_y+\id$. Assume that $F(u)$ is given by $y\mapsto f(y,u(y))$ where $f(y,u)u\leq 0$ for $|y|\geq 1$. We obtain a classical reaction-diffusion equation 
$$\partial_t u(x,y,t)=\Delta u(x,y,t) - u(x,y,t) + f(y,u(x,y,t))~~,~~~(x,y)\in\RR^2.$$
The cylindrical framework follows from the fact that the reaction part $-\id+f$ is unfavorable outside the cylinder $\RR\times [-1,1]$. If it is sufficiently favorable to the propagation of $u$ inside this cylinder, we expect the existence of traveling fronts describing the propagation of $u$ in the $x$-direction as in Figure \ref{fig-corridor}. Notice that these fronts should not vanish outside the cylinder $\RR\times [-1,1]$ but they are expected to have a fast decay as $y$ goes to $\pm\infty$. In contrast, as $x$ goes to $-\infty$ or $+\infty$, the front should converge to a non-zero profile. This leads to one of the main concerns of this article: to be able to set a Cauchy problem for solutions with finite energy in the transverse variable $y$ (for each $x$, $y\mapsto u(x,y)$ belongs to $L^2(\RR)$) but infinite energy in the longitudinal variable ($|u(x,\cdot)|_Y$ does not go to zero as $x\rightarrow \pm\infty$), see Figure \ref{fig-corridor}.
In this case, we can say that the anisotropy associated to the ``cylindrical geometry" follows from the difference of topology in each direction: the uniformly local spaces with respect to $x$ and the classical $L^2$ with respect to $y$.

We can obviously apply this type of model to consider the propagation of a biological species in a corridor shaped environment. Less obviously, this also models the propagation of spreading depressions, a possibly harmful depolarization of neurons following strokes or migraines, see \cite{Chapuisat1,Chapuisat2,Chapuisat3}. In this case, the ``favorable'' environment is the layer of gray matter and it is more realistic to consider $x\in\RR^2$ to model this layer. Higher dimensional spaces will be allowed in our more general setting after this introduction.

\begin{figure}[ht]
\begin{center}
\includegraphics[width=8cm]{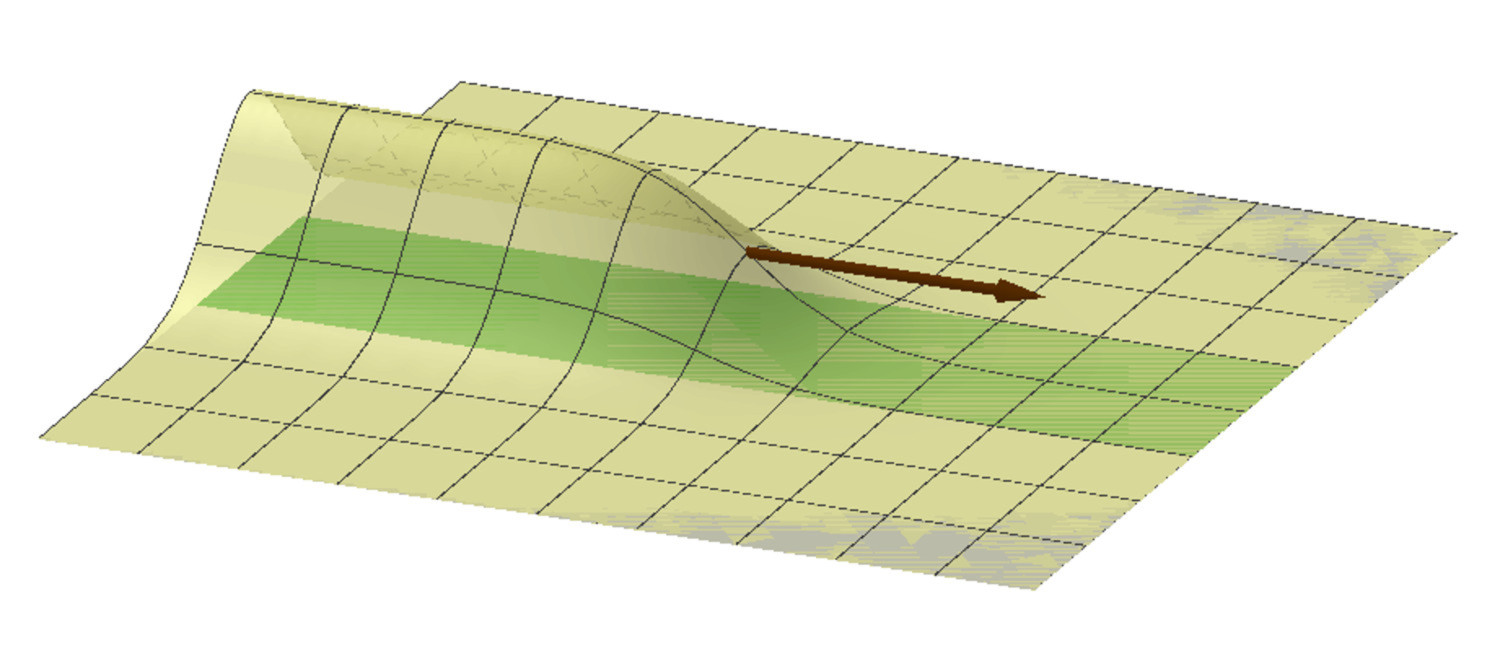} 
\end{center}
\caption{\it A front-like solution is invading a corridor of favorable environment surrounded by an unfavorable one. We aim at considering a general Cauchy problem including this type of solutions, which have finite energy in each transverse section but infinite energy along the propagation variable.}\label{fig-corridor}
\end{figure}

\item {\it Anisotropic diffusion.} We can also consider operators $B$ that are not operators of the Laplacian type, for example $B=(-\Delta_y)^{\sigma/2}$ is possible for any $\sigma>0$. In this case, the transverse diffusion is of different nature, for example generated by a bilapacian operator if $\sigma=4$. For $\sigma\in(0,2)$, such $B$ corresponds to an ``anomalous diffusion'' and is related to L\'evy processes, see \cite{Bertoin} for example. This type of anomalous diffusion appears in several models, from mechanics to finance, we refer to the introduction of \cite{anomalous} for a list of many applications. Being able to consider different types of diffusion shows the flexibility of our abstract framework.

\item {\it Parabolic systems.} A bounded operator $B$ can be artificially added in or separated from $F$. Thus, up to introducing an artificial $B=\id$, we can include  parabolic systems in our framework by choosing $Y=\RR^d$. Then, the fact that $u(x)$ is valued in $Y$ simply means that $u(x)$ is a vector with $d$ component. More interestingly, we can also use a transverse variable $y=\lambda$ as a parameter, for example a genetic trait or the age of individuals in a biological model. Set $Y=L^2([0,1])$ and, up to change $F$ in $F+\id$, take $B=\id$. We can consider for example a nonlinearity $F(u)$ of the form $\lambda\mapsto f(\lambda,u(\lambda))-u(\lambda)+\int_0^1 K(\lambda,\lambda')u(\lambda')\d \lambda'$, leading to a system of the type
$$\partial_t u(x,\lambda,t)=\partial^2_{xx} u(x,\lambda,t) + f(\lambda,u(x,\lambda,t))+\int_0^1 K(\lambda,\lambda')u(x,\lambda',t)\d \lambda'.$$
The kernel $K$ is not necessarily symmetric (since included in $F$ and not in $B$) and may describe the mutations of the trait. Even if the boundness of $B$ makes this case much simpler, it seems that the status of the Cauchy problem in uniformly local Sobolev spaces was still unclear until now.
\end{itemize}

\subsection{The linear semigroup in $L^2(\RR,Y)$}

Before considering the uniformly local Sobolev spaces, we can simply study the linear part of Equation \eqref{eq_intro} in $L^2(\RR,Y)$. Even if this context is very classical in all the above concrete examples, where the main linear operator is a Laplacian operator, it is less obvious when we consider the linear abstract operator $-\partial^2_{xx}+B$.
\begin{prop}\label{prop_intro1}
Let $B:D(B)\rightarrow Y$ be a positive self-adjoint operator with dense domain. 
The operator $-\partial^2_{xx}+B$ is a well-defined self-adjoint operator from $D(-\partial_{xx}^2+B)=H^2(\RR,Y)\cap L^2(\RR,D(B))$ into $L^2(\RR,Y)$. In particular, it generates an analytic semigroup $e^{(\partial_{xx}^2-B)t}$ on $L^2(\RR,Y)$.
\end{prop}
The above result seems completely natural, but it is more subtle that it may look. Actually, the crucial point concerns the domain: its is clear that $(-\partial^2_{xx}+B)u$ belongs to $L^2(\RR^d,Y)$ if $u$ belongs to $H^2(\RR^d,Y)\cap L^2(\RR^d,D(B))$. However, it could be a priori possible that $(-\partial^2_{xx}+B)u$ belongs to $L^2(\RR^d,Y)$ without having both $\partial^2_{xx}u$ and $Bu$ in $L^2(\RR^d,Y)$. Even in  the case of the Laplacian operator, i.e. $Y=L^2(\RR)$ and $B=-\partial_{yy}^2$, proving the closedness of the operator in $H^2(\RR^2)$ uses a non trivial characterization of $H^1$ via a continuity estimate of the translation. When working in a general space $Y$ instead of $\RR$, a similar characterization holds if $Y$ satisfies the Radon-Nikodym property, see Section \ref{section_RN}. In our Hilbertian framework, $Y$ satisfies this property, but this is not the case if $Y=L^1(\RR)$ for example. 

This problem belongs to a broader class of questions: if $A$ and $B$ are closed operators, is $(A+B)$ a closed operator with domain $D(A)\cap D(B)$? If not, is it at least closable? This type of questions is related to the question of ``maximal regularity'' or ``mixed derivative'' and $A$ is often viewed as a differential operator with respect to the time, so that the whole PDE is $(A+B)u=0$. There have been a lot of works studying this kind of problems, in particular the 
the pioneer works of da Prato and Grisvard \cite{da-Prato-Grisvard} and of Dore and Venni \cite{DV}, see also \cite{HvNvW2, Kalton-Lancien, Kalton-Weis,Roidos}. The frameworks of these articles are much more general than ours and they use more involved concepts as $H^\infty-$calculus, UMD spaces or $R-$sectorial operators. 
Some counter-examples have been constructed, underlining how the question of maximal regularity can be delicate, see \cite{BC,Fackler2,Venni}.

\subsection{The heat equation in the uniformly local Sobolev spaces}
The above mentioned difficulty may seem artificial since it appears due to our choice of working in an abstract framework. On the opposite, working with uniformly local Sobolev spaces leads to several issues, even in the simplest concrete case of the heat equation. Let us explicit them before considering the general abstract parabolic case. 

We first recall that we would like to consider solutions with possibly unbounded energy in the $x$ direction. Then, it is natural to consider the following norm
\begin{equation}\label{def-norme-intro}
\|u\|_{L^2_\ul(\RR,Y)}= \sup_{a\in\RR} \left(\int_{a}^{a+1} |u(x)|_Y^2 \d x\right)^{1/2}
\end{equation}
and to introduce the ``weak'' uniformly local Lebesgue space $L^2_\ulw(\RR,Y)$ defined as the space of functions with finite norm $\|\cdot\|_{L^2_\ul(\RR,Y)}$.  We can also consider the associated Sobolev spaces, see Section \ref{section_ul} for details and additional discussions. Spaces as $L^2_\ulw$ are sometimes called ``locally uniform'' and have been first considered by Kato in \cite{Kato}. For the readers that are not familiar with this kind of spaces, we can motivate them by the following remarks:
\begin{itemize}
\item The uniformly local Sobolev spaces contain solutions as traveling fronts or more generally states with different limits when $x\rightarrow\pm\infty$. These are very important patterns when studying diffusion and invasion phenomena but none has finite energy and none is contained in a $L^p$-space for $p<+\infty$.
\item The $L^\infty$-topology is sensitive to very localized perturbations, even ones with small energy and it is not very relevant in diffusion phenomena. Also notice that a space as $(\Cc^0_b(\RR),\|\cdot\|_\infty)$ does not contain functions of Heaviside type, that are common for initial data when studying invasion phenomena. Anyway, we will see that working with the $L^\infty-$topology brings the same difficulties as working in the uniformly local Sobolev spaces.
\item When considering parabolic equation, the symmetry of the Laplacian operator or the existence of an energy are important features. Working with the norm \eqref{def-norme-intro} enables to keep track of the $L^2-$structure while opening the possibility of solutions with infinite energy.
\end{itemize}

\begin{figure}[ht]
\begin{center}
\includegraphics[width=\textwidth]{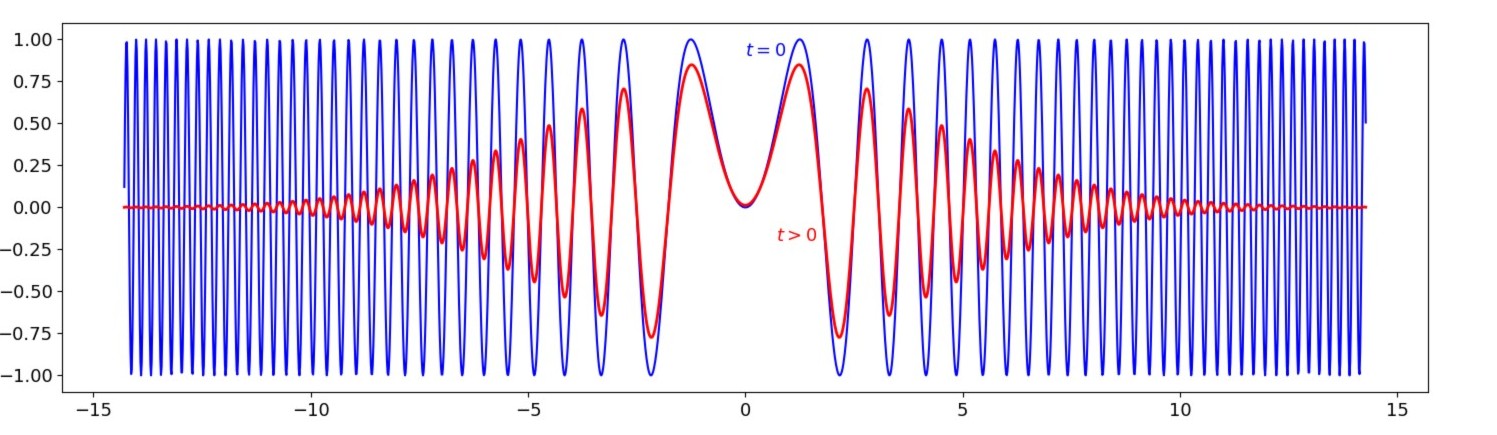} 
\end{center}
\caption{\it The graphic shows the initial data $\sin(x^2)$ (blue) and its evolution by the heat flow after some time (red). As soon as $t$ is positive, the fast oscillations of $\sin(x^2)$ are uniformized and the solution of the heat equation jumps at distance $1$ from the initial data in the $L^\infty$-topology. In this sense, the heat equation is ill-posed in this topology.}\label{fig-heat1}
\end{figure}

\begin{figure}[ht]
\begin{center}
\includegraphics[width=\textwidth]{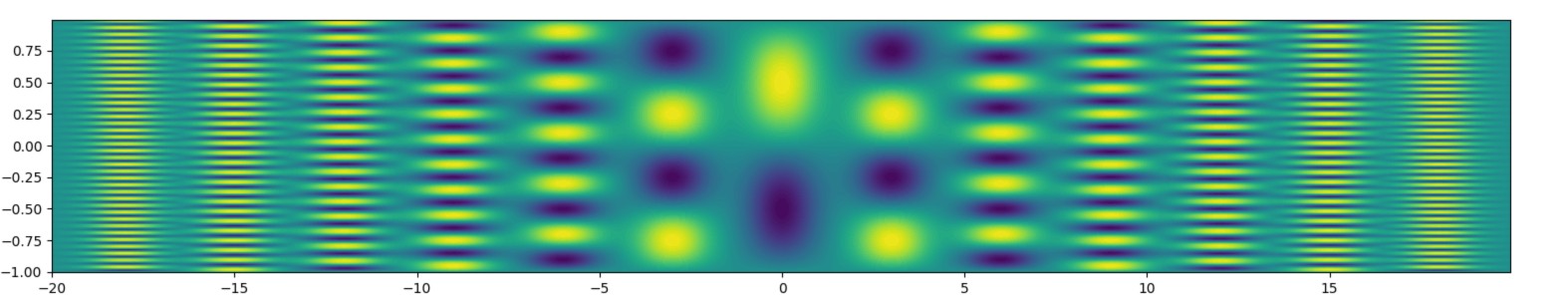} 
\end{center}
\caption{\it An example of initial data in a cylinder that would lead to the same problem as illustrated in Figure \ref{fig-heat1}. Notice that the horizontal dependence is uniformly smooth, so this data belongs to the strong version $L^2_\uls(\RR,L^2((-1,1)))$ of the uniformly local Lebesgue space, valued in the transverse section space $L^2(-1,1)$. Hence, when considering parabolic equations in vector-valued uniformly local space, using the strong version does not help.}\label{fig-heat2}
\end{figure}

These remarks motivate the main purpose of the present paper: setting a Cauchy problem for \eqref{eq_intro} in $L^2_\ulw(\RR,Y)$. To understand the difficulty behind this purpose, we discuss the simplest case: the heat equation $\partial_t u =\Delta u$.
\begin{enumerate}
\item {\it The heat equation is ill-posed in $L^2_\ulw(\RR,\RR)$}. Indeed, an initial data as $u_0(x)=\sin(x^2)$ exhibits very fast oscillation at infinity. These oscillations are homogenized by the heat semigroup (defined via the convolution by its kernel) and, as soon as $t$ is positive, $u(t)=e^{\Delta t}u_0$ must converge to $0$ as $x$ goes to $\pm\infty$, see Figure \ref{fig-heat1}. Thus $\|u(t)-u_0\|_{L^2_\ul}$ does not converge to zero as $t$ goes to $0^+$ and the solution fails to be continuous in the uniformly local space. Notice that the problem does not come from our ``fancy'' choice of functional space: the same initial data shows that the heat equation is ill-posed in $L^\infty(\RR,\RR)$ or in $\Cc^0(\RR,\RR)$.
\item {\it The heat equation is well-posed in $L^2_\uls(\RR,\RR)$.} The space $L^2_\uls(\RR,\RR)$ is a stronger version of the uniformly local space where the translation $\xi\mapsto u-u(\cdot-\xi)$ is continuous for all functions $u$, see Section \ref{section_ul}. This excludes the problematic $x \mapsto\sin(x^2)$ and can be seen as replacing $\Cc^0(\RR,\RR)$ by the space of uniformly continuous functions.
It is known that the heat semigroup is well defined and continuous in $L^2_\uls(\RR,\RR)$, as shown by Arrieta, Cholewa, Dlotko and Rodr\'\i{}guez-Bernal in \cite{Arrieta1}. Other PDEs have been considered, see for example \cite{Arrieta2,HIOS,Matos-Souplet,Michalek,Suguro}, still in similar strong versions of the uniformly local spaces. Until now, choosing this strong version $L^2_\uls$ has been considered as the right way to solve the Cauchy problem. 
\item {\it The heat equation is ill-posed in $L^2_\uls(\RR,L^2((-1,1)))$.} In our cylindrical setting, choosing the strong version of the uniformly local spaces is not sufficient to define a proper Cauchy problem, even for the heat equation in $\RR\times(-1,1)$. Indeed, we can use the oscillations in the $y$ direction to construct an initial data for which the heat semigroup is not continuous at $t=0^+$, see Figure \ref{fig-heat2}. 
This shows that the strong version of the uniformly local spaces is not so helpful in a general context of vector valued functions.
\end{enumerate}
To our point of view, one of the main interests of the present paper will be to clarify the issues raised in the above remarks. 


\subsection{The linear operators: a problem of density of domains}

Working in an abstract framework, where we cannot use heat kernels or other specific tools, forces us to understand the fundamental problem behind the above facts concerning the heat equation. It turns out that the ill-posedness is related to the lack of density of the domain of the linear operator. To our knowledge, this is the first time that this fact is explicitly stated. 

To the ill-posedness of the heat equation on the weak version of the uniformly local Sobolev spaces, we can associate this more general result.
\begin{prop}\label{prop_intro22}
The operator 
\[(-\partial_{xx}^2+B):H^2_\ulw (\RR,Y)\cap L^2_\ulw (\RR,D(B))\longrightarrow L^2_\ulw (\RR,Y)\]
is a closed operator satisfying sectorial-like resolvent estimates, but its domain is not dense. There exist functions $u_0\in L^2_\ulw(\RR,Y)$ such that the associated semigroup trajectory $e^{(\partial_{xx}^2-B)t}u_0$ is not continuous at $t=0$ in $L^2_\ulw(\RR,Y)$.
\end{prop}
As a generalization of the result of \cite{Arrieta1}, we prove that the strong version of the uniformly local spaces yields the well-posedness of \eqref{eq_intro}, but only for bounded $B$.  
\begin{prop}\label{prop_intro23}
If $B$ is a bounded operator on $Y$, then the operator 
\[(-\partial_{xx}^2+B):H^2_\uls (\RR,Y)\cap L^2_\uls (\RR,D(B))\longrightarrow L^2_\uls (\RR,Y)\]
is a sectorial operator which generates an analytic semigroup on $L^2_\uls(\RR,Y)$.
\end{prop}
However, for general $B$, the strong version is not more helpful.
\begin{prop}\label{prop_intro2}
If $Y$ is infinite-dimensional and if $B$ has compact resolvent, the operator 
\[(-\partial_{xx}^2+B):H^2_\uls (\RR,Y)\cap L^2_\uls (\RR,D(B))\longrightarrow L^2_\uls (\RR,Y)\]
is a closed operator satisfying sectorial-like resolvent estimates, but its domain is not dense. There exist functions $u_0\in L^2_\uls(\RR,Y)$ such that the associated semigroup trajectory $e^{(\partial_{xx}^2-B)t}u_0$ is not continuous at $t=0$ in $L^2_\uls(\RR,Y)$.
\end{prop}
In some sense, the three propositions above contain most of the interest of the present paper: a new perspective on the comparison between both versions of the uniformly local spaces and a new natural example of differential operators with non-dense domain.

Once we have identified that the problem comes from the lack of density of the domain,
we need a suitable theory to set a proper Cauchy problem. This leads us back to the theory 
introduced by Sinestrari and da Prato in \cite{da-Prato-Sinestrari,Sinestrari}, where the authors have introduced a suitable notion of ``integrated semigroup'' generated by such non-densely defined operators. Related Cauchy problems for quasilinear parabolic equation have been studied by Lunardi in \cite{Lunardi}. Recently, the interest on this concept has raised due to applications in age structured biological models, see the book \cite{Magal-Ruan} of Magal and Ruan and the associated articles as \cite{Ducrot-et-al,Ducrot-et-al2,Magal-Seydi-Wang}.

\subsection{Examples of statements for the Cauchy problem}

In the context of analytic semigroups and parabolic equations, it is usual to consider a nonlinearity that could be defined not on $X$ but on a fractional space $D(A^\alpha)$. In particular, we can consider polynomial nonlinearities that are not Lipschitz continuous on $X$. In the present paper, we would like to consider similar PDEs but with our operators with non-dense domain. This requires some study to characterize the ``intermediate spaces'' related to these operators and to study nonlinear functions on them. Once this is done, the theory of integrated semigroups of Sinestrari, da Prato and Lunardi can be directly applied. To give a quick example of our aim, we state the following result in the simplified context of \eqref{eq_intro}. 
\begin{prop}\label{prop_intro3}
Let $\alpha\in [0,1)$ and let $F:D(B^\alpha)\rightarrow Y$ be a function satisfying
\[\forall u,v\in D(B^\alpha)~,~~|F(u)-F(v)|_Y \leq C (1+|u|_{D(B^\alpha)}^\gamma+|v|_{D(B^\alpha)}^\gamma) |u- v|_{D(B^\alpha)}\]
with $\gamma\geq 0$ such that $\alpha+\frac {\gamma}{4(1+\gamma)} < 1$. For any $u_0\in H^2_\ulw(\RR,Y) \cap L^2_\ulw (\RR,D(B))$, there exists a well-defined solution $u(t)$ of \eqref{eq_intro} with initial data $u(0)=u_0$ and in a time interval $[0,T)$. It satisfies the equation in the strong sense and converges to $u_0$ in $L^2_\ulw(\RR,Y)$, when $t\rightarrow 0$. Moreover, if $\gamma=0$ (the globally Lipschitz case), then we can take $T=+\infty$.
\end{prop}
One of the interests to work in the uniformly local space $L^2_\ulw(\RR,Y)$ is to keep track of the $L^2-$structure, naturally associated with the Lyapunov structure of the parabolic equation and its related energy. To give an explicit application, assume that $F=-\grad_u V$ where $\grad_u V$ is the gradient of a coercive quadratic potential $V$ (see Section \ref{section_gradient} for precise definitions). We can associate to \eqref{eq_intro} an energy
$$\Ec(t)=\int_\RR \frac 12 |\partial_x u(x,t)|^2_Y + \frac 12 |B^{1/2}u(x,t)|^2_Y + V(u(x,t)) \d x.$$
The fact that this energy is non-increasing has important qualitative dynamical consequence, as precluding the blow-up of solutions. However, this is true in the classical context of solutions of finite energy. As already said, we are interested in solutions as fronts or invading patterns, for which the above energy is infinite. However, the framework of the uniformly local spaces as $H^1_\ulw(\RR,Y)\cap L^2_\ulw(\RR,D(B^{1/2}))$ enables to keep track of this formal Lyapunov structure, even for solutions for which $\Ec(t)$ is infinite. For example, we have the following result.
\begin{prop}\label{prop_intro4}
Consider the framework above with $F=-\grad_u V$ and assume in addition that $\alpha<1/2$. Then, the solution $u(t)$ defined in Proposition \ref{prop_intro3} exists for all times $t\geq 0$ and stays bounded in $H^2_\ulw(\RR,Y) \cap L^2_\ulw (\RR,D(B))$. 

If moreover, $B$ has compact resolvents, then for any sequence $(t_n)$ going to $+\infty$, there exists a subsequence $(t_{\varphi(n)})$ and a function $u_\infty\in H^2_\ulw (\RR,Y)\cap L^2_\ulw(\RR,D(B))$ such that, for all $M>0$ and $\beta\in [0,1)$,
\[u(t_{\varphi(n)})\xrightarrow[~~n\longrightarrow+\infty~~]{}u_\infty~~\text{ in }L^2((-M,M),D(B^\beta)).\]
\end{prop}

\subsection{Plan of the article}
The first sections of this article are mainly recalling known facts to introduce the basic notions and tools that are used in the paper. We provide there some quick proofs when we have not found suitable references or when we would like to explicit some arguments. But we are not claiming that these sections contain real novelties. 
In Section \ref{section_deux}, we  recall the definitions and basic properties of vector-valued Sobolev spaces, the points that are less important for the global understanding being postponed to Appendix \ref{section_spaces}. We also study the density of inclusions of the type $L^p(\RR,Y_2)\subset L^p(\RR,Y_1)$ when $Y_2$ is dense in $Y_1$ and $p\in [1,+\infty]$. In Section \ref{section_ul}, we introduce the uniformly local Sobolev spaces and provide some basic properties. Section \ref{section_nondense} is devoted to the theory of semigroups generated by non-densely defined operators, as introduced by Sinestrari, da Prato and Lunardi. We gather there the essential tools we will need later on.

Sections \ref{section_A} and \ref{section_Aul} construct the linear differential operator related to our framework. Section \ref{section_A} studies auxiliary differential operators in more classical $L^p$ spaces, yielding for example Proposition \ref{prop_intro1} above. Then, Section \ref{section_Aul} defines the linear operators in the uniformly local Sobolev spaces and integrates them into the theory of Sinestrari and da Prato. Propositions \ref{prop_intro22}, \ref{prop_intro23} and \ref{prop_intro2} above are corollaries of the results of this section. In this sense, Section \ref{section_Aul} may be seen as the chore of the present article.

Then, Section \ref{section_Cauchy} applies all the previous material to provide proper settings for the Cauchy problems, as the one stated in Proposition \ref{prop_intro3} for example. Section \ref{section_gradient} shows how to take advantage of the formal gradient structure and to obtain results as Proposition \ref{prop_intro4} even for solutions for which the natural parabolic energy is infinite. Section \ref{section_compact} simply recalls some compactness results concerning vector-valued Sobolev spaces to apply them to get convergences as the last statement of  Proposition \ref{prop_intro4}.

Finally, in Section \ref{section_appli}, we go back to simple applications of our abstract framework, as done in the beginning of this introduction. Doing so, we will briefly discuss the relevance of some of the hypotheses assumed in our abstract results.


\section{Lebesgue and Sobolev spaces of vector valued functions}\label{section_deux}

Our simple cylindrical framework consists in considering functions $(x,y)\mapsto u(x,y)$, whose sections $u(x,\cdot)$ belongs to a functional space $Y$ independent of $x$. It naturally leads to consider spaces as $L^2(\RR^d,Y)$, where $Y$ is a Banach space. This kind of vector valued spaces is usual, in particular in the analysis of PDE, considering time-dependent functions valued in a Banach space of spatially dependent functions. Their definitions are based on the Bochner integral and most of the basic properties of the classical spaces as $L^p(\RR^d,\RR)$ remain true. In this section, we simply recall the basic notations and enhance some particular results. We refer to Appendix \ref{section_spaces} for additional details.

\subsection{Basic definitions}\label{section_basic_spaces}
 Let $\Omega$ be a smooth open subset of $\RR^d$ and $Y$ a Banach space. Let $p\in [1,+\infty]$ and let $\rho\in\Cc^0(\Omega,\RR_+)$ be a non-negative weight. We introduce the Banach space $L^p_\rho(\Omega,Y)$ of measurable functions $u:\Omega\rightarrow Y$ such that the following associated norm is finite:
$$\|u\|_{L^p_\rho(\Omega,Y)}=\left(\int_\Omega \rho(x)|u(x)|_Y^p \d x\right)^{1/p}~~~~\text{if }p<\infty,$$
or 
$$\|u\|_{L^\infty_\rho(\Omega,Y)}=\sup_{x\in\Omega} \rho(x)|u(x)|_Y \d x~~~~\text{if }p=\infty.$$
In the ``flat'' case $\rho\equiv 1$, the subscript will be omitted.

The classical way to extend the notion of weak derivative is to say that $\partial_{x_i}u$ is the integrable function such that, for any $\varphi\in\Cc^\infty_c(\Omega,\RR)$,
$$\int_\Omega \partial_{x_i}u(x) \varphi(x) \d x ~= ~-~\int_\Omega u(x) \partial_{x_i}\varphi(x) \d x.$$
Notice that the test function is assumed to be real valued. The natural way to extend the notion of Sobolev spaces is to set 
\[W^{k,p}_\rho(\Omega,Y)=\{u:\Omega\rightarrow Y\text{ measurable }, \forall \nu\in\NN^d \text{ with }|\nu|\leq k~,~~ \partial_{x}^{\nu} u \in L^p_\rho(\Omega,Y)\},\]
where we use the multi-index notation $|\nu|=\nu_1+\ldots+\nu_d$ and $\partial_x^\nu=\partial_{x_1}^{\nu_1}\ldots \partial_{x_d}^{\nu_d}$, and to embed this space with the norm 
\[ \|u\|_{W^{k,p}_\rho(\Omega,Y)} = \left(\sum_{|\nu|\leq k} \|\partial_x^{\nu} u\|^2_{L^p_\rho(\Omega,Y)}\right)^{1/2}~. \]
If $Y$ is a Hilbert space, then $H^k_\rho(\Omega,Y):=W^{k,2}_\rho(\Omega,Y)$ is a Hilbert space associated with the scalar product
\[ \langle u|v\rangle_{H^k_\rho(\Omega,Y)} = \sum_{|\nu|\leq k} \int_\Omega \rho(x) \langle \partial_x^\nu u|\partial_x^\nu v\rangle_{Y}~.\]
It is also possible to define Sobolev spaces of fractional orders by using either the Fourier transform with respect to $x$, or the Sobolev-Slobodecki method, but we will not use this kind of spaces in this paper.

\subsection{Some results concerning dense embeddings}
Most of the classical results about dense embeddings of real-valued functions extend to the case of vector-valued functions, see Appendix \ref{section_spaces}. Indeed, the mollification technique is still valid. Let $\eta_\varepsilon=\frac{1}{\varepsilon^d}\eta(\cdot/\varepsilon)\in\Cc^\infty(\RR^d,\RR)$ be a smooth real-valued approximation of the identity with support in the ball $B_{\RR^d}(0,\varepsilon)$. We can check that, for any $u\in L^p(\RR^d,Y)$ with $p\in [1,+\infty)$, 
\[\big(\eta_\varepsilon \ast u\big)(x)~:=~\frac{1}{\varepsilon^d} \int_{\RR^d} \eta({x'}/ \varepsilon) u(x-x')\d x'\] 
is a smooth function belonging to $L^p(\RR^d,Y)$ and converging to $u$ in $L^p(\RR^d,Y)$ when $\varepsilon$ goes to $0$. This is shown as in the classical case, see for example \cite{Meyers-Serrin} for $Y=\RR$. For Bochner spaces, we refer for example to Section 1.3 of \cite{Arendt-Batty-et-al} or Section 4.2 of \cite{Kreuter} where the local regularity around any $x\in\RR^d$ is shown. A detailed estimation is also given in the proof of Proposition \ref{prop_dense_1} below. In particular, the mollification shows that $\Cc^\infty_c(\Omega,Y)$ is dense in $L^p(\Omega,Y)$ and we can extend this density when considering two Banach spaces $Y_2\hookrightarrow Y_1$. 
\begin{prop}\label{prop_dense_3}
Let $Y_1$ and $Y_2$ be two Banach spaces with $Y_2\hookrightarrow Y_1$ densely. For any $p\in [1,+\infty)$, the space $\Cc^\infty_c(\Omega,Y_2)$ of smooth and compactly supported functions with values in $Y_2$ is dense in $L^p(\Omega,Y_1)$.
\end{prop}
\begin{demo}
Let $u\in L^p(\Omega,Y_1)$. By definition of the Bochner integral (see Appendix \ref{section_spaces}), the function $u$ can be approximated in $L^p(\Omega,Y_1)$ by a simple function $v=\sum_{k=1}^n y_k \Un_{x\in A_k}$, where $(A_k)$ is a finite decomposition of a compact subdomain of $\Omega$ in measurable sets and $y_k$ belongs to $Y_1$. By density, we can approximate each $y_k$ by $y'_k$ belonging to $Y_2$. Since there is only a finite number of them, $w=\sum_{k=1}^n y'_k \Un_{x\in A_k}$ can be constructed as close as wanted from $u$ in $L^p(\Omega,Y_1)$. Then, we can use the mollification to obtain a smooth version $\eta_\varepsilon\ast w$ of $w$. To conclude, simply notice that, since $w$ is valued in $Y_2$ which is a vector space, $\eta_\varepsilon\ast w$ is also valued in $Y_2$.
\end{demo}

If the previous result is natural, the case of $p=+\infty$ leads to difficulties. 
Of course, we cannot expect $\Cc^\infty_c(\Omega,Y)$ to be dense in $L^\infty(Y)$ since this is not true even for $Y=\RR$. The following counter-example shows that we can neither expect in general that $L^\infty(\RR^d,Y_2)$ is dense in $L^\infty(\RR^d,Y_1)$ even if $Y_2$ is dense in $Y_1$.
\begin{lemme}\label{lemme_suites_approx}
Let $Y_2\subset Y_1$ be two Banach spaces. Assume that $Y_1$ has infinite dimension and that the embedding $Y_2\hookrightarrow Y_1$ is compact. Then, there exist $\varepsilon>0$ and a sequence $(y_n)$ bounded in $Y_1$ such that any sequence $(z_n)$ of $Y_2$ with $|y_n-z_n|_{Y_1}<\varepsilon$ is unbounded in $Y_2$. 
\end{lemme}
\begin{demo}
Since $Y_1$ has infinite dimension, Riesz Lemma shows that there is a sequence $(y_n)$ in its unit ball and $\varepsilon>0$ such that $|y_p-y_q|_{Y_1}\geq 3\epsilon >0$ for all $p\neq q$. Assume that $(z_n)$ is a sequence of $Y_2$ such that $|y_n-z_n|_{Y_1}<\varepsilon$ and $|z_n|_{Y_2}$ is bounded. Then, by compactness, we can find a subsequence $(z_{\varphi(n)})$ converging in $Y_1$. A fortiori, this subsequence should be a Cauchy sequence in $Y_1$, but this is precluded by
\[|z_p-z_q|_{Y_1} \geq |y_p-y_q|_{Y_1} - |y_p-z_p|_{Y_1} - |z_q-y_q|_{Y_1} > \varepsilon.\]
\end{demo}  
\begin{coro}\label{coro_not_dense}
Let $Y_2\subset Y_1$ be two Banach spaces. Assume that $Y_1$ has infinite dimension and that the embedding $Y_2\hookrightarrow Y_1$ is compact. Then $L^\infty(\Omega,Y_2)$ is \underline{not} dense in $L^\infty(\Omega,Y_1)$ even if $Y_2$ is dense in $Y_1$.
\end{coro}
\begin{demo}
We consider the sequence $(y_n)$ and $\varepsilon>0$ as in Lemma \ref{lemme_suites_approx}. We easily construct $u\in L^\infty(\Omega,Y_1)$ such that, for each $n$, $u$ is constant equal to $y_n$ in some small open ball $B_n\subset \Omega$. Due to Lemma \ref{lemme_suites_approx}, there is no function in $L^\infty(\Omega,Y_2)$ which is $\varepsilon-$close to $u$ in $L^\infty(\Omega,Y_1)$.
\end{demo}

\subsection{The importance of Radon-Nikodym property}\label{section_RN}
As recalled in Appendix \ref{section_spaces}, most of the basic properties of Sobolev spaces of real-valued functions extend to the case of vector-valued functions. But a difficulty appears when we want to characterize the differentiable functions: there, the fact that $Y$ satisfies the Radon-Nikodym property becomes important, see \cite{Arendt-Batty-et-al,Arendt-Kreuter,Diestel-Uhl,Kreuter}. 
\begin{defi}\label{defi_RNP}
A Banach space $Y$ satisfies the {\bf Radon-Nikodym property} if every Lipschitz continuous map $f:\RR\rightarrow Y$ is almost everywhere differentiable. 
\end{defi}
The Radon-Nikodym property is a delicate matter that has been the subject of much study. Fortunately, it will be trivially satisfied in the present paper since we will assume $Y$ to be a Hilbert space. 
\begin{theorem}[Dunford-Pettis]\label{th_RNP}
If $Y$ is separable and if $Y$ is the dual space of a Banach space, then $Y$ satisfies the Radon-Nikodym property. In particular, any reflexive Banach space $Y$ satisfies the Radon-Nikodym property.
\end{theorem}
\begin{demo}
See for example Theorem 1.2.6 and Corollary 1.2.7 of \cite{Arendt-Batty-et-al}. 
\end{demo}

The importance of Radon-Nikodym property in the study of Sobolev spaces of vector-valued functions appears in the following result.
\begin{theorem}\label{th_equiv_W1p}
Let $p\in (1,\infty]$. Assume that $Y$ satisfies the Radon-Nikodym property and assume that $u\in L^p(\RR^d,Y)$. Then $u$ belongs to $W^{1,p}(\RR^d,Y)$ if and only if there exists $C>0$ such that, for all $\tau\in\RR^d$,  
\[ \|u(\cdot - \tau)-u\|_{L^p(\RR^d,Y)} \leq C |\tau|~.\]
If this is the case, then we can take $C=\sqrt{d}\max_{i=1\ldots d} \|\partial_{x_i} u\|_{L^p(\RR^d,Y)}$.
\end{theorem}
\begin{demo}
See Lemma 2.1 and Theorem 2.2 of \cite{Arendt-Kreuter} or Theorem 3.20 of \cite{Kreuter}.
\end{demo}

Note that Theorem 2.5 of \cite{Arendt-Kreuter} and Proposition 2.81 of \cite{HvNvW} show that the Radon-Nikodym property is necessary to obtain this characterization. It fails for example if $Y=L^1(\RR)$ and $u\in L^2(\RR,Y)$ defined by $u(x)(y)=\Un_{[0,1]}(x)\Un_{[0,x]}(y)$ since, obviously, $|\Un_{[0,x+\tau]}-\Un_{[0,x]}|_{Y}\leq|\tau|$ but the derivative of $x\mapsto u(x,\cdot)\in Y$ in the sense of distribution is the Dirac mass, not belonging to $Y$. This implies in particular that $L^1(\RR)$ does not satisfy the Radon-Nikodym property. In the same way, the space $\Cc^0([0,1])$ endowed with the $L^\infty$-norm neither satisfies the Radon-Nikodym property.

\section{The uniformly local Sobolev spaces}\label{section_ul}

Let $Y$ be a Banach space and let $d\in\NN^*$. We introduce the uniformly local norm $L^p_\ul(\RR^d,Y)$ defined by 
\begin{equation}\label{def_norme_unif_loc}
\|u\|_{L^p_\ul(\RR^d,Y)}= \sup_{a\in\RR^d} \left(\int_{B_{\RR^d}(a,1)} |u(x)|_Y^p \d x\right)^{1/p}.
\end{equation}
Two versions of the space associated with this norm can be found: the weak one is simply
\begin{equation}\label{def_espace_unif_loc_weak}
L^p_\ulw(\RR^d,Y):=\{u\in L^p_\loc(\RR^d,Y)~|~\|u\|_{L^p_\ul(\RR^d,Y)}<\infty\}.
\end{equation}
The most usual stronger version of the space includes the continuity with respect to translations: 
\begin{equation}\label{def_espace_unif_loc_strong}
L^p_\uls(\RR^d,Y):=\{u\in L^p_\ulw(\RR^d,Y)~|~\|u-u(\cdot-\xi)\|_{L^p_\ul(\RR^d,Y)}\xrightarrow[~\xi\rightarrow 0~]{} 0\}.
\end{equation}
As an alternative definition, $L^p_\uls(\RR^d,Y)$ is sometimes introduced as the closure of the bounded uniformly continuous functions $BUC(\RR^d,Y)$ in $L^p_\ulw(\RR^d,Y)$.

For any $k\in \NN$, we can define in a similar way the uniformly local Sobolev spaces $W^{k,p}_\ulx(\RR^d,Y)$ (where $*$ is w or s) associated with the norm
\begin{equation}\label{def_norme_Sob_unif_loc}
\|u\|_{W^{k,p}_\ul(\RR^d,Y)}= \sum_{|\nu|\leq k} \|\partial_x^\nu u\|_{L^p_\ul(\RR^d,Y)}.
\end{equation}
Notice that, due to the presence of a supremum in the definition of the norm, the spaces $W^{k,p}_\ulx$ are Banach spaces but are neither reflexive nor separable. 
The continuity with respect to translations is often required, that is that the spaces $W^{k,p}_\uls$ are more frequently used than the spaces $W^{k,p}_\ulw$. Indeed, it is necessary to recover the density of smooth functions, see Propositions \ref{prop_dense_1_bis} and \ref{prop_dense_1} below. However, in the present article, 
using the strong version will not be so helpful because the semigroups associated to our parabolic equation fail to be continuous in general in both cases. 

\medskip

To make the study of uniformly local spaces easier, we notice the following usual trick. For any positive $\mu$ and any $x_\sharp\in\RR^d$, we set  
\begin{equation}\label{eq_rho_mu_x_sharp}
\rho_{\mu,x_\sharp}(x)=\frac 1{\cosh(\mu \sqrt{1+|x-x_\sharp|^2})}.
\end{equation}
We will be able to translate the results obtained in weighted spaces to the uniformly local spaces due to the following fact.
\begin{lemme}\label{lemme_equiv_norms}
Let $p\in[1,+\infty]$. For any $\mu>0$, the uniformly local norm defined by \eqref{def_norme_unif_loc} is equivalent to the norm
\begin{equation}
\|u\|_{L^p_{\mu,\ul}(\RR^d,Y)}:= \sup_{x_\sharp\in\RR^d} \|u\|_{L^p_{\rho_{\mu,x_\sharp}}(\RR^d,Y)}~.
\end{equation}
\end{lemme}
\begin{demo}
Since $\rho_{\mu,0}$ is uniformly positive on the ball $B_{\RR^d}(0,1)$, by translations, the estimate 
\[\|u\|_{L^p_{\ul}(\RR^d,Y)} \leq C \|u\|_{L^p_{\mu,\ul}(\RR^d,Y)}\]
is obvious. The reverse estimate is obtained by writing 
\begin{align*}
\int_{\RR^d} \rho_{\mu,x_\sharp}(x) |u(x)|_Y^p \d x &\leq \sum_{j\in\frac 12 \ZZ^d} \Big( \max_{x\in B_{\RR^d}(j,1)}\rho_{\mu,x_\sharp}(x)\Big) \int_{B_{\RR^d}(j,1)} |u(x)|_Y^p \d x  \\
&\leq \sum_{j\in\frac 12 \ZZ^d} \Big( \max_{x\in B_{\RR^d}(j,1)}\rho_{\mu,x_\sharp}(x)\Big) \|u\|^p_{L^p_{\ul}(\RR,Y)}.
\end{align*}
Then, we simply notice that $\rho_{\mu,x_\sharp}(j)\sim 2e^{-\mu |j-x_\sharp|}$ when $|j|\rightarrow +\infty$, so that  
the sum $\sum_{j\in\frac 12 \ZZ^d} \Big( \max_{x\in B_{\RR^d}(j,1)}\rho_{\mu,x_\sharp}(x)\Big)$ is finite. In addition, we notice that this sum can be bounded uniformly with respect to $x_\sharp$. This concludes the proof by showing that 
\[\|u\|_{L^p_{\mu,\ul}(\RR^d,Y)}^p\leq  \|u\|^p_{L^p_{\ul}(\RR,Y)} \sup_{x_\sharp\in\RR^d} \sum_{j\in\frac 12 \ZZ^d} \Big( \max_{x\in B_{\RR^d}(j,1)}\rho_{\mu,x_\sharp}(x)\Big).\]
\end{demo}

The reason why the continuity of the translations is usually required in the definition of the uniformly local spaces, is to get the density of smooth functions, which does not hold in the weak version.
\begin{prop}\label{prop_dense_1_bis}
If $Y\neq \{0\}$, $k\geq 1$ and $p\in[1,+\infty)$, the space ${W^{k,p}_\ulw(\RR^d,Y)}$ is \underline{not} dense in ${L^p_\ulw(\RR^d,Y)}$. 
\end{prop}
\begin{demo}
First consider the case $d=1$. Pick a non-zero vector $y\in Y$ and consider $u(x)=\sin(x^2)y$. By local boundedness, $u$ belongs to $L^p_\ulw(\RR,Y)$. If we consider $v\in W^{1,p}_\ulw(\RR,Y)$, then, for $|x'-x|\leq 1$, 
\[|v(x')-v(x)|=\left|\int_x^{x'} \partial_x v(\xi)\d\xi \right|\leq |x'-x|^{1-1/p} \|v\|_{W^{1,p}_\ul(\RR,Y)}.\]
Thus, $v$ is uniformly continuous and for large $x$, it cannot be close everywhere to the fast oscillations of $u(x)=\sin(x^2)y$. Thus, $v$ cannot approximate $u$ as precisely as wanted. 

For $p=1$, the above pointwise estimate is replaced by its integrated version:
for $|\tau|\leq 1$,
\[
\int_a^{a+1}|v(x-\tau)-v(x)|_Y\d x\leq \int_0^{|\tau|}\!\!\int_{a-1}^{a+2}|\partial_x v(\xi)|_Y\d\xi\d s
\leq 2|\tau|\,\|v\|_{W^{1,1}_{\ul}(\RR,Y)} ,
\]
so that $v$ belongs to $L^1_{\uls}(\RR,Y)$. Since $L^1_{\uls}(\RR,Y)$ is closed in
$L^1_{\ulw}(\RR,Y)$ and does not contain $u$, $v$ cannot approximate $u$ as precisely as wanted.

The general case of $\RR^d$ is obtained by arguing in the same way with $u(x)=\sin(x_1^2)y$, where $x_1$ is the first component of $x$.
\end{demo}

By contrast, the density holds if we add the requirement of the continuity of translations, bringing some uniform continuity.
\begin{prop}\label{prop_dense_1}
For any $k\geq 0$ and $p\in[1,+\infty)$, ${W^{k,p}_\uls(\RR^d,Y)}\hookrightarrow {L^p_\uls(\RR^d,Y)}$ densely.
\end{prop}
\begin{demo}
The embedding is obvious. To show the density, the simplest idea is to apply the mollification to $u\in L^p_\uls(\RR^d,Y)$. Let $\eta_\varepsilon=\frac{1}{\varepsilon^d}\eta(\cdot/\varepsilon)\in\Cc^\infty(\RR^d,\RR)$ be a smooth approximation of the identity with support in the ball $B_{\RR^d}(0,\varepsilon)$. We argue as in Proposition \ref{prop_dense_4} in appendix. We simply have to notice that the estimates of the norms can be made uniformly with respect to $x\in\RR^d$ because the estimates on $u$ are also uniform. 

To provide the more crucial detail, let us show that $\eta_\varepsilon \ast u$ converges to $u$ in the uniformly local norms. Consider $\rho_{\mu,x_\sharp}$ as in \eqref{eq_rho_mu_x_sharp}, using H\"older inequality, the fact that $\int_{\RR^d} \eta_\varepsilon=1$ and that the support of $\eta_\varepsilon$ is small, we get the bound
\begin{align*}
\| \eta_\varepsilon \ast u - u&\|_{L^p_{\rho_{\mu,x_\sharp}}}^p = \int_{\RR^d} \rho_{\mu,x_\sharp}(x) \left|\int_{\RR^d}\eta_\varepsilon(x')(u(x-x')-u(x)) \d x' \right|^p \d x\\
& \leq  \int_{\RR^d} \rho_{\mu,x_\sharp}(x) \left(\int_{\RR^d} \eta_\varepsilon(x')\d x'\right)^{p-1} \int_{\RR^d}\eta_\varepsilon(x')|u(x-x')-u(x)|^p \d x' \d x\\
& \leq \int _{\RR^d} \eta_\varepsilon(x') \left( \int_{\RR^d}\rho_{\mu,x_\sharp}(x)|u(x-x')-u(x)|^p \d x \right) \d x'\\
&\leq \max_{x'\in B_{\RR^d}(0,\varepsilon)} \|u(\cdot-x')-u\|_{L^p_{\rho_{\mu,x_\sharp}}}^p.
\end{align*}
Due to the continuity with respect to translations assumed in the definition \eqref{def_espace_unif_loc_strong}, $\eta_\varepsilon \ast u$ converges to $u$ in $L^p_{\rho_{\mu,x_\sharp}}(\RR^d,Y)$, uniformly with respect to $x_\sharp$, concluding the proof.
\end{demo}

However, we have to be aware of the following negative result, which holds for both versions of the spaces.
\begin{prop}\label{prop_dense_2}
Let $Y_2\subset Y_1$ be two Banach spaces. Assume that $Y_1$ has infinite dimension and that the embedding $Y_2\hookrightarrow Y_1$ is compact. Then, for any $p\in [1,+\infty]$, ${L^p_\ulx(\RR^d,Y_2)}$ is \underline{not} dense in ${L^p_\ulx(\RR^d,Y_1)}$, where $*$ is either w or s, even if $Y_2$ is dense in $Y_1$.
\end{prop}
\begin{demo}
We use the sequence $(y_n)\subset Y_1$ and $\varepsilon>0$ provided by Lemma \ref{lemme_suites_approx} to exhibit a counterexample of the density. Consider a function $u$ such that, for each $n$, $u$ is constant equal to $y_n$ on a cube $Q_n$ of size $1$ somewhere in $\RR^d$. If the $Q_n$ are far away from each other, it is not difficult to do this in a smooth and uniform way, so that $u$ belongs to ${L^p_\uls(\RR^d,Y_1)}$, and thus also to ${L^p_\ulw(\RR^d,Y_1)}$.
Assume that $v$ belongs to ${L^p_\ulw(\RR^d,Y_2)}$ (or even in the stronger version of the space) and that $v$ is $\eta$-close to $u$ for the norm ${L^p_\ul(\RR^d,Y_1)}$, with $\eta=2^{-(1+1/p)}\varepsilon$. Let us focus on one of the squares $Q_n$. In $Q_n$, the function $v$ can be approximated by a simple function $\sum_{k=1}^K a_k \Un_{x\in A_k}$ in ${L^p(Q_n,Y_2)}$ which is $2\eta$-close to $u$ in ${L^p(Q_n,Y_1)}$ and of norm less than $M=2\|v\|_{L^p_\ul(\RR^d,Y_2)}$ in ${L^p(Q_n,Y_2)}$. We thus have
\[ \sum_{k=1}^K |a_k-y_n|_{Y_1}^p|A_k| \leq (2\eta)^p ~~\text{ and }~~ \sum_{k=1}^K |a_k|_{Y_2}^p|A_k| \leq M^p.\]
Let $\Kc_1=\{k\,|\,|a_k-y_n|_{Y_1}^p > 2 (2\eta)^p\}$ and $\Kc_2=\{k\,|\,|a_k|_{Y_2}^p > 2 M^p\}$. Both sets $\Kc_i$ should be such that $\sum_{k\in\Kc_i} |A_k|<\frac 12$ whereas $\sum_{k=1}^K |A_k|=1$. Thus, there is at least one index $k$ not belonging to $\Kc_1\cup\Kc_2$ and satisfying $|a_k-y_n|_{Y_1} \leq 2^{1+1/p} \eta$ and $|a_k|_{Y_2}\leq 2^{1/p}M$. We denote by $z_n$ this vector $a_k\in Y_2$. Doing this in any $Q_n$, we construct a sequence $(z_n)$, which is bounded by $2^{1/p}M$ and satisfies $|z_n-y_n|_{Y_1} \leq 2^{1+1/p} \eta=\varepsilon$. This yields a contradiction with Lemma \ref{lemme_suites_approx}, meaning that $v$ cannot be taken as close as wanted to $u$. 
\end{demo}


\section{Abstract parabolic operators with non-dense domain}\label{section_nondense}
Nowadays, the theory of sectorial operators and their associated analytic semigroup is well known. In the present paper, we will consider operators satisfying sectorial-like resolvent estimates but with non-dense domain. A more involved semigroup theory has been constructed in this case. The purpose of this section is to recall part of this theory, in particular some of the results of the fundamental works of Sinestrari, Da Prato and Lunardi, see \cite{da-Prato-Sinestrari,Lunardi,Sinestrari}. 

More precisely, we consider the following type of operators.
\begin{defi}\label{defi_op_nondense}
Let $X$ be a Banach space and let $A:D(A)\subset X\rightarrow X$ a closed linear operator. We say that $A$ is an {\bf abstract parabolic operator} if there exist $\omega\in\RR$ and $\phi\in (0,\pi/2)$ such that the sector 
$$S_{\omega,\phi}=\{ z \in \CC\setminus\{\omega\}~,~~ \phi \leq |\arg(z-\omega)| \leq \pi\}$$
is in the resolvent set of $A$ and there exists $M>0$ such that 
$$\forall z\in S_{\omega,\phi}~,~~\|(z\id-A)^{-1}\|_{\Lc(X)} \leq \frac {M}{|z-\omega|}~. $$
\end{defi}
Notice that, compared to the sectorial operators, the only missing property is the density of the domain, see the brief recalls in Appendix \ref{section_sectoriel} or  \cite{Henry}. We also underline that the name ``abstract parabolic operator'' is not conventional, since we do not find a common name in the literature. But the associated semigroup is called ``parabolic abstract evolution'' in \cite{da-Prato-Sinestrari} and the associated equation is called ``abstract parabolic equation'' in \cite{Lunardi}. In the whole section, $A$ and $X$ are assumed to be as above.

\subsection{The integrated semigroup}
Exactly as in the classical sectorial case, we can introduce an associated semigroup by setting 
\begin{equation}\label{defi_semigroup}
e^{-At}=\frac 1{2i\pi}\int_\Gamma (\lambda+A)^{-1}e^{\lambda t}\d \lambda
\end{equation}
where $\Gamma$ is a suitable contour, see \cite{Sinestrari}. Notice that we choose the convention opposite to the historical papers where $\tilde A=-A$ and the semigroup is therefore $e^{\tilde At}$. This semigroup is called the {\bf integrated semigroup} generated by $A$. It has all the expected properties, similar to the ones of the classical analytic semigroups generated by sectorial operators, except the continuity at $t=0$. Let us recall some of them.
\begin{prop}\label{prop_integrated_semigroup}
Let $A$ be an abstract parabolic operator of parameters $\omega\in\RR$ and $\phi\in (0,\pi/2)$ as in Definition \ref{defi_op_nondense} and let $e^{-At}$ be its integrated semigroup. Then:
\begin{enumerate}
\item[(i)] There exists $M>0$ such that $\|e^{-At}\|_{\Lc(X)}\leq M e^{-\omega t}$ for all $t\geq 0$. 
\item[(ii)] For all $t>0$ and $x\in X$, $e^{-At}x$ belongs to $D(A)$ and there exists $M$ such that 
$$\forall t>0~,~~\|A e^{-At}\|_{\Lc(X)}\leq \frac {Me^{-\omega t}}{t}.$$
\item[(iii)] For all $t>0$ and $x\in D(A)$, $Ae^{-At}x=e^{-At}Ax$.
\item[(iv)] The function $t\mapsto e^{-At}$ can be extended analytically in a sector containing the positive real axis and, for all $t>0$, $\frac{\textrm{d}}{\textrm{d}t} e^{-At}=-A e^{-At}$.
\end{enumerate}
\end{prop}
\begin{demo}
We refer to \cite{da-Prato-Sinestrari,Sinestrari} for most of the proof and details. The only point to notice is that these articles only consider the case $\omega=0$. Thus, we simply have to go back to $\tilde A=A-\omega \id$, to define the integrated semigroup of $A$ by $e^{-At}=e^{-\omega t}e^{-\tilde A t}$ and to apply the results of \cite{da-Prato-Sinestrari,Sinestrari} to $\tilde A$. 
\end{demo}

The main consequence of the possible non-density of the domain is the behaviour when $t\rightarrow 0^+$.
\begin{prop}
If $x\in \overline{D(A)}$, then $e^{-At}x\rightarrow x$ when $t\rightarrow 0^+$. Conversely, if the limit $\lim_{t\rightarrow 0^+} e^{-At}x$ exists, then this limit is $x$ and $x$ belongs to $\overline{D(A)}$.  
\end{prop}
As a consequence, if $D(A)$ is not dense, there exists initial data $x$ such that $e^{-At}x$ does not converge to $x$, or anything else, when $t$ goes to $0$. This difficulty leads to introduce a weak notion of solutions of equations as $\partial_t x(t)=Ax(t)$ with initial data $x_0$: the ``integral solutions'', which satisfy $x(t)=x_0+A\int_0^t x(s)\d s$, see \cite{da-Prato-Sinestrari,Magal-Ruan}. 
However, we will not use this concept in the present paper, because we will consider regular initial data, which is nevertheless natural for nonlinear parabolic problems. 
We simply enhance the regularity of the integral solution, see Proposition 1.2 of \cite{Sinestrari}
\begin{prop}\label{prop_regu_integrated}
For all $x\in X$ and $t>0$, $\int_0^t e^{-As}x\d s$ belongs to $D(A)$ and $$A\int_0^t e^{-As}x\d s=e^{-At}x-x.$$
\end{prop}

\subsection{The associated sectorial operator}
We briefly notice that, to come back to the classical parabolic theory of operator with dense domain, we can restrict $A$ to smaller spaces to make its domain dense. 
\begin{prop}[Propositions 1.2 and 1.3 of \cite{Sinestrari}]\label{prop_A0}$~$\\
Let $X_0:=\overline{D(A)}$ be the closure of $D(A)$ in $X$ and let $A_0$ be the restriction of $A$ to
$$D(A_0)~:=~\{\,x\in D(A)\,,~ Ax \in X_0\,\}.$$
Then $A_0$ is a sectorial operator in $X_0$ and, in particular, has a dense domain. Moreover, $e^{-A_0 t}$ and $e^{-A t}$ coincide for all $x\in X_0$.
\end{prop}
We should be able to use this classical analytic semigroup to set a Cauchy problem with regular initial data. However, in order to apply the classical theory to $A_0$, we have to ensure that a given nonlinear term is valued in $X_0$. A concrete description of this space $X_0$ is not so simple in our abstract context. 

\subsection{The intermediate spaces}\label{section_intermediate}
In the classical parabolic theory, the fractional spaces $D(A^\alpha)$ are very useful. However, it is not easy to define the fractional powers of abstract parabolic operators with non-dense domain, see \cite{MS} for example. Moreover, it could be difficult to describe them, as well as the space $D(A_0^\alpha)$ of the sectorial operator $A_0$ of Proposition \ref{prop_A0}. 
The key idea is to consider the intermediate spaces $D_A(\theta,\infty)$ introduced by Sinestrari, see \cite{da-Prato-Sinestrari,Sinestrari}. 
\begin{defi}\label{defi_inter_spaces}
For any $\theta\in(0,1)$, we set 
$$D_{A}(\theta,\infty):=\Big\{x\in X~,~~\sup_{t\in(0,1]} \frac {\|e^{-A t}x-x\|_{X}}{t^\theta}<\infty\Big\}$$
endowed with the norm
$$\|x\|_{D_A(\theta,\infty)}=\|x\|_X+\sup_{t\in(0,1]} \frac {\|e^{-A t}x-x\|_{X}}{t^\theta}.$$
\end{defi}
We recall that they are Banach spaces. Moreover, these spaces exactly correspond to the 
real interpolation method of Lions and Peetre in \cite{Lions,Lions-Peetre}, that is 
$$D_A(\theta,\infty)~=~\big(X,D(A)\big)_{\theta,\infty}$$
with the notations of Appendix \ref{section_interpolation}, see \cite{da-Prato-Grisvard}. Using the interpolation and the estimates of Proposition \ref{prop_integrated_semigroup}, it is not difficult to obtain the following estimations. Notice that the last ones correspond to a regularization effect of the parabolic semigroup.
\begin{prop}\label{prop_estimation_semigroup_interm}
Let $A$ be an abstract parabolic operator of parameters $\omega\in\RR$ and $\phi\in (0,\pi/2)$ as in Definition \ref{defi_op_nondense}. For all $\theta\in(0,1)$, there exists $M>0$ such that the following estimates hold:
\begin{itemize}
\item The space $D_A(\theta,\infty)$ is invariant by $e^{-At}$ and, for all $x\in D_A(\theta,\infty)$,
$$\forall t>0~,~~\|e^{-At}x\|_{D_A(\theta,\infty)}\leq Me^{-\omega t}\|x\|_{D_A(\theta,\infty)}.$$
\item For all $x\in D_A(\theta,\infty)$, 
$$\forall t\in(0,1]~,~~\|e^{-At}x-x\|_X\leq Mt^\theta \|x\|_{D_A(\theta,\infty)}.$$
\item For all $x\in X$,     
$$\forall t>0~,~~\|e^{-At}x\|_{D_A(\theta,\infty)}\leq M\left(1+\frac 1{t^\theta}\right)e^{-\omega t} \|x\|_X.$$
\item For all $\alpha \in (0,\theta)$ and $x\in D_A(\alpha,\infty)$,    
$$\forall t>0~,~~\|e^{-At}x\|_{D_A(\theta,\infty)}\leq M\left(1+\frac 1{t^{\theta-\alpha}}\right)e^{-\omega t} \|x\|_{D_A(\alpha,\infty)}.$$
\item For all $x\in D_A(\theta,\infty)$, $e^{-At}x$ belongs to $D(A)$ and 
$$\forall t>0~,~~\|Ae^{-At}x\|_X\leq Mt^{\theta-1}e^{-\omega t}\|x\|_{D_A(\theta,\infty)}.$$
\end{itemize}
\end{prop}
\begin{demo}
See (1.8), (1.14), (1.23), (1.36) and (1.37) of \cite{Sinestrari}.
\end{demo}

\subsection{An abstract Cauchy problem}\label{section_lunardi}
In the present article, we aim at setting a Cauchy problem for a nonlinear parabolic equation. The first historical papers of da Prato and Sinestrari only consider linear problems with forcing terms. Da Prato and Grisvard consider nonlinear equations in \cite{da-Prato-Grisvard}, but with strong assumptions. In \cite{Lunardi}, Lunardi considers a nonlinear abstract parabolic problem with assumptions close to the ones of the classical parabolic problems. Actually, her setting is more general than what we need here since it involves time-dependent quasilinear operators. We will only state here a simplified version of her results, we refer to \cite{Lunardi} for the general statements.

Let $A$ be an abstract parabolic operator on $X$. Let $\theta\in(0,1)$ and let $f:D_A(\theta,\infty)\rightarrow X$ a function which is locally Lipschitz continuous in the following sense: for all $x_*\in D_A(\theta,\infty)$, there exist $r>0$ and $K>0$ such that
$$\forall x,x'\in B_\theta(x_*,r)~,~~\|f(x)-f(x')\|_X \leq K \|x-x'\|_{D_A(\theta,\infty)}$$
where $B_\theta(x_*,r)$ is the ball of center $x_*$ and radius $r$ in the intermediate space $D_A(\theta,\infty)$. We consider the abstract parabolic evolution equation
\begin{equation}\label{eq_abstract}
\frac{\textrm{d}}{\textrm{d}t} x(t)=-Ax(t) + f(x(t)). 
\end{equation}
\begin{theorem}[simplified version of Theorem 2.1 of \cite{Lunardi}]\label{th_lunardi}
Let $\theta\in(0,1)$, let $A$ be an abstract parabolic operator and let $f$ be a locally Lipschitz function from $D_A(\theta,\infty)$ to $X$ as above. For all $x_*\in D_A(\theta,\infty)$, there exists a ball $B_\theta(x_*,r)$ in $D_A(\theta,\infty)$ and $T>0$ such that, for all $x_0\in B_\theta(x_*,r)$, there exists a unique solution of \eqref{eq_abstract} with initial data $x(0)=x_0$ and time interval $[0,T]$. Moreover, this solution is classical in the following sense:
\begin{itemize}
\item the function $t\mapsto x(t)$ is continuous from $[0,T]$ into $D_A(\theta,\infty)$ and we have $x(0)=x_0$,
\item the function $t\mapsto x(t)$ is continuous $(0,T]$ to $D(A)$ and of class $\Cc^1$ from $(0,T]$ to $X$,
\item for all $t>0$, the equation \eqref{eq_abstract} is satisfied in the strong sense.
\end{itemize}
\end{theorem}
Notice that, in \cite{Lunardi}, there is a small omission since the notion of classical solution does not include any continuity at $t=0$, so that the condition $x(0)=x_0$ has no meaning. But the proof of the result involves a classical fixed point theorem and the characterization of the solution by the Duhamel's formula 
\begin{equation}\label{eq_duhamel}
x(t)=e^{-tA}x_0+\int_0^t e^{-A(t-s)}f(x(s))\d s. 
\end{equation}
Due to estimates as the ones of Proposition \ref{prop_estimation_semigroup_interm}, the continuity at $t=0$ for $x_0\in D_A(\theta,\infty)$ holds. Notice the importance of having $\theta>0$, contrarily to the classical sectorial case, otherwise, even for $f\equiv 0$, the convergence to the initial data may not hold. 

In the present paper, we would like to consider the problem of global existence. Once all the above theory is constructed, the arguments are simply the same as the classical ones, as in \cite{Henry} for example. For sake of completeness, we provide a short proof.
\begin{prop}\label{prop_lunardi_glob}
Let $\theta\in(0,1)$, let $A$ be an abstract parabolic operator and let $f$ be a globally Lipschitz function from $D_A(\theta,\infty)$ to $X$. For all $x_0\in D_A(\theta,\infty)$, there exists a unique global solution of \eqref{eq_abstract} with initial data $x(0)=x_0$ and time interval $[0,+\infty)$. Moreover, this solution is classical in the sense of Theorem \ref{th_lunardi} and $\|x(t)\|_{D_A(\theta,\infty)}$ growth at most exponentially fast.
\end{prop}
\begin{demo}
Let $\lambda>0$ to be fixed later and let $\omega\in\RR$ and $\phi\in (0,\pi/2)$ be as in Definition \ref{defi_op_nondense}. Let $K$ be the Lipschitz constant of $f$. For any continuous function $x(t)$ from $\RR_+$ to $D_A(\theta,\infty)$, we define the norm 
$$\|x\|_\lambda=\sup_{t\geq 0} e^{-\lambda t}\|x(t)\|_{D_A(\theta,\infty)}$$
and the Banach space $\Cc_\lambda$ of the functions of $\Cc^0(\RR_+,D_A(\theta,\infty))$ with finite norm. 

We introduce the function $\psi: \Cc_\lambda\rightarrow \Cc_\lambda$ defined by $y=\psi(x)$ with 
$$y(t)= e^{-At}x_0+\int_0^t e^{-A(t-s)}f(x(s))\d s.$$
We aim to show that $\psi$ is a well-defined contraction on $\Cc_\lambda$. This is done by using Proposition \ref{prop_estimation_semigroup_interm}. Let us simply detail the estimation of the Lipschitz constant. Let $x$ and $\tilde x$ in $\Cc_\lambda$ and let $y=\psi(x)$ and $\tilde y=\psi(\tilde x)$. Let $t\geq 0$, we have
\begin{align*}
e^{-\lambda t}\|y(t)-\tilde y(t)&\|_{D_A(\theta,\infty)} = e^{-\lambda t}\left\|\int_0^t e^{-A(t-s)}\big(f(x(s))-f(\tilde x(s))\big) \d s\right\|_{D_A(\theta,\infty)}\\
& \leq  \int_0^t e^{-\lambda t} M e^{-\omega(t-s)}\Big(1+\frac 1{(t-s)^\theta}\Big)\big\|f(x(s))-f(\tilde x(s))\big\|_X \d s \\
& \leq  \int_0^t e^{-(\lambda+\omega) (t-s)} M \Big(1+\frac 1{(t-s)^\theta}\Big)e^{-\lambda s} K \|x(s)-\tilde x(s)\|_{D_A(\theta,\infty)} \d s \\
& \leq  MK \|x-\tilde x\|_\lambda \int_0^t e^{-(\lambda+\omega) \tau }  \Big(1+\frac 1{\tau^\theta}\Big) \d \tau  \\
&\leq  MK \|x-\tilde x\|_\lambda \int_0^{+\infty} e^{-(\lambda+\omega) \tau } \Big(1+\frac 1{\tau^\theta}\Big) \d \tau 
\end{align*}
It remains to notice that, if $\lambda>-\omega$, $\int_0^{+\infty} e^{-(\lambda+\omega) \tau }\d \tau=\frac 1{\lambda+\omega}$ and, remembering that $\theta\in (0,1)$, that 
$$\int_0^{+\infty} \frac 1{\tau^\theta} e^{-({\lambda+\omega}) \tau }\d \tau  = \frac 1{({\lambda+\omega})^{1-\theta}} \int_0^{+\infty} \frac {1}{\sigma^\theta} e^{-\sigma }\d \sigma=\frac C{({\lambda+\omega})^{1-\theta}}.$$
Finally, we obtain that $\|y-\tilde y\|_\lambda\leq \frac 12 \|x-\tilde x\|_\lambda$ for $\lambda>0$ large enough. The fact that the range of $\psi$ indeed belongs to $\Cc_\lambda$ is shown by similar estimates since the term $e^{-At}x_0$ belongs to $D_A(\theta,\infty)$ with an at most exponential growth as shown by Proposition \ref{prop_estimation_semigroup_interm}. 

The unique fixed point of $\psi$ provides a mild solution by definition. 
To check that this mild solution is classical, this can be done in the same way as \cite{Lunardi} for example, see also \cite{da-Prato-Grisvard,da-Prato-Sinestrari}. 
Also notice that the uniqueness follows from Theorem \ref{th_lunardi}: the interval of times where two solutions coincide is closed by continuity and is open due to the local Cauchy problem stated in Theorem \ref{th_lunardi}. A posteriori, this shows that the only possible mild solution belongs to $\Cc_\lambda$ and has a norm growing at most as fast as $e^{\lambda t}$.   
\end{demo}

\begin{coro}\label{coro_lunardi_glob}
Let $\theta\in(0,1)$, let $A$ be an abstract parabolic operator and let $f$ be a function from $D_A(\theta,\infty)$ to $X$, which is Lipschitz continuous in any bounded set of $D_A(\theta,\infty)$. Let $x_0\in D_A(\theta,\infty)$, let $x(t)$ be the unique solution of \eqref{eq_abstract} with initial data $x(0)=x_0$ given by Theorem \ref{th_lunardi} and let $[0,T)$ be the maximal time interval of existence. If $T<+\infty$, then $x(t)$ blows up in the sense that
$$\|x(t)\|_{D_A(\theta,\infty)}~\xrightarrow[~t\longrightarrow T^-~]{}+\infty.$$
\end{coro}
\begin{demo}
If $T<+\infty$ and $x(t)$ does not blow up, it must remain in a bounded ball, since estimates similar to the ones obtained in the previous proof show that the variations of the norm are locally bounded. Thus, we can truncate $f$ in the ball of double radius  to get a globally Lipschitz function and apply Proposition \ref{prop_lunardi_glob}. This global solution coincides with the original one and even extends it with the same nonlinearity at least for times $t\geq T$ sufficiently close to $T$, leading to a contradiction.    
\end{demo}


\section{Auxiliary linear operators}\label{section_A}
To be able to define our linear differential operator in the uniformly local Sobolev spaces, we first consider the same differential operator in different functional spaces.

\subsection{The self-adjoint case}\label{section_As}

We consider an unbounded linear operator of the type 
$$-\div \big( a(x)\grad \cdot \big) + B~~\text{ with }~~D(B)\subset Y$$ 
in $L^2(\RR,Y)$. The (possibly unbounded) operator $B$ is positive, self-adjoint and defined on a dense domain $D(B)$. We assume that the matrix of diffusion coefficients $a$ is of class $\Cc^1(\RR^d,\Lc(Y^d))$, where $\Lc(Y^d)$ are the bounded linear operators from $Y^d$ into itself, with a uniform bound $M_a$:
\begin{equation}\label{hyp_a1}
\forall x\in \RR^d~, \|a(x)\|_{\Lc(Y^d)}+\sum_{i=1}^d \|\partial_{x_i} a(x)\|_{\Lc(Y^d)}\leq M_a.
\end{equation}
Moreover, we assume that $a(x)$ is symmetric, that is that $\langle a(x)\xi|\zeta\rangle_{Y^d}=\langle \xi|a(x)\zeta\rangle_{Y^d}$ for all $x\in\RR^d$ and $\xi,\zeta\in Y^d$, and that there exists $m_a>0$ such that
\begin{equation}\label{hyp_a2}
\forall x\in \RR^d~,~~\forall \xi \in Y^d~,~~ \langle a(x)\xi |\xi\rangle \geq m_a |\xi|^2.
\end{equation}

\begin{lemme}\label{lemme_quadratic}
Consider the above setting. Let $Q=H^1(\RR^d,Y)\cap L^2(\RR^d,D(B^{1/2}))$. The bilinear form 
$$q:\left(\begin{array}{ccc}  Q \times Q & \longrightarrow & \RR \\ (u,v) & \longmapsto & \int_{\RR^d} \langle a(x)\grad u | \grad v\rangle_{Y^d} + \int_{\RR^d}  \langle B^{1/2} u |  B^{1/2} v\rangle_{Y} \end{array}\right)
$$
is closed symmetric positive with dense domain in $L^2(\RR^d,Y)$. Thus, it can be uniquely associated with a self-adjoint operator $A_\s$ with domain $D(A_\s)\subset Q$ dense in $L^2(\RR^d,Y)$, satisfying 
\begin{equation}\label{eq_lemme_quadratic}
\forall u\in D(A_\s)~,~~\forall v\in Q~,~~ \langle A_\s u | v \rangle_{L^2(\RR^d,Y)} = q(u,v)~.
\end{equation}
\end{lemme}
\begin{demo}
We follow here Section VIII.6 of \cite{Reed-Simon}. The symmetry and the positivity are obvious. The closedness follows from the fact that the bilinear form is equivalent to the natural scalar product of $Q$ due to Assumptions \eqref{hyp_a1} and \eqref{hyp_a2}. 
The density of $Q$ in $L^2(\RR^d,Y)$ can be obtained by noticing that $H^1(\RR^d,D(B^{1/2}))$ is dense in $L^2(\RR^d,Y)$, see Proposition \ref{prop_dense_3}. 
To conclude, we use Theorem VIII.15 of \cite{Reed-Simon}: since $q$ is a closed positive bilinear form, it defines a unique self-adjoint operator $A_\s$ in $L^2(\RR^d,Y)$. The proof of Theorem VIII.15 of \cite{Reed-Simon} also shows that $D(A_\s)$ is dense in $Q$.
\end{demo}

We naturally expect that the operator $A_\s$ defined by the previous lemma is $-\div\big( a(x)\grad \cdot \big) + B$ in some sense. The main difficulty is to obtain its domain $D(A_\s)$ to ensure that both $B$ and $-\div\big(a(x)\grad \cdot \big)$ are well defined on it, meaning that $-\div\big( a(x)\grad \cdot \big) + B$ is a well defined sum and not only a formal notation. As explained in the introduction of this paper, this type of problem is less obvious as it may look. The following result can be seen as a consequence of \cite[Theorem 3.14]{da-Prato-Grisvard}. However, we choose to present an alternative stand-alone proof, which follows the classical method of construction of the Laplacian operators. 
\begin{prop}\label{prop_operateur_plat}
The linear operator $A_\s$ defined in Lemma \ref{lemme_quadratic} is the operator 
\begin{equation}\label{def-As}
A_\s=-\div\big( a(x)\grad \cdot \big) + B
\end{equation}
defined from $D(A_\s)=H^2(\RR^d,Y)\cap L^2(\RR^d,D(B))$ into $L^2(\RR^d,Y)$ and the domain norm $\|\cdot\|_{D(A_\s)}:=\|A_\s\cdot\|_{L^2(\RR^d,Y)}$ is equivalent to the norm $\|\cdot\|_{H^2(\RR^d,Y)\cap L^2(\RR^d,D(B))}$. Moreover, the operator $A_\s$ is a closed positive self-adjoint operator of dense domain in $L^2(\RR^d,Y)$. 
\end{prop}
\begin{demo}
First notice that the claimed domain is included in $Q$ and that, on this domain, the claimed definition \eqref{def-As} and the claimed domain indeed satisfy \eqref{eq_lemme_quadratic}. To check this, we recall that the integration by parts is possible in this abstract setting, see Proposition \ref{prop_inte_parts} in appendix. It is also obvious that  $\|\cdot\|_{D(A_\s)}$ is dominated by the norm $\|\cdot\|_{H^2(\RR^d,Y)\cap L^2(\RR^d,D(B))}$ for all $u\in H^2(\RR^d,Y)\cap L^2(\RR^d,D(B))$.

To identify the domains and to show that $-\div\big( a(x)\grad \cdot \big)  + B$ is the correct expression of $A_\s$, it remains to show that any $u$ belonging to the domain of the operator $A_\s$ defined in Lemma \ref{lemme_quadratic} belongs to $H^2(\RR^d,Y)\cap L^2(\RR^d,D(B))$. We follow the classical arguments for proving the regularity of the Laplacian operator, see \cite{Brezis} for example.
Let $u\in D(A_\s)$. For all $\tau\in\RR^d$, we introduce the difference operator defined in $L^2(\RR^d,Y)$ by
$$D_\tau u = \big(u(\cdot-\tau)-u\big)~.$$
We notice that $(D_\tau)^*=D_{-\tau}$ and that $D_\tau$ commutes with the spatial derivatives $\partial_{x_i}$ as well as with $B$. In particular, $v:=D_{-\tau}(D_\tau u)$ belongs to the domain $Q$ of the quadratic form. By definition of  $A_\s$, we have 
\begin{align}
\int_{\RR^d} \langle A_\s u | v \rangle &= \int_{\RR^d} \langle a(x) \grad u | \grad v\rangle_{Y^d} + \int_{\RR^d}  \langle B^{1/2} u |  B^{1/2} v\rangle_Y \nonumber\\
&= \int_{\RR^d}  \langle a(x) \grad u | \grad D_{-\tau}(D_\tau u)\rangle_{Y^d} + \int_{\RR^d}  \langle B^{1/2} u |  B^{1/2} D_{-\tau}(D_\tau u)\rangle_Y \nonumber\\
&= \int_{\RR^d} \langle D_\tau (a(x)\grad u) | D_\tau \grad u\rangle_{Y^d} + \int_{\RR^d}  \langle D_\tau B^{1/2}  u | D_\tau B^{1/2} u\rangle_Y \nonumber\\
&= \int_{\RR^d}  \langle a(x-\tau)D_\tau \grad  u | D_\tau \grad u\rangle_{Y^d} + \int_{\RR^d} \langle (D_\tau a)(x)\grad  u | D_\tau \grad u\rangle_{Y^d}\nonumber\\ &~~~~ + \int_{\RR^d}  \langle D_\tau B^{1/2}  u | D_\tau B^{1/2} u\rangle_Y. \label{eq_demo_operateur_plat}
\end{align}
Due to the bound \eqref{hyp_a2}, we know that 
$$\forall \tau\in\RR^d~,~~  \|D_\tau (\grad u)\|_{L^2(\RR^d,Y^d)}^2  ~ \leq~  C \int_{\RR^d}  \langle a(x-\tau)D_\tau \grad  u | D_\tau \grad u\rangle_{Y^d}.$$
We bound the right-hand side of the previous inequality by using \eqref{eq_demo_operateur_plat}, Cauchy-Schwarz inequality (see Proposition \ref{prop_Holder} in appendix) and the bound \eqref{hyp_a1}. We obtain that, for all $\tau\in\RR^d$,
\begin{align}
\|D_\tau (\grad u)\|_{L^2(\RR^d,Y^d)}^2  ~
\leq~ & \int_{\RR^d} \langle A_\s u | v \rangle -  \int_{\RR^d} \langle (D_\tau a)(x)\grad  u | D_\tau \grad u\rangle_{Y^d}\nonumber \\
& -\int_{\RR^d}  \langle D_\tau B^{1/2}  u | D_\tau B^{1/2} u\rangle_Y \nonumber \\
\leq~ &  C \|A_\s u\|_{L^2(\RR^d,Y)} \|v\|_{L^2(\RR^d,Y)} \nonumber \\
& + C |\tau|\cdot \|\grad u\|_{L^2(\RR^d,Y^d)} \|D_\tau \grad u\|_{L^2(\RR^d,Y^d)} ~. \label{eq_demo_operateur_plat_1}
\end{align}
Due to Theorem \ref{th_equiv_W1p}, we have 
\begin{equation}\label{eq_demo_operateur_plat_2}
\|v\|_{L^2(\RR^d,Y)} = \|D_{-\tau}D_\tau u\|_{L^2(\RR^d,Y)} \leq |\tau| \cdot \|\grad D_\tau u\|_{L^2(\RR^d,Y^d)}~.
\end{equation}
The above estimates \eqref{eq_demo_operateur_plat_1} and \eqref{eq_demo_operateur_plat_2} show that 
\[\forall \tau\in\RR^d~,~~  \|D_\tau \grad u\|_{L^2(\RR^d,Y^d)} \leq  C |\tau| \left( \|A_\s u\|_{L^2(\RR^d,Y)} + \|u\|_{H^1(\RR^d,Y)}\right)~.\]
By assumption, $u$ belongs to $D(A_\s)\subset H^1(\RR^d,Y)$, so the norms of the right side are bounded. Applying again Theorem \ref{th_equiv_W1p}, we conclude that each spatial derivative $\partial_{x_i} u$ belongs to $H^1(\RR^d,Y)$ and thus that $u$ belongs to $H^2(\RR^d,Y)$. We also get the estimate 
\[\forall u\in D(A_\s)~,~~ \|u\|_{H^2(\RR^d,Y)} \leq C\|A_\s u\|_{L^2(\RR^d,Y)}.\]
It remains to show that $u$ belongs to $L^2(\RR^d,D(B))$. To this end, we go back to the definition of $q$. For any $\varphi \in Q=D(q)$, we have that 
\[\int_{\RR^d} \langle A_\s u | \varphi \rangle = \int_{\RR^d}  \langle a(x)\grad  u | \grad \varphi\rangle_{Y^d} + \int_{\RR^d}  \langle B^{1/2} u |  B^{1/2} \varphi\rangle_Y~.\]
We know that $u\in H^2(\RR^d,Y)$ and, thus, we can integrate by parts: 
\[\int_{\RR^d} \langle a(x) \grad u | \grad \varphi\rangle_{Y^d}= - \int_{\RR^d}  \langle \div a(x) \grad u | \varphi\rangle_Y \]
This yields 
\begin{align*}
\left| \int_{\RR^d}  \langle B u |  \varphi \rangle_{D(B^{-1/2}),D(B^{1/2})} \right|&=\left| \int_{\RR^d}  \langle B^{1/2} u |  B^{1/2} \varphi \rangle_Y \right|\\& \leq \big(\|A_\s u \|_{L^2(\RR^d,Y)} + 2C\|u\|_{H^2(\RR^d,Y)}\big) \|\varphi\|_{L^2(\RR^d,Y)}
\end{align*}
for all $\varphi \in Q$. Since $Q$ is dense in $L^2(\RR^d,Y)$, this shows that $Bu$ belongs to $L^2(\RR^d,Y)$ and provides the bound $\|u\|_{L^2(\RR^d,D(B))} \leq C \|A_\s u\|_{L^2(\RR^d,Y)}$ for all $u\in D(A_\s)$. This concludes the proof.
\end{demo}

Since $A_\s$ and $B$ are self-adjoint positive operators with dense domains, $-A_\s$ generates a classical analytic semigroup $e^{-A_\s t}$ on $L^2(\RR^d,Y)$ and we can define the fractional powers $A_\s^\alpha$ and $B^\alpha$, see \cite{Henry} for example. We finish this subsection by a study of the spaces $D(A_\s^\alpha)$. 
\begin{prop}\label{prop_domain_A}
The domain of $A_\s$ satisfies
\[D(A_\s)\hookrightarrow H^1(\RR^d,D(B^{1/2})).\]
Moreover, 
\begin{itemize}
\item if $d=1$, then $D(A_\s)\hookrightarrow L^\infty(\RR,D(B^{3/4}))$.
\item if $d=2$, then $D(A_\s)\hookrightarrow L^\infty(\RR^2,D(B^\alpha))$ for all $\alpha\in [0,1/2)$.
\item if $d=3$, then $D(A_\s)\hookrightarrow L^\infty(\RR^3,D(B^\alpha))$ for all $\alpha\in [0,1/4)$.
\end{itemize}
\end{prop}
\begin{demo}
To obtain the first embedding, we simply notice that, if $u\in D(A_\s)=H^2(\RR^d,Y)\cap L^2(\RR^d,D(B))$, then
\[\int_{\RR^d}|\partial_{x_i} B^{1/2}u|_Y^2 = - \int_{\RR^d} \langle \partial_{x_ix_i}^2 u | B u \rangle_Y \leq \|u\|_{H^2(\RR^d,Y)}\|u\|_{ L^2(\RR^d,D(B))}.\]
The next embeddings depend on the spatial dimension. We know that any $u\in D(A_\s)$ belongs to $H^1(\RR^d,D(B^{1/2}))\cap L^2(\RR^d,D(B))$. For $d=1$, we can directly apply Propositions \ref{prop_interpolation_Amann} and \ref{prop_interpolation_autoadjoint} in appendix to get
\[H^1(\RR,D(B^{1/2}))\cap L^2(\RR,D(B)) \hookrightarrow L^\infty(\RR, D(B^{3/4})).\]
Notice that we use here the ``reiteration principle'' to deduce from the interpolation between $Y$ and $D(B)$ the interpolation between $D(B^\alpha)$ and $D(B^\beta)$, see Sections I.2.8 and I.2.9 of \cite{Amann-livre}.  

Assume now that $d=2$ and that $u\in D(A_\s)$. Fix $\alpha\in [0,1/2)$ and $q>2$. Since $u$ belongs to $H^1(\RR^2,D(B^{1/2}))$,  the embeddings of Theorem \ref{th_embeddings} in appendix show that $u$ belongs to $L^q(\RR^2,D(B^{1/2}))\subset L^q(\RR^2,D(B^\alpha))$. Since $\partial_{x_i}u$ belongs to $H^1(\RR^2,Y)$, $\partial_{x_i}u$ belongs to $L^{q'}(\RR^2,Y)$ for all $q'\geq 2$, as large as needed. On the other hand, $\partial_{x_i}u$ belongs to $L^2(\RR^2,D(B^{1/2}))$. So, choosing $q'$ large enough and applying Propositions \ref{prop_interpolations}, \ref{prop_interpolation_Lp} and \ref{prop_interpolation_autoadjoint} in appendix, we obtain that $\partial_{x_i}u$ belongs to $L^q(\RR^2,D(B^{\alpha}))$. At the end, we have shown that $u$ belongs to $W^{1,q}(\RR^2,D(B^\alpha))$ which is embedded in $L^\infty(\RR^d,D(B^\alpha))$ since $q>2$.

Finally, for $d=3$, the arguments are similar. Since $u$ belongs to $H^1(\RR^3,D(B^{1/2}))$, the embeddings of Theorem \ref{th_embeddings} show that $u$ belongs to $L^6(\RR^3,D(B^{1/2}))$. Since $\partial_{x_i}u$ belongs to $H^1(\RR^3,Y)$, $\partial_{x_i}u$ belongs to $L^6(\RR^3,Y)$. On the other hand, $\partial_{x_i}u$ belongs to $L^2(\RR^3,D(B^{1/2}))$. Applying the interpolation recalled in Proposition \ref{prop_interpolation_Lp}, $\partial_{x_i}u$ belongs to $L^{p_\theta}(\RR^3,(Y,D(B^{1/2}))_{\theta,{p_\theta}})$ with $1/p_\theta=\theta/2 + (1-\theta)/6$. Due to Propositions \ref{prop_interpolations} and \ref{prop_interpolation_autoadjoint}, this shows that 
$\partial_{x_i}u$ belongs to $L^{p_\theta}(\RR^3,D(B^\alpha))$ for any $\alpha<\theta/2$. 
The same interpolation shows that $u$ belongs to $L^{p_\theta}(\RR^3,(D(B^{1/2}),D(B))_{\theta,{p_\theta}})\subset L^{p_\theta}(\RR^3,D(B^{\alpha}))$. Thus, $u$ belongs to $W^{1,p_\theta}(\RR^3,D(B^{\alpha}))$ for all $\alpha<\theta/2$. If we choose $\theta<1/2$, then we have $p_\theta>3$ and the embeddings of Theorem \ref{th_embeddings} concludes.
\end{demo}
\begin{prop}\label{prop_domain_Asalpha}
The domains of the fractional power $A_\s^\alpha$ of $A_\s$ satisfy
\[D(A_\s^{1/2})=H^1(\RR^d,Y)\cap L^2(\RR^d,D(B^{1/2}))\]
and, for all $\alpha \in [0,1]$, 
\begin{equation}\label{eq_prop_sectoriel}
D(A_\s^\alpha) \hookrightarrow  L^2(\RR^d,D(B^{\alpha})).
\end{equation}
\end{prop}
\begin{demo}
The domain of $D(A_\s^{1/2})$ is the domain $Q$  of the quadratic form of Lemma \ref{lemme_quadratic}. According to Proposition \ref{prop_interpolation_autoadjoint}, we have
\[D(A_\s^\alpha)=(L^2(\RR^d,Y),D(A_\s))_{\alpha,2}~.\]
As $D(A_\s)\hookrightarrow  L^2(\RR^d,D(B))$, using Propositions \ref{prop_interpolation_Lp} and \ref{prop_interpolation_autoadjoint} in appendix, we get
\[D(A_\s^\alpha) \hookrightarrow (L^2(\RR^d,Y),L^2(\RR^d,D(B)))_{\alpha,2} = L^2(\RR^d,(Y,D(B))_{\alpha,2})= L^2(\RR^d, D(B^{\alpha}))~.\]
\end{demo}

\begin{coro}\label{coro_domain_Asalpha_bis}
Assume that $d=1$, $2$ or $3$. For all $\alpha\in \big[0,1-\frac d4\big)$ and any $\beta \in \big[\frac 12,1\big]$ with $\beta>\alpha+\frac d4$, we have  
\[D(A_\s^\beta) \hookrightarrow L^\infty(\RR^d,D(B^\alpha)).\]
\end{coro}
\begin{demo}
The case $\beta=1$ is shown in Proposition \ref{prop_domain_A}, let us assume below that $\beta<1$.

Assume first that $d=1$. Using the above results, we know that $D(A_\s)\subset L^\infty(\RR,D(B^{3/4}))$. We also have that $D(A_\s^{1/2})\subset H^1(\RR,Y)\cap L^2(\RR,D(B^{1/2}))$. Due to Proposition \ref{prop_interpolation_Amann}, this yields $D(A_\s^{1/2})\subset L^\infty(\RR,D(B^{1/4}))$. For any $\beta\in [1/2,1)$, using the interpolation results of Section \ref{section_interpolation} with $\theta=2\beta-1$, we have $D(A_\s^\beta)\subset L^\infty(\RR,D(B^\alpha))$ for any $\alpha< \frac 34 (2\beta-1)+\frac 14 (2-2\beta) =\beta - \frac 14$.

Assume now that $d=3$. We follow the proof of Proposition \ref{prop_domain_A}. We have shown there that $D(A_\s)\subset W^{1,p_\theta}(\RR^3,D(B^{\alpha'}))$ for all $\alpha'<\theta/2$ with $\theta\in [0,1]$. We also know that $D(A_\s^{1/2})\subset H^1(\RR^3,Y)$. According to the results of Section \ref{section_interpolation}, the interpolation between both embeddings shows that $D(A_\s^\beta)\subset W^{1,q}(\RR^3,D(B^{\alpha}))$ with $\frac 1q=\frac{2\beta-1}{p_\theta} + \frac {2-2\beta}{2}$ and $\alpha < (2\beta-1)\alpha'<\theta (\beta-\frac 12)$. To conclude using the Sobolev embeddings of Theorem \ref{th_embeddings}, we need to be able to choose $q>3$, that is 
\[\beta\in \left[\frac 12,1\right]~,~~\theta\in [0,1]~,~~2\alpha<\theta(2\beta-1) ~\text{ and }~(2\beta-1)\left(\frac\theta 2 + \frac {1-\theta}6\right)+(1-\beta) < \frac 13.\]
The last inequality can be written
\[(2\beta-1)\left(\frac16 + \frac \theta 3 \right)- \frac 12 (2\beta-1) + \frac 12 < \frac 13,\]
leading to $(2\beta-1)(1-\theta)>\frac 12$, which is equivalent to $\theta < \frac{2\beta -\frac 32}{2\beta - 1}$.
Together with $2\alpha<\theta(2\beta-1)$, this yields the constraint $\alpha<\beta-\frac 34$. We can check that, if this last inequality holds, it is always possible to find a suitable $\theta\in [0,1]$.

Finally, assume that $d=2$ and fix $\alpha\in [0,1/2)$. We can follow the same  arguments above, except that we simply need the weaker constraint $q>2$. We only have to ensure that $2\alpha<\theta(2\beta-1)$, which can be done with $\theta$ close to $1$ as soon as $\beta>\alpha+\frac 12$.     
\end{demo}

\begin{coro}\label{coro_domain_Asalpha_ter}
Assume that $d=1$, $2$ or $3$. For all $\alpha,\beta \in \big[0,1]$ and $p\geq 2$ with $\beta>\alpha+\frac d4-\frac d{2p}$, we have 
\[D(A_\s^\beta) \hookrightarrow L^p(\RR^d,D(B^\alpha)).\]
\end{coro}
\begin{demo}
We know that $D(A_\s)\hookrightarrow L^\infty(\RR^d,D(B^\gamma))$ for all $\gamma\in[0,1-d/4)$, as shown in Proposition \ref{prop_domain_A}. Since we also have $D(A_\s)\hookrightarrow L^2(\RR^d,D(B))$, the interpolation results of Appendix \ref{section_interpolation} show that 
\[D(A_\s)\hookrightarrow  \big( L^\infty(\RR^d,D(B^\gamma)), L^2(\RR^d,D(B))\big)_{\theta,p_\theta}\hookrightarrow L^q(\RR^d,D(B^{\gamma'}))\] 
with $\theta =2/q$ and $\gamma'$ smaller than $(1-\frac 2q)(1-\frac d4)+\frac 2q=1-\frac d4 + \frac d{2q}$ (showing the result for $\beta=1$). We interpolate again the above embedding with $D(A_\s^0)=L^2(\RR^d,Y)$ and the power $\beta\in(0,1)$ satisfying $\frac 1p=\frac {1-\beta}2 + \frac{\beta}q$ to get that 
\[D(A_\s^\beta )\hookrightarrow  L^p(\RR^d,D(B^{\alpha}))\]
as soon as 
\[\alpha<\beta\big(1-\frac d4 + \frac d{2q}\big)=\beta + \frac d2 \big(\frac 1p - \frac 12 \big) .\]
\end{demo}

\subsection{The linear operators in weighted spaces}\label{section_Aw}
In this subsection, $\rho$ denotes a smooth weight in $\Cc^2(\RR^d,\RR_+^*)$ satisfying 
\begin{equation}\label{hyp_rho}
\exists C_\rho>0~,~~\forall x\in\RR^d~,~~ \sum_{\nu \in \NN^d, |\nu|\leq 2} |\partial_x^\nu \rho(x)|\leq C_\rho \,\rho(x). 
\end{equation}
For example, $\rho\equiv 1$, $\rho(x)=e^{\lambda \sqrt{1+|x|^2}}$ or 
\[\rho~:~x\in\RR^d~\longmapsto \frac 1{\cosh(|x|)}\]
are possible choices. 
Let $\alpha \geq 0$ and $k\in\NN$, we consider the weighted Sobolev spaces $H^k_\rho(\RR^d,D(B^\alpha))$ as introduced in Section \ref{section_basic_spaces}. If the estimate \eqref{hyp_rho} holds, then $\phi:u\mapsto \sqrt\rho u$ is an isomorphism from $H^k_\rho(\RR^d,D(B^\alpha))$ into $H^k(\RR^d,D(B^\alpha))$ for $k\leq 2$ (and even an isometry for $k=0$). We can use this isomorphism to transpose the study of the ``flat'' case in the weighted spaces.

To the main terms of the operator of the previous section, we may add lower-order terms. We include this possibility by considering a family $L(x)$ of linear operators uniformly bounded from $D(B^{1/2})\times Y^{d}$ into $Y$ and assuming that there exists $C_L>0$ such that
\begin{equation}\label{hyp_L}
\forall x\in\RR^d,~\forall (u,v) \in D(B^{1/2})\times Y^{d},~|L(x)(u,v)|_Y \leq C_L \big(|u|_{D(B^{1/2})}+|v|_{Y^d}\big).
\end{equation}

The operator in the weighted spaces is constructed as follows.
\begin{prop}\label{prop_Aw_sectoriel}
Assume that the weight $\rho$ satisfies \eqref{hyp_rho}, that $a\in\Cc^1(\RR^d,\Lc(Y^d))$ is a symmetric matrix of diffusion coefficients satisfying \eqref{hyp_a1} and \eqref{hyp_a2} and that $L(x)$ is a family of linear operators satisfying \eqref{hyp_L}. Then, the linear operator defined by
\[A_\w u= -\div \big( a(x)\grad_x u \big) + Bu + L(x)(u,\grad_x u)\] 
with domain 
\[D(A_\w)=H^2_\rho(\RR^d,Y)\cap L^2_\rho(\RR^d,D(B))\]
is a sectorial operator on $L^2_\rho(\RR^d,Y)$.

Moreover, the domain norm $\|\cdot\|_{D(A_\w)}:=\|A_\w\cdot\|_{L^2_\rho(\RR^d,Y)}$ is equivalent to the norm $\|\cdot\|_{H^2_\rho(\RR^d,Y)\cap L^2_\rho(\RR^d,D(B))}$. In addition, the constants of this equivalence, the sector $S_{z_0,\phi}$ and the constant $M$ of Definition \ref{defi_secto} in appendix, only depend on  $M_a$ and $m_a$ of \eqref{hyp_a1} and \eqref{hyp_a2}, on $C_L$ of \eqref{hyp_L} and on $C_\rho$ of \eqref{hyp_rho}. 
\end{prop}
\begin{demo}
We use the conjugacy by $\phi :u\mapsto \sqrt\rho u$ to go back to the ``flat'' spaces. We notice that $\tilde A_\s := \phi \circ A_\w \circ \phi^{-1}$ is an operator of the form
\begin{align*}
 \tilde A_\s u ~=~&  -\div \big( a(x)\grad u \big)~+~ Bu ~+~ \frac 1\rho a(x)\grad \rho  \grad u ~+~  L(x)(u,\grad u-\frac 1{2\rho} u\grad \rho)\\
   & +  \left(\frac 1{2\rho} \div(a(x)\grad \rho)   - \frac 3{4\rho^2} (a(x)\grad \rho)\cdot\grad \rho \right)u
\end{align*}
defined from $H^2(\RR,Y)\cap L^2(\RR,D(B))$ to $L^2(\RR,Y)$. The operator $A_\s$ introduced in Proposition \ref{prop_operateur_plat} is a positive self-adjoint operator with dense domain, thus a sectorial operator (Proposition \ref{prop_sa_secto1} in appendix). The above operator $\tilde A_\s$ is a perturbation of $A_\s$ by lower-order terms with bounded coefficients (recall Assumption \eqref{hyp_rho} on $\rho$). Indeed, using Corollary \ref{coro_H2} in appendix and the equivalence of norms stated in Proposition \ref{prop_operateur_plat}, we obtain that 
$$\|u\|_{H^1(\RR,Y)} \leq \varepsilon \|u\|_{H^2(\RR,Y)} + \frac {d}{4\varepsilon}\|u\|_{L^2(\RR,Y)} \leq  C \varepsilon \|A_\s u \|_{L^2(\RR,Y)} + \frac {d}{4\varepsilon}\|u\|_{L^2(\RR,Y)} $$
with $\varepsilon>0$ as small as wanted. This shows that the differential terms of order one are simply perturbative. 
In the same way, Assumption \eqref{hyp_L} implies
\begin{align*}
\|L(x)(u,\grad u&-\frac 12 u\grad \rho)\|_{L^2(\RR,Y)} \leq C_L \left(\|u\|_{L^2(\RR,D(B^{1/2}))} + \|u\|_{H^1(\RR,Y)} + \frac {C_\rho}2 \|u\|_{L^2(\RR,Y)}\right)\\
& \leq  C_L \varepsilon \left(\|u\|_{H^2(\RR,Y)}+\|u\|_{L^2(\RR,D(B))}\right) +
\left(\frac d{4\varepsilon} + \frac {C_\rho}2 + \frac 1{4\varepsilon} \right) \|u\|_{L^2(\RR,Y)}\\
& \leq  C \varepsilon \|A_\s u \|_{L^2(\RR,Y)} + \frac {C}{\varepsilon}\|u\|_{L^2(\RR,Y)}
\end{align*}

Thus, a classical perturbation argument, as Proposition \ref{prop_sa_secto2} recalled in appendix, yields that $\tilde A_\s$ is also a sectorial operator with same domain. By conjugacy, $A_\w$ is also a sectorial operator. 

To ensure the uniformities stated in the second part of Proposition \ref{prop_Aw_sectoriel}, we simply have to consider more carefully the above arguments. All the coefficients in the expression of $\tilde A_\s-A_\s$ computed above can be bounded in terms of $C_\rho$, $M_a$ and $m_a$. As noticed in the proof of Proposition \ref{prop_sa_secto2}, this shows that the sector $S_{z_0,\phi}$ and the constant $M$ of Definition \ref{defi_secto} may be chosen to depend on these numbers only. 

We also notice that the mapping $\phi :u\mapsto \sqrt\rho u$ is an isometry from $L^2_\rho(\RR^d,D(B))$ to $L^2(\RR^d,D(B))$. If we consider it from $H^2_\rho(\RR^d,Y)$ to $H^2(\RR^d,Y)$, its norm of linear mapping and the one of its inverse only depends on $C_\rho$. Thus, the equivalence between the domain norm and $\|\cdot\|_{H^2_\rho(\RR^d,Y)\cap L^2_\rho(\RR^d,D(B))}$ is also uniform with respect to $C_\rho$, $M_a$, $m_a$ and $C_L$. Finally, the definitions of the semigroup  and of the fractional power, recalled in \eqref{defi_semigroup2} and \eqref{defi_fracto} in appendix, only involve the sector $S_{z_0,\phi}$ and the inverse $(\tilde A_\s + \lambda \id)^{-1}$. Therefore, every estimate on $\tilde A_\s$ can be made only depending on $C_\rho$, $M_a$, $m_a$ and $C_L$. By conjugacy, this is also the case for the spectrum of $A_\w$.
\end{demo}

With the low-order terms represented by $L$, our main operator $-\div(a(x)\grad\cdot)+B+L(x)(\cdot,\grad\cdot)$ is somehow general if we are considering second order PDE with respect to $x$. Of course, we could also add more general perturbations terms as fractional powers of $B^{\alpha}$ with $\alpha > 1/2$. But the purpose of this paper is to present a simple whole picture of the method rather than providing the most accurate results. It seems the right place to underline that there is a main feature, which is rigid: the fact that the main order terms commute. For example, our proofs fail if we simply change $B$ into $h(x)B$ even with smooth $h$. It is surely possible to study this case, but some arguments and strong hypotheses have to be added. This is not surprising since the fact that two operators $A$ and $B$ commute is an important feature when studying the sum $A+B$ and its domain, see \cite{da-Prato-Grisvard} and \cite{DV}.

We recover the same embeddings for the domains of the powers of $A_\w$ as the ones for the case of $A_\s$.
\begin{prop}\label{prop_domain_Awalpha}
Assume that $\lambda\in\RR$ is large enough, such that the spectrum of $(A_\w +\lambda \id)$ has positive real part. 
The domains of the fractional powers $(A_\w+\lambda \id)^\alpha$ of $(A_\w+\lambda \id)$ satisfy
\[D((A_\w+\lambda \id)^{1/2})=H^1_\rho(\RR^d,Y)\cap L^2_\rho(\RR^d,D(B^{1/2}))\]
and for all $\alpha \in (0,1)$, 
\begin{equation}\label{eq_prop_sectoriel2}
D((A_\w+\lambda \id)^\alpha) \hookrightarrow  L^2_\rho(\RR^d,D(B^{\alpha})).
\end{equation}
Moreover, if $d=1$, $2$ or $3$, then for all $\alpha\in \big[0,1-\frac d4\big)$ and for any $\beta \in \big(\frac 12,1\big)$ with $\beta>\alpha+\frac d4$,
\[D((A_\w+\lambda \id)^\beta) \hookrightarrow L^\infty_{\sqrt{\rho}}(\RR^d,D(B^\alpha)).\] 
Finally, for all $\alpha,\beta \in[0,1]$ and $p\geq 2$ with $\beta>\alpha+\frac d4-\frac d{2p}$,
\[D((A_\w+\lambda \id)^\beta) \hookrightarrow L^p_{{\rho}^{p/2}}(\RR^d,D(B^\alpha)),\] 
\end{prop}
\begin{demo}
We use the conjugacy by $\phi :u\mapsto \sqrt\rho u$ and Propositions \ref{prop_domain_A} and \ref{prop_domain_Asalpha} to obtain the properties of the domain of $D((A_\w+\lambda \id)^\alpha)$. To this end, notice that the general definitions of $e^{-At}$ and $A^\alpha$, recalled in Appendix \ref{section_sectoriel}, commute with $\phi$, that is that
\[ (\tilde A_\s+\lambda \id)^\alpha := (\phi \circ (A_\w + \lambda \id) \circ \phi^{-1})^\alpha = \phi \circ (A_\w + \lambda \id)^\alpha  \circ \phi^{-1}\]
and that the perturbation arguments used in the proof of Proposition \ref{prop_Aw_sectoriel} show that $D((\tilde A_\s+\lambda \id)^\alpha)=D((A_\s)^\alpha)$, see Theorem 1.4.8 of \cite{Henry}.

Translating Corollaries \ref{coro_domain_Asalpha_bis} and \ref{coro_domain_Asalpha_ter} is also possible. We simply have to notice that the application of $\phi :u\mapsto \sqrt\rho u$ maps $L^\infty_{\sqrt{\rho}}(\RR^d,D(B^\alpha))$ into $L^\infty(\RR^d,D(B^\alpha))$ and $L^p_{{\rho}^{p/2}}(\RR^d,D(B^\alpha))$ into $L^p(\RR^d,D(B^\alpha))$.
\end{demo}


\section{The linear operators in uniformly local Sobolev spaces}\label{section_Aul}

Consider the uniformly local Sobolev spaces introduced in Section \ref{section_ul}. 
We introduce the operator 
\[A_\ulx = -\div \big( a(x)\grad \cdot \big) + B -L(x)(\cdot,\grad\cdot)\] 
defined from the domain
\[D(A_\ulx)=H^2_\ulx(\RR^d,Y)\cap L^2_\ulx(\RR^d,D(B))\]
into $L^2_\ulx(\RR^d,Y)$, where $*$ is either s or w, depending if the continuity of translations is assumed or not in the definition of the uniformly local spaces.  We will see in the present section that this operator is not necessarily a sectorial operator in the standard sense due to the lack of density of the domain. Actually, in most of our cases, $A_\ulx$ will only be an abstract parabolic operator in the sense of Definition \ref{defi_op_nondense}. 

\begin{prop}\label{prop_Aul_sectoriel_w}
Let $a\in\Cc^1(\RR^d,\Lc(Y^d))$ be a symmetric matrix of diffusion coefficients satisfying \eqref{hyp_a1} and \eqref{hyp_a2} and let $L(x)$ be a family of linear operators satisfying \eqref{hyp_L}. The operator 
$A_\ulw$ introduced above is an abstract parabolic operator in the sense of Definition \ref{defi_op_nondense}. The domain $D(A_\ulw)$ is not dense in $L^2_\ulw(\RR^d,Y)$ and thus $A_\ulw$ is never a classical sectorial operator. 
\end{prop}
\begin{demo}
Let $\mu>0$ be fixed. For any $x_\sharp\in\RR^d$, the weight $\rho_{\mu,x_\sharp}$ of \eqref{eq_rho_mu_x_sharp} satisfies the assumption \eqref{hyp_rho} since, for all $x,x_\sharp$ in $\RR^d$,
\begin{equation}\label{bound_rho_mu}
\big|\partial_{x_i} \rho_{\mu,x_\sharp}(x)\big| \leq \mu  \rho_{\mu,x_\sharp}(x)~~\text{ and }~~ \big|\partial_{x_i}\partial_{x_j} \rho_{\mu,x_\sharp}(x)\big| \leq (2\mu^2+\mu)  \rho_{\mu,x_\sharp}(x) , 
\end{equation}
so that the results of Section \ref{section_Aw} may be applied. In particular, for any $x_\sharp$, the operator $A_{\w,x_\sharp}$ defined by Proposition \ref{prop_Aw_sectoriel} for the weight $\rho=\rho_{\mu,x_\sharp}$ is a sectorial operator from $H^2_{\rho_{\mu,x_\sharp}}(\RR^d,Y)\cap L^2_{\rho_{\mu,x_\sharp}}(\RR^d,D(B))$ to $L^2_{\rho_{\mu,x_\sharp}}(\RR^d,Y)$ and the sector $S_{z_0,\phi}$ and the constant $M$ of Definition \ref{defi_secto} can be chosen independently of $x_\sharp$ because the estimate \eqref{bound_rho_mu} is independent of $x_\sharp$.

Consider any $u\in H^2_\ulw(\RR^d,Y)\cap L^2_\ulw(\RR^d,D(B))$. Due to Lemma \ref{lemme_equiv_norms}, $u$ belongs to $H^2_{\rho_{\mu,x_\sharp}}(\RR^d,Y)\cap L^2_{\rho_{\mu,x_\sharp}}(\RR^d,D(B))$ for any $x_\sharp$ and its norm can be bounded uniformly with respect to $x_\sharp$. Thus, $A_{\w,x_\sharp}u$ is well defined and belongs to $L^2_{\rho_{\mu,x_\sharp}}(\RR,Y)$ with a norm uniformly bounded with respect to $x_\sharp$. Obviously, all the functions $A_{\w,x_\sharp}u$ are the same and can be defined as $A_\ulw u$ and the uniform bound implies that $A_\ulw u$ is well defined and belongs to $L^2_\ulw(\RR^d,Y)$ due to Lemma \ref{lemme_equiv_norms}. 

Now consider any $v\in L^2_\ulw(\RR^d,Y)$ and $z\in S_{z_0,\phi}$. Since $v$ belongs to $L^2_{\rho_{\mu,x_\sharp}}(\RR,Y)$ (with a norm uniformly bounded with respect to $x_\sharp$), for all $x_\sharp\in\RR^d$, we can set $u_{x_\sharp}=(z\id-A_{\w,x_\sharp})^{-1}v$. Notice that, the different weights are equivalent in the sense that 
\[\forall x_\sharp\in\RR^d~,~~\exists M_{x_\sharp} >0 ~,~~ \frac 1{M_{x_\sharp}}  \rho_{\mu,x_\sharp} \,\leq\,\rho_{\mu,0}\,\leq\, {M_{x_\sharp}}  \rho_{\mu,x_\sharp}.\]
This implies that the solutions $u_{x_\sharp}$ all belong to the same weighted spaces and are thus the same by uniqueness, let us denote them simply $u$. Following the above considerations about uniformity, the estimate 
\[\|u\|_{L^2_{\rho_{\mu,x_\sharp}}(\RR,Y)}\leq \frac M{|z-z_0|} \|v\|_{L^2_{\rho_{\mu,x_\sharp}}(\RR,Y)}\]
holds uniformly with respect to $x_\sharp$. Lemma \ref{lemme_equiv_norms} concludes that 
\[\|u\|_{L^2_{\ul}(\RR,Y)}\leq \frac M{|z-z_0|} \|v\|_{L^2_{\ul}(\RR,Y)}.\]

We next show that, with the domain $D(A_\ulw)=H^2_\ulw(\RR^d,Y)\cap L^2_\ulw(\RR^d,D(B))$, the operator $A_\ulw$ is closed. Indeed, if $(u_n,A_\ulw u_n)$ converges to $(u_\infty,v_\infty)$ in $L^2_\ulw(\RR^d,Y)^2$, then the convergence holds in all the translations $L^2_{\rho_{\mu,x_\sharp}}(\RR^d,Y)$ of the weighted space. Proposition \ref{prop_Aw_sectoriel} shows that $A_{\w,x_\sharp}$ is closed. Thus, for all $x_\sharp\in\RR^d$, $u_\infty$ belongs to $H^2_{\rho_{\mu,x_\sharp}}(\RR^d,Y)\cap L^2_{\rho_{\mu,x_\sharp}}(\RR^d,D(B))$ and $v_\infty=A_\w u_\infty$ in $L^2_{\rho_{\mu,x_\sharp}}(\RR^d,Y)$. In Proposition \ref{prop_Aw_sectoriel}, it is also shown that the domain norm of $A_{\w,x_\sharp}$ is equivalent to the norm of $H^2_{\rho_{\mu,x_\sharp}}(\RR^d,Y)\cap L^2_{\rho_{\mu,x_\sharp}}(\RR^d,D(B))$ with constants independent of $x_\sharp$. Thus, $u_\infty$ belongs to $H^2_\ulw(\RR^d,Y)\cap L^2_\ulw(\RR^d,D(B))$ and $A_\ulw$ is closed.

To obtain a sectorial operator, we would need to show that the domain of $A_\ulw$ is dense. 
But, due to Proposition \ref{prop_dense_1_bis}, $H^2_\ulw(\RR^d,Y)$ is not dense in $L^2_\ulw(\RR^d,Y)$ (except in the trivial case $Y=\{0\}$) and thus $D(A_\ulw)$ is not dense in $L^2_\ulw(\RR^d,Y)$. 
\end{demo}

When dealing with the stronger version of the uniformly local spaces, we assume that the coefficients are invariant with respect to translations. This is not mandatory, but if we would like to include variable coefficients, some additional estimates have to be shown. For sake of simplicity, we only discuss here the case of constant coefficients. 
\begin{prop}\label{prop_Aul_sectoriel_s}
Let $a\in\Lc(Y^d)$ be a symmetric matrix of diffusion coefficients and let $L\in\Lc(D(B^{1/2})\times Y^d,Y)$ be a bounded linear operator. The operator 
$$A_\uls = -\div( a\grad\cdot) + B - L(\cdot,\grad \cdot)$$
with domain $D(A_\uls)=H^2_\uls(\RR^d,Y)\cap L^2_\uls(\RR^d,D(B))$ is an abstract parabolic operator in the sense of Definition \ref{defi_op_nondense}. The domain $D(A_\uls)$ is not necessarily dense in $L^2_\uls(\RR^d,Y)$. In particular, it is not dense if $Y$ is infinite-dimensional and $B$ has compact resolvent and, in this case, $A_\uls$ fails to be a classical sectorial operator.
\end{prop}
\begin{demo}
We consider the operator $A_\ulw$ defined by the previous Proposition \ref{prop_Aul_sectoriel_w}.
We can introduce the restriction $A_\uls$ of $A_\ulw$ to the strong version $L^2_\uls(\RR^d,Y)$ of the space. All the previous arguments estimating norms and convergences are still valid. We only have to show that the operator $A_\uls$ and its resolvent $(z\id-A_\uls)^{-1}$ preserve the continuity with respect to translations. To this end, we simply use the linearity and the invariance of the coefficients with respect to translations:  
\[(A_{\uls}u)(\cdot)-(A_{\uls}u)(\cdot-\xi)=A_{\uls}(u(\cdot))-A_{\uls}(u(\cdot-\xi))=A_{\uls}(u(\cdot)-u(\cdot-\xi))\]
and the transport of the condition on translations follows from the estimations of the norms obtained in the proof of Proposition \ref{prop_Aul_sectoriel_w}. 

We finally notice that the density of the domain fails in general. If $B$ has compact resolvent, $D(B)\hookrightarrow Y$ compactly and, if $Y$ is infinite-dimensional, then Proposition \ref{prop_dense_2} shows that $A_\uls$ is not a densely defined operator because $L^2_\uls(\RR^d,D(B))$ is not dense in $L^2_\uls(\RR^d,Y)$.
\end{demo}

The real case $Y=\RR$ was already well understood. In Theorem 2.1 of \cite{Arrieta1}, it is shown that the heat equation defines a semigroup in the uniformly local spaces $L^p_\ulw(\RR^d,\RR)$ but that the continuity at $t=0$ may fail, corresponding to Proposition \ref{prop_Aul_sectoriel_w} above. It is also shown in \cite{Arrieta1} that we can recover the continuity at $t=0$ of the heat semigroup if we work in the strong version $L^p_\uls(\RR^d,\RR)$. This corresponds to cases of Proposition \ref{prop_Aul_sectoriel_s} where the domain of $A_\uls$ is dense. We can generalize the case $Y=\RR$ studied in \cite{Arrieta1} as follows.
\begin{coro}\label{coro_Aul_sectoriel_s}
If $B$ is a linear bounded operator, i.e. $D(B)=Y$, then $A_\uls$ is a densely defined sectorial operator and thus generates a strongly continuous analytic semigroup.
\end{coro}
\begin{demo}
If $D(B)=Y$, then the domain of $A_\uls$ is simply $H^2_\uls(\RR^d,Y)$, which is dense in $L^2_\uls(\RR^d,Y)$ as shown in Proposition \ref{prop_dense_1}. This is the only property which was missing in Proposition \ref{prop_Aul_sectoriel_w} to obtain that $A_\uls$ is sectorial in the classical sense, see Appendix \ref{section_sectoriel}. 
\end{demo}

As recalled in Section \ref{section_lunardi}, the intermediate spaces $D_A(\theta,\infty)$ introduced by Sinestrari plays an important role to define solutions of abstract parabolic evolution equations, when the operator is not sectorial, see \cite{da-Prato-Sinestrari,Sinestrari}. As we have just seen, working in the strong version of the uniformly local spaces is not helpful in general. So, we focus on the weak version of these spaces.

To be able to describe the spaces $D_{A_\ulw}(\theta,\infty)$, we come back to the association with the weighted spaces. Fix $\mu>0$ and consider, for all $x_\sharp\in\RR^d$, the weight $\rho_{\mu,x_\sharp}$ of \eqref{eq_rho_mu_x_sharp}. All the results of Section \ref{section_Aw} apply. We can choose $\lambda$ large enough such that, for all $x_\sharp\in\RR^d$ and $\alpha\in\RR$, $(A_\w+\lambda \id)^\alpha$ is a well defined fractional power of the sectorial operator $(A_\w+\lambda \id)$. For $\alpha\geq 0$, we introduce the space $X_{\alpha,\infty}$ as all the functions of $X=L^2_\ulw(\RR^d,Y)$ having a finite $\|\cdot\|_{\alpha,\infty}$-norm, where 
$$\|u\|_{\alpha,\infty}=\sup_{x_\sharp\in\RR^d} \|(A_\w+\lambda \id)^\alpha u(\cdot+x_\sharp)\|_{L^2_\rho(\RR^d,Y)}.$$
Following the previous results and discussions, we have $X_{0,\infty}=X=L^2_\ulw(\RR^d,Y)$ and $X_{1,\infty}=D({A_\ulw})=H^2_\ulw(\RR^d,Y)\cap L^2_\ulw(\RR^d,D(B))$.
\begin{prop}\label{prop_domain_inter}
For all $1>\theta>\theta-\varepsilon>0$, we have 
$$D({A_\ulw})\hookrightarrow X_{\theta,\infty} \hookrightarrow D_{A_\ulw}(\theta,\infty)\hookrightarrow X_{\theta-\varepsilon,\infty} \hookrightarrow X.$$
\end{prop}
\begin{demo}
We already noticed that $X_{0,\infty}=X$ and $X_{1,\infty}=D({A_\ulw})$. Since the inclusions of the type $X_{\theta_1,\infty} \hookrightarrow X_{\theta_2,\infty}$ are trivial, we only need to deal with the central inclusions. To this end, we use some interpolations, see Appendix \ref{section_interpolation}. To interpolate between the spaces $X_{\alpha,\infty}$, we notice that they can be seen as the closed subspaces of $L^\infty(\RR^d,D((A_\w+\lambda \id)^\alpha))$ consisting of functions such that $u(x_\sharp,\cdot)=u(0,\cdot+x_\sharp)$. Using Proposition \ref{prop_interpolation_Lp}, we know that 
\[
\Big(L^\infty(\RR^d,Y),L^\infty(\RR^d,D((A_\w+\lambda \id)))\Big)_{\theta,\infty}
=L^\infty\big(\RR^d,\big(Y,D((A_\w+\lambda \id))\big)_{\theta,\infty}\big).\nonumber
\]
Using Propositions \ref{prop_interpolations} and \ref{prop_interpolation_autoadjoint}, we obtain that, for any small $\varepsilon>0$,
$$X_{\theta,\infty} \hookrightarrow \big(X_{0,\infty},X_{1,\infty}\big)_{\theta,\infty}=\big(X,D({A_\ulw})\big)_{\theta,\infty} \hookrightarrow X_{\theta-\varepsilon,\infty}.$$
This concludes the proof since, by definition,  
$D_{A_\ulw}(\theta,\infty)=\big(X,D({A_\ulw})\big)_{\theta,\infty}$.
\end{demo}


\section{Examples of Cauchy problems}\label{section_Cauchy}
In the previous sections, we have set almost all the framework to be able to construct a solution to an abstract parabolic equation in cylinders and uniformly local Sobolev spaces. Very general results could be stated, but probably to the price of heavy notations. Thus, we choose to focus on some particular settings illustrating the method rather than risking becoming too general. Also notice that, as just said above, the strong version of the uniformly local spaces does not help in general. That is why we state all the following results with the weak version.

\subsection{Global Lipschitz nonlinearity}

Let $Y_1$ and $Y_2$ be two Banach spaces. A function $F:\RR^d\times Y_1\rightarrow Y_2$ is naturally extended into a function $\tilde F$ between cylindrical spaces by setting $\tilde F(u)(x)=F(x,u(x))$. 
\begin{prop}\label{prop_nonlin_globlip}
Let $Y_1$ and $Y_2$ be two Banach spaces. Assume that $F:\RR^d\times Y_1\rightarrow Y_2$ is a uniformly Lipschitz function in the sense that there exists $K>0$ such that 
$$\forall x\in\RR^d~,~~\forall y,y'\in Y_1~,~~|F(x,y)-F(x,y')|_{Y_2}\leq K |y-y'|_{Y_1}.$$
Also assume that, for each fixed $y\in Y_1$, the map $x\mapsto F(x,y)$ belongs to $L^2_\ulw(\RR^d,Y_2)$. Then, the canonical extension 
\[\tilde F ~:~ u\in L^2_\ulw(\RR^d,Y_1) \longmapsto \big(\,x\mapsto F(x,u(x))\,\big) \in L^2_\ulw(\RR^d,Y_2)\]
is a well defined Lipschitz function. 
\end{prop}
\begin{demo}
We use the same trick as in Section \ref{section_Aul}: this result is equivalent to obtaining the same result in the weighted spaces $L^2_{\rho_{\mu,x_\sharp}}(\RR^d,Y_i)$, for the weight $\rho_{\mu,x_\sharp}$ defined by \eqref{eq_rho_mu_x_sharp}, and with all the estimates being uniform with respect to the translations $x_\sharp$. To simplify, set $\rho:=\rho_{\mu,x_\sharp}$. All the estimates in the remaining part of this proof will be obviously uniform with respect to $x_\sharp$.

We use here some properties of the Bochner integral recalled in Appendix \ref{section_spaces}.
Let $u\in L^2_\rho(\RR^d,Y_1)$. The function $F(\cdot,u):\RR^d\rightarrow Y_2$ is well defined. Since $u$ is measurable, it is the almost everywhere limit of simple functions $(u_n)$. For each $n$, $u_n$ is the finite sum of plateau functions of the type $y \Un_{x\in A}$ and thus $F(\cdot,u_n)$ is the finite sum of functions of the type $x\mapsto F(x,y)\Un_{x\in A}$. By assumption, this last function is measurable and can be approximated by simple functions. Gathering these approximations shows that $F(\cdot,u):\RR^d\rightarrow Y_2$ is a measurable function. The uniform Lipschitz property of $F$ yields that $|F(x,u(x))|_{Y_2}\leq |F(x,0)|_{Y_2} + K |u(x)|_{Y_1}$ and thus,
\[\int_{\RR^d} \rho(x)|F(x,u(x))|_{Y_2}^2 \d x \leq  C\left( \int_{\RR^d} \rho(x) |u(x)|_{Y_1}^2 \d x +   \int_{\RR^d} \rho(x) |F(x,0)|^2_{Y_2} \d x\right).\]
By assumption, the last term is a finite constant. Due to Bochner's theorem, the measurability of $F(\cdot,u)$ and the above estimation, $F(\cdot,u)$ belongs to $L^2_\rho(\RR,Y_2)$. The Lipschitz property of $\tilde F$ is obtained directly by the estimate 
\[\int_{\RR^d} \rho(x)|F(x,u(x))-F(x,v(x))|^2_{Y_2} \d x \leq  \int_{\RR^d} K^2 \rho(x)|u(x)-v(x)|_{Y_1}^2 \d x~.\]
This concludes the proof.
\end{demo}

We recall that the operator $A_\ulw$ has been introduced by Proposition \ref{prop_Aul_sectoriel_w} and that $D_{A_\ulw}(\theta,\infty)$ are the intermediate spaces introduced in Section \ref{section_intermediate}
\begin{coro}\label{coro_globLip}
Let $\alpha\in [0,1]$. Assume that $F:\RR^d\times D(B^\alpha) \rightarrow Y$ is a uniform Lipschitz function in the sense that there exists $K>0$ such that 
$$\forall x\in\RR^d~,~~\forall y,y'\in D(B^\alpha)~,~~|F(x,y)-F(x,y')|_{Y}\leq K |y-y'|_{D(B^\alpha)}.$$
Also assume that, for each fixed $y\in D(B^\alpha)$, the map $x\mapsto F(x,y)$ belongs to $L^2_\ulw(\RR^d,Y)$. Then, for any $\beta\in(\alpha,1)$, we can extend $F$ to a well defined Lipschitz function
\[\tilde F ~:~ u\in D_{A_\ulw}(\beta,\infty) \longmapsto \big(\,x\mapsto F(x,u(x))\,\big) \in L^2_\ulw(\RR^d,Y).\]
\end{coro}
\begin{demo}
The previous proposition shows that the extension is well defined in the space $L^2_\ulw(\RR^d,D(B^\alpha))$. Due to \eqref{eq_prop_sectoriel2}, we have $X_{\alpha,\infty}\hookrightarrow L^2_\ulw(\RR^d,D(B^\alpha))$, where $X_{\alpha,\infty}$ is the auxiliary space introduced above Proposition \ref{prop_domain_inter}. It remains to use this last proposition to restrict $\tilde F$ to $D_A(\beta,\infty)$ for $\beta>\alpha$.
\end{demo}

We have now set all the framework to consider the following Cauchy problem.
\begin{theorem}\label{th_Cauchy_ul_w}
Let $d\in\NN^*$, let $a\in\Cc^1(\RR^d,\Lc(Y^d))$ be a symmetric matrix of diffusion coefficients satisfying \eqref{hyp_a1} and \eqref{hyp_a2}, let $L(x)$ be a family of linear operators satisfying \eqref{hyp_L} and let $B$ be a positive self-adjoint operator on $Y$. Let $F:\RR^d\times D(B^\alpha)\rightarrow Y$ be a globally Lipschitz function in the sense of Proposition \ref{prop_nonlin_globlip} and assume that, for each $y\in D(B^\alpha)$, $x\mapsto F(x,y)$ belongs to $L^2_{\ulw}(\RR^d,Y)$. 

Then, for all $\beta\in \big(\alpha,1\big)$ and for any initial data $u_0$ in $D_{A_\ulw}(\beta,\infty)$, or more simply in $H^2_\ulw(\RR^d,Y) \cap L^2_\ulw (\RR^d,D(B))$, the equation
\begin{equation}\label{eq_Cauchy_globLip} 
\left\{\begin{array}{l} \partial_t u(x,t)=\div(a(x) \grad u(x,t))-Bu(x,t)+L(x)(u(x,t),\grad u(x,t))+F(x,u(x,t))\\ u(t=0)=u_0  \end{array}\right. 
\end{equation}
admits a unique global classical solution in the following sense:
\begin{itemize}
\item the function $t\mapsto u(t)$ is a continuous function from $[0,+\infty)$ into $L^2_\ulw(\RR^d,Y)$ and $u(t=0)=u_0$,
\item the function $t\mapsto u(t)$ belongs to $\Cc^0((0,+\infty),H^2_\ulw(\RR^d,Y)\cap L^2_\ulw(\RR^d,D(B)))$ and to $\Cc^1((0,+\infty),L^2_\ulw(\RR^d,Y))$,
\item for all $t>0$, the first equation of \eqref{eq_Cauchy_globLip} is satisfied in the strong sense in $L^2_\ulw(\RR^d,Y)$.
\end{itemize}
\end{theorem}
\begin{demo}
We simply apply Proposition \ref{prop_lunardi_glob} to $X=L^2_\ulw(\RR^d,Y)$, $A=A_\ulw$ (see Proposition \ref{prop_Aul_sectoriel_w}) and $f=\tilde F$ (see Corollary \ref{coro_globLip}). Simply notice that the space $D_A(\theta,\infty)$ of Proposition \ref{prop_lunardi_glob} is here the space $D_{A_\ulw}(\beta,\infty)$ of Corollary \ref{coro_globLip}.
\end{demo}

\subsection{Nonlinearity of polynomial type}

When considering a parabolic equation, it is usual that the nonlinearity is not a globally Lipschitz function but is a Lipschitz function on any bounded set. To provide an example in this case, we extend the Proposition \ref{prop_nonlin_globlip} in the following context. Notice that we skip the $x$ dependence of $F$ only to lighten the discussion.

\begin{prop}\label{prop_nonlin_loclip2}
Assume that $d=1$, $2$ or $3$. Let $\alpha\in \big[0,1)$ and let $F:\RR^d\times D(B^\alpha)\rightarrow Y$ be a function such that there exist $C>0$ and $\gamma\geq 0$ such that, for all $x\in\RR^d$,
\begin{equation}\label{eq_prop_nonlin_loclip2}
\forall u,v\in D(B^\alpha),~|F(u)-F(v)|_Y \leq C (1+|u|_{D(B^\alpha)}^\gamma+|v|_{D(B^\alpha)}^\gamma) |u- v|_{D(B^\alpha)}.
\end{equation}
Then, for any $\beta \in (0,1)$ such that $\beta>\alpha+\frac d4 \frac \gamma{1+\gamma}$, the canonical extension 
\[\tilde F ~:~ u\in D_{A_\ulw}(\beta,\infty)  \longmapsto \big(\,x\mapsto F(u(x))\,\big) \in L^2_\ulw(\RR^d,Y)\]
is a well defined function and a Lipschitz map from any bounded set of $D_{A_\ulw}(\beta,\infty)$ to $L^2_\ulw(\RR^d,Y)$. 
\end{prop}
\begin{demo}
Up to translate the function, assume that $F(0)=0$. The measurability of $\tilde F(u)$ is shown as in Proposition \ref{prop_nonlin_globlip}. Its integrability will follow from the one of $\tilde F(u)-\tilde F(v)$ for general $u$ and $v$ and we focus on the corresponding estimation.

We aim at obtaining estimates in an auxiliary space $X_{\theta,\infty}$ introduced above Proposition \ref{prop_domain_inter}. Again, this is equivalent to obtaining estimates in the spaces weighted by $\rho:=\rho_{\mu,x_\sharp}$ defined by \eqref{eq_rho_mu_x_sharp}, and uniformly with respect to the translations $x_\sharp$.
Choose any $\theta$ such that $\beta>\theta>\alpha+\frac d4 \frac \gamma{1+\gamma}$.
We have  
\begin{align}
\|\tilde F&(u)-\tilde F(v)\|_{L^2_\rho(\RR^d,Y)}^2 = \int_{\RR^d} \rho(x) |F(u(x))-F(v(x))|_Y^2 \d x \nonumber\\
&\leq K \int \rho(x) \big(1+ |u(x)|_{D(B^\alpha)}^{2\gamma} +|v(x)|_{D(B^\alpha)}^{2\gamma}\big) |u(x)- v(x)|_{D(B^\alpha)}^2 \d x \label{eq_demo_prop_nonlin}
\end{align}
Due to Proposition \ref{prop_domain_Awalpha}, $D((A_\w+\lambda \id)^\theta) \hookrightarrow L^2_\rho(\RR^d,D(B^\theta))$. Since $\theta>\alpha$, we have 
\[ \int \rho(x)|u(x)- v(x)|_{D(B^\alpha)}^2 \d x  \leq \|u-v\|_{D((A_\w+\lambda \id)^\theta)}^2\leq \|u-v\|_{X_{\theta,\infty}}^2.\]
To deal with both other terms of \eqref{eq_demo_prop_nonlin}, we set $p=2(1+\gamma)$ and write
\begin{align*}
\int \rho(x) |u(x)|_{D(B^\alpha)}^{2\gamma} &|u(x)- v(x)|_{D(B^\alpha)}^2 \d x\\
&= \int \rho(x)^{2\gamma/p}|u(x)|_{D(B^\alpha)}^{2\gamma} \rho(x)^{2/p}|u(x)- v(x)|_{D(B^\alpha)}^2 \d x \\
&\leq \left(\int \rho(x) |u|_{D(B^\alpha)}^p\right)^{2\gamma/p} \left(\int \rho(x)|u-v|_{D(B^\alpha)}^p\right)^{2/p}\\
&\leq \|u\|_{L^p_{\rho}(\RR^d,D(B^\alpha))}^{2\gamma} \|u- v\|_{L^p_{\rho}(\RR^d,D(B^\alpha))}^{2}
\end{align*}
Then, we notice that Proposition \ref{prop_domain_Awalpha} shows that $D((A_\w+\lambda \id)^\theta)$ is embedded in $L^p_{{\rho}^{p/2}}(\RR^d,D(B^\alpha))$ if $\theta>\alpha+\frac d4-\frac d{2p}$. For $p=2(1+\gamma)$, this last condition is the assumed one $\theta>\alpha+\frac d4 \frac \gamma{1+\gamma}$. Notice that we cannot globally compare the weight ${\rho}^{p/2}$ to $\rho$. However, we can do it locally and then use the control by all the translations $x_\sharp$. Doing so, we obtain that 
\[\int \rho(x) |u(x)|_{D(B^\alpha)}^{2\gamma} |u(x)- v(x)|_{D(B^\alpha)}^2 \d x \leq C \|u\|_{X_{\theta,\infty}}^{2\gamma} \|u- v\|_{X_{\theta,\infty}}^{2}.\]
Coming back to \eqref{eq_demo_prop_nonlin} and gathering all the estimates and their translations, we obtain that 
\[\|\tilde F(u)-\tilde F(v)\|_{L^2_\ulw(\RR^d,Y)}^2 \leq C (1+\|u\|_{X_{\theta,\infty}}^{2\gamma} +\|v\|_{X_{\theta,\infty}}^{2\gamma})\|u- v\|_{X_{\theta,\infty}}^{2}\]
It remains to remember that Proposition \ref{prop_domain_inter} implies $D_{A_\ulw}(\beta,\infty)\hookrightarrow X_{\theta,\infty}$, showing that $\tilde F$ is a Lipschitz map from any bounded set of $D_{A_\ulw}(\beta,\infty)$ to $L^2_\ulw(\RR^d,Y)$. 
\end{demo}

We have set all the framework to apply Theorem \ref{th_lunardi}.
\begin{theorem}\label{th_Cauchy_ul_loc_w}
Assume that $d=1$, $2$ or $3$. Let $a\in\Cc^1(\RR^d,\Lc(Y^d))$ be a symmetric matrix of diffusion coefficients satisfying \eqref{hyp_a1} and \eqref{hyp_a2}, let $L(x)$ be a family of linear operators satisfying \eqref{hyp_L} and let $B$ be a positive self-adjoint operator on $Y$. Let $\alpha\in [0,1)$ and assume that $F:\RR^d\times D(B^\alpha)\rightarrow Y$ is such that \eqref{eq_prop_nonlin_loclip2} holds for some $\gamma\geq 0$ such that $\alpha+\frac d4\frac \gamma{1+\gamma}<1$. 

Then, for all $\beta\in \big(\alpha+\frac d4\frac \gamma{1+\gamma},1\big)$ and for any initial data $u_0$ in $D_{A_\ulw}(\beta,\infty)$, or more simply in $H^2_\ulw(\RR^d,Y) \cap L^2_\ulw (\RR^d,D(B))$, there exists $T>0$ such that the equation 
\begin{equation}\label{eq_Cauchy_locLip} 
\left\{\begin{array}{l} \partial_t u(x,t)=\div(a(x) \grad u(x,t))-Bu(x,t)+L(x)(u(x,t),\grad u(x,t)) + F(u(x,t))\\ u(t=0)=u_0  \end{array}\right. 
\end{equation}
admits a unique local classical solution in the following sense:
\begin{itemize}
\item the function $t\mapsto u(t)$ is a continuous function from $[0,T)$ into $L^2_\ulw(\RR^d,Y)$ and $u(t=0)=u_0$,
\item the function $t\mapsto u(t)$ belongs to $\Cc^0((0,T),H^2_\ulw(\RR^d,Y)\cap L^2_\ulw(\RR^d,D(B)))$ and to $\Cc^1((0,T),L^2_\ulw(\RR^d,Y))$,
\item for all $t\in(0,T)$, the first equation of \eqref{eq_Cauchy_locLip} is satisfied in the strong sense in $L^2_\ulw(\RR^d,Y)$.
\end{itemize}
Moreover, if the maximal time of existence $T$ is finite, then
$$\|u(t)\|_{D_{A_\ulw}(\beta,\infty)}~\xrightarrow[~t\longrightarrow T^-~]{}+\infty~~\text{ and } \|u(t)\|_{X_{\beta,\infty}}~\xrightarrow[~t\longrightarrow T^-~]{}+\infty.$$
\end{theorem}
\begin{demo}
This is a direct application of Theorem \ref{th_lunardi}, Corollary \ref{coro_lunardi_glob} and Proposition \ref{prop_nonlin_loclip2}. 
We work in the space $D_{A_\ulw}(\beta,\infty)$ for $\beta\in \big(\alpha+\frac d4\frac \gamma{1+\gamma},1\big)$. Simply notice that the second blow-up condition follows from Proposition \ref{prop_domain_inter}.  
\end{demo}


\section{The gradient structure}\label{section_gradient}

Many parabolic equations admit an energy, that is a Lyapunov functional, whose value decreases along the trajectories. In the framework of uniformly local spaces, the natural energy is often infinite. However, it is possible to obtain qualitative properties from the gradient structure, see \cite{Gallay,Muratov_Novaga_1,Risler} for example.

The purpose of this section is to show some arguments using the formal gradient structure. To fix the idea, we consider the equation 
\begin{equation}\label{eq_gradient}
\partial_t u = \Delta u - B u - \grad_u V(u)~,~~t>0
\end{equation}
where:
\begin{itemize}
\item As in the whole paper, $B$ is a positive self-adjoint operator densely defined on a Hilbert space $Y$. 
\item The nonlinear term $\grad_u V$ is defined from $D(B^\alpha)$ into $Y$ for some $\alpha\in (0,1/2]$ and is globally Lipschitz function in the sense that there exists $C>0$ such that, for all $u$ and $v$ in $D(B^\alpha)$,
\begin{equation}\label{hyp_gradV}
\big|\grad_u V(u) -\grad_u V(v)\big|_Y ~\leq~ C |u-v|_{D(B^\alpha)}.
\end{equation}
\item There exists a potential $V$ defined from $D(B^\alpha)$ into $\RR$ which is the antiderivative of $\grad_u V$ in the following sense. For all curves $\Gamma(\tau)\in\Cc^0([0,1],D(B^\alpha))\cap \Cc^1([0,1],Y)$, we have $\partial_\tau V(\Gamma(\tau))=\langle \grad_u V(\Gamma(\tau))|\partial_\tau \Gamma(\tau)\rangle_Y$.
\item The potential $V$ is coercive in the sense that there exist positive $\kappa$ and $\delta$ such that 
\begin{equation}\label{hyp_V_coercive}
\forall x\in\RR^d~,~~\forall u\in D(B^\alpha)~,~~ \langle \grad_u V(u) | u \rangle_Y \geq \kappa |u|_Y^2- \delta|u|_Y.
\end{equation}
\end{itemize}
Fix $\beta \in (1/2,1)$ and 
consider $u_0$ in $D_{A_\ulw}(\beta,\infty)$, or more simply in $H^2_\ulw(\RR^d,Y)\cap L^2_\ulw(\RR^d,D(B))$. Theorem \ref{th_Cauchy_ul_w} or Theorem \ref{th_Cauchy_ul_loc_w}, applied with $F=-\grad_u V$, ensures the existence of a solution $u(t)$ for $t$ in an interval $[0,T)$. In this section, the fact that $u(t)$ solves \eqref{eq_gradient} has to be understood in the sense of these theorems. We acknowledge that the above hypotheses may be weakened. However, the purpose of this section is to show how to make use of the formal energy structure, so we choose a simple setting to avoid technicalities. 

Formally, to the parabolic equation \eqref{eq_gradient} is associated the energy 
\begin{equation}\label{def_energy_formal}
\Ec(t)=\int_{\RR^d} \frac 12 |\grad_x u(x,t)|^2_{Y^d} + \frac 12 |B^{1/2}u(x,t)|_Y^2 + V(u(x,t))\d x,
\end{equation}
which is expected to dissipate since a formal computation provides 
\[\partial_t \Ec(t)= - \int_{\RR^d} |\partial_t u(x,t)|_Y^2 \d x.\]
The existence of this formal decreasing energy is what can be called a {\it formal gradient structure}. However, even for $V\equiv 0$, this energy $\Ec$ can be infinite because no spatial decay is assumed for $u$. The uniformly local spaces simply provide bounds on localized versions of this energy. The purpose of this section is to show that we still can use this formal structure to obtain qualitative properties. 

\begin{theorem}\label{th_gradient}
Consider the above framework.
For any $R_0>0$ and $\epsilon>0$, there exists a bound $R_1>0$ such that the following holds. For any $u_0\in H^2_\ulw(\RR^d,Y)\cap L^2_\ulw(\RR^d,D(B))$ with norm less than $R_0$, the solution $u(t)$ of \eqref{eq_gradient} with $u(t=0)=u_0$ exists for all times $t\geq 0$. Moreover, $u(t)$ stays in the ball of $H^2_\ulw(\RR^d,Y)\cap L^2_\ulw(\RR^d,D(B))$ of radius $R_1$ for all $t\in [\epsilon,+\infty)$.
\end{theorem}
\begin{proof}
In this proof, when no confusion is possible, we lighten the notations by setting $X:=L^2_\ulw(\RR^d,Y)$, $A:=A_\ulw$ and $D(A):=H^2_\ulw(\RR^d,Y)\cap L^2_\ulw(\RR^d,D(B))$. First notice that the Cauchy problem is well posed in any space $D_{A}(\beta,\infty)$ with $\beta>1/2\geq \alpha$ due to Theorem \ref{th_Cauchy_ul_w} and that the solution exists globally. The only property that has to be shown is the boundedness of the solution.  
The key idea to use the gradient structure is to introduce a truncated version of the energy. We split the proof in several steps.

\smallskip

\noindent{\it $\bullet$ Step 1: definition of the truncated energy.} 
Let $\mu>0$ to be fixed later. Let $x_\sharp\in\RR^d$ and let $\rho_{\mu,x_\sharp}$ be as in \eqref{eq_rho_mu_x_sharp}, often simply denoted $\rho$. 
Up to an addition of a constant to the potential, we assume that $V(0)=0$ in order to simplify the future estimates.
For any solution $u$, we set 
\begin{align}
 \Ec_{\mu,x_\sharp}(t)=\int_{\RR^d} \rho_{\mu,x_\sharp}(x)&\left(\frac 12|\grad_x u(x,t)|_{Y^d} + \frac 1 2 |u(x,t)|_{Y}^2\right.\nonumber \\ & \quad\quad\quad\quad \left.+ \frac 12 |B^{1/2}u(x,t)|_Y^2 + V(u(x,t))\right)\d x.\label{def_E_tronquee}
\end{align}
Remember that $\beta>1/2$, so that the 
terms $\int \rho |\grad u|^2$ and $\int \rho |B^{1/2}u|^2$ of the above energy are well defined and continuous when $t\rightarrow 0^+$ since $D_{A}(\beta,\infty)\hookrightarrow X_{1/2,\infty}:=H^1_\ulw(\RR^d,Y)\cap L^2_\ulw(\RR^d,D(B^{1/2}))$. Also notice that, applying the definition of $V$ with $\Gamma(\tau)=\tau u$ and recalling that $V(0)=0$, we have
\begin{equation}\label{eq_V_primitive}
V(u)=\int_0^1 \langle \grad_u V(\tau u)|u\rangle_Y \d \tau.
\end{equation}
Due to Assumption \eqref{hyp_gradV}, we have that $|\grad V(v)|_Y\leq |\grad V(0)|_Y + C|v|_{D(B^\alpha)}$, implying that $V(u)$ is well defined for all $u\in D(B^\alpha)$ and that 
\begin{equation}\label{eq_V_estim}
|V(u)|\leq |\grad_u V(0)|_{Y^d}|u|_Y + \frac C2 |u|_{D(B^\alpha)}|u|_Y 
 \leq C'\left(\frac 1\kappa + (1+\kappa)|u|_Y^2 + |u|_{D(B^\alpha)}^2\right)
\end{equation}
Since $\int\rho$ is finite, the term $\int \rho V(u)$ is well defined for all $u\in D(B^\alpha)$. Finally, since the Cauchy problem is set in $D_{A}(\beta,\infty)\hookrightarrow L^2_\ulw(\RR^d,D(B^\alpha))$ (see Proposition \ref{prop_domain_inter}), we obtain that the energy $\Ec$ is well defined and continuous when $t\rightarrow 0^+$.

\smallskip

\noindent{\it $\bullet$ Step 2: differentiation of the energy.} 
By assumption on $V$ and due to the regularity of the solutions stated in Theorems \ref{th_Cauchy_ul_w} and \ref{th_Cauchy_ul_loc_w}, we can differentiate the terms $\int \rho |u|^2$ and $\int \rho V(u)$. This is less direct for the two other terms. For example, nothing states that $t\mapsto B^{1/2}u(t)$ is differentiable in $Y$. The trick is as follows. Set $\Phi(t)=\int_{\RR^d} \rho(x) |B^{1/2}u(t,x)|_Y^2 \d x$. We have 
\begin{align}
\Phi(t+\delta t)-\Phi(t) &=  \int_{\RR^d} \rho  \langle B^{1/2} u(t+\delta t) |B^{1/2} u(t + \delta t) - B^{1/2} u(t)\rangle_Y \d x  \label{eq-derivation}\\&\quad   + \int_{\RR^d} \rho  \langle B^{1/2} u(t + \delta t ) - B^{1/2} u(t) |B^{1/2}u(t)\rangle_Y \d x~. \nonumber
\end{align}
Consider the first term, we have
\begin{align*}
\int_{\RR^d} \rho  \langle B^{1/2} u(t+\delta t) &|B^{1/2} (u(t + \delta t) - u(t))\rangle_Y \d x\\ &= \int_{\RR^d}\rho \langle B u(t+\delta t)) |u(t + \delta t) - u(t)\rangle_Y \d x.
\end{align*}
Now, we can use the regularities of the solution: $t\in (0,T]\mapsto u(t)$ is continuous in $L^2_\ulw(\RR^d,D(B))$ and $\partial_t u(t)$ exists and belongs to $L^2_\ulw(\RR^d,Y)$. Thus, when $\delta t\rightarrow 0$,
$$\frac 1{\delta t} \int_{\RR^d} \rho  \langle B^{1/2} u(t+\delta t) |B^{1/2} u(t + \delta t) - B^{1/2} u(t)\rangle_Y \d x ~~\longrightarrow~~ \int_{\RR^d} \rho \langle B u(t) |\partial_t u(t)\rangle_Y \d x.$$
The second term of \eqref{eq-derivation} is similar and we obtain the natural computation:
$$\partial_t \int_{\RR^d} \rho(x)|B^{1/2}u(x,t)|_Y^2 \d x = 2 \int_{\RR^d}\rho(x)\langle Bu(x,t)|\partial_t u(x,t) \rangle_Y.$$
The term $\int \rho|\grad_x u|^2$ is differentiated in the same way, using the integration by parts stated in Proposition \ref{prop_inte_parts} in appendix.
Following these arguments, we can compute 
\[\partial_t \Ec_{\mu,x_\sharp}(t) = - \int_{\RR^d} \rho|\partial_t u|_{Y}^2 - \int_{\RR^d}\langle \grad_x u|(\grad_x \rho)\partial_t u\rangle_{Y^d} + \int_{\RR^d}\rho \langle u|\partial_t u\rangle_{Y^d}.\]
Then, we use the evolution equation \eqref{eq_gradient} to replace $\partial_t u$ in the last term and we obtain
\begin{align}
\partial_t \Ec_{\mu,x_\sharp}(t) = & - \int_{\RR^d} \rho|\partial_t u|_{Y}^2 - \int_{\RR^d}\langle \grad_x u|(\grad_x \rho)\partial_t u\rangle_{Y^d}  -   \int_{\RR^d}\rho| \grad_x u|^2_{Y^d} \label{deriv_energy} \\
& -  \int_{\RR^d}\rho|B^{1/2}u|_Y^2  -  \int_{\RR^d}\langle \grad_x u|(\grad_x \rho)u\rangle_{Y^d}  - \int_{\RR^d}\rho \langle u|\grad_u V(u)\rangle_Y .\nonumber
\end{align}

\smallskip

\noindent{\it $\bullet$ Step 3: bound of the energy} 
Now, we use the coercivity of $V$, assumed in \eqref{hyp_V_coercive}, and the control \eqref{bound_rho_mu} of the derivative of the weight. We obtain that 
\begin{align}
\partial_t &\Ec_{\mu,x_\sharp}(t) \leq 
- \int_{\RR^d} \rho_{\mu,x_\sharp}|\partial_t u|_{Y}^2 + \mu  \int_{\RR^d} \rho_{\mu,x_\sharp} |\grad_x u|_{Y^d} |\partial_t u|_{Y}  -  \int_{\RR^d}\rho_{\mu,x_\sharp}|\grad_x u|_{Y^d}^2   
\nonumber\\ &- \int_{\RR^d}\rho_{\mu,x_\sharp}|B^{1/2}u|_Y^2 +  \mu  \int_{\RR^d} \rho_{\mu,x_\sharp} |\grad_x u|_{Y^d} |u|_{Y} -  \kappa  \int_{\RR^d}\rho_{\mu,x_\sharp} |u|_Y^2 +\delta \int_{\RR^d}\rho_{\mu,x_\sharp}|u|_Y .\nonumber
\end{align}
We  can bound the last term $\delta\int\rho |u|_Y$ by $\int\rho(\frac {2\delta^2}\kappa+\frac{\kappa}2|u|_Y^2)$ since $\int\rho$ is finite (even if blowing-up when $\mu\rightarrow 0$).
Since, $\alpha\leq 1/2$, we have $X_{1/2,\infty}:=H^1_\ulw(\RR^d,Y)\cap L^2_\ulw(\RR^d,D(B^{1/2}))\hookrightarrow L^2_\ulw(\RR^d,D(B^\alpha))$. From \eqref{eq_V_estim}, we then obtain 
\begin{equation}\label{eq_V_estim2}
\int_{\RR^d}\rho_{\mu,x_\sharp} V(u) \leq C' \int_{\RR^d} \rho_{\mu,x_\sharp} \left(1+|\grad_x u|_{Y^d}^2 + |B^{1/2}u|_Y^2 \right)\d x.
\end{equation}
We obtain that, up to choose $\mu>0$ small enough, there exists $\nu>0$ such that
\begin{equation*}\label{dec_Ecal}
\partial_t \Ec_{\mu,x_\sharp}(t) \leq - \nu \Ec_{\mu,x_\sharp}(t) + \frac 1\nu.
\end{equation*}
Gr\"onwall lemma and the continuity of $\Ec_{\mu,x_\sharp}(t)$ at $t=0$ yield
\begin{equation}\label{dec_Ecal2}
\forall t\geq 0~,~~\Ec_{\mu,x_\sharp}(t) \leq e^{-\nu t} \Ec_{\mu,x_\sharp}(0) + \frac 1{\nu^2}.
\end{equation}
On the other hand, using \eqref{eq_V_primitive}, we have that 
$$V(u)= \int_0^1 \frac 1\tau \langle \grad_u V(\tau u)|\tau u\rangle_Y \d \tau \geq  \frac 12 \kappa |u|_Y^2 - \delta |u|_Y,$$
showing that $V$ is bounded from below. Thus,  
\[\Ec_{\mu,x_\sharp}(t) \geq  \int_{\RR^d}  \rho_{\mu,x_\sharp}(x) \left(\frac 12 |\grad_x u(x,t)|_{Y^d}^2 + \frac 1 2 |u(x,t)|_{Y}^2 + \frac 12 |B^{1/2}u(x,t)|_Y^2 + \min V \right)\d x.\]
Using the fact that all the previous computations are uniform with respect to the translation vector $x_\sharp$, we obtain the existence of $m_\Ec>0$ such that 
\begin{equation}\label{E_et_H1}
\sup_{x_\sharp\in\RR^d}\Ec_{\mu,x_\sharp}(t) \geq m_{\Ec} \left(\|u(t)\|_{H^1_\ul(\RR^d,Y)\cap L^2_\ul(\RR^d,D(B^{1/2}))} - 1\right).
\end{equation}
The bounds \eqref{dec_Ecal2} and \eqref{E_et_H1} imply that the $X_{1/2,\infty}$-norm of the solution remains bounded by a constant only depending on $\sup_{x_\sharp\in\RR^d}\Ec_{\mu,x_\sharp}(0)$, which is bounded.

\smallskip

\noindent{\it $\bullet$ Step 4: use of the regularization effect.} We just proved that the $X_{1/2,\infty}$-norm of the solutions with $\|u_0\|_{D_A(\beta,\infty)}\leq R_0$ is bounded by a quantity, let us say $R_{1/2}$. It remains to obtain the boundness of the $X_{1,\infty}-$norm, that is the norm of $D(A_\ulw)=H^2_\ulw(\RR^d,Y)\cap L^2_\ulw(\RR^d,D(B))$.
Let $\epsilon>0$, we aim at using the regularization effect during this time $\epsilon$. For all $t\geq \epsilon$,
\begin{equation}\label{duhamel}
u(t)=e^{-A\epsilon}u(t-\epsilon) +\int_0^{\epsilon} e^{-A s}\tilde F(u(t-s))\d s
\end{equation}
where $\tilde F$ is the canonical extension of $-\grad_u V$ provided by Corollary \ref{coro_globLip}. Note that the uniform bound $\|u(t)\|_{X_{1/2,\infty}}\leq R_{1/2}$ implies a uniform bound $\|\tilde F(u(t-s))\|_{Y}\leq R_F$ due to \eqref{hyp_gradV}.
Using Propositions \ref{prop_regu_integrated} and \ref{prop_estimation_semigroup_interm}, \eqref{duhamel} implies that $u(t)$ belongs to $D(A)$ and 
$$\forall t\geq\epsilon~,~~\|A u(t)\|_{X}\leq C \epsilon^{-1/2} R_{1/2} + C R_F=:R_1.$$
\end{proof}

\section{Compactness properties}\label{section_compact}

In this last section, we discuss the compactness properties. The compact embeddings in vector-valued Sobolev spaces has been studied by Amann, see \cite{Amann} for example. We will simply apply them to our context. In this sense, this section is simply here for sake of completeness and for gathering materials for possible future uses. 

In many possible applications, the positive self-adjoint operator $B$ has compact resolvents. In this case, $D(B^\alpha)\doublehookrightarrow D(B^\beta)$ as soon as $\alpha>\beta$, where in this section $Y_2\doublehookrightarrow Y_1$ means that the embedding is compact. Then, we may expect some compactness properties for the solutions of the abstract parabolic equation. Since we consider an unbounded space variable $x\in\RR^d$, the compactness is only expected as a spatially local property, but may be still of interest in several applications.

We have not introduced the Sobolev spaces with fractional regularity and we will not apply the result in these general cases. However, we recall here the exact theorem for completeness.

\begin{theorem}[Theorem 5.2 of \cite{Amann}]\label{th_Amann}
Let $\Omega\subset \RR^d$ be a smooth bounded domain. Let $E_0$ and $E_1$ be two Banach spaces such that $$E_1\doublehookrightarrow E_0~.$$ 
Let $s_0,s_1\in\RR$ and let $p_0,p_1\in [1,\infty]$. Let $\theta\in (0,1)$, we set $p_\theta$ as in Proposition \ref{prop_interpolation_Lp} and  $s_\theta:=(1-\theta)s_0 + \theta s_1$. Finally, assume that $E$ is a Banach space satisfying 
$$(E_0,E_1)_{\theta,p_\theta}\hookrightarrow E \hookrightarrow E_0~.$$

Then, for all $s<s_\theta$ with $s-d/p<s_\theta - d/p_\theta$, we have 
\begin{equation}\label{eq_th_Amann1}
W^{s_0,p_0}(\Omega,E_0) \cap W^{s_1,p_1}(\Omega,E_1) ~ \doublehookrightarrow ~ W^{s,p}(\Omega,E)~. 
\end{equation}
Moreover, if $0\leq s < s_\theta-d/p_\theta$ then 
\begin{equation}\label{eq_th_Amann2}
W^{s_0,p_0}(\Omega,E_0) \cap W^{s_1,p_1}(\Omega,E_1) \doublehookrightarrow \Cc^s (\overline \Omega,E)
\end{equation}
where $\Cc^s$ is the space of bounded continuous functions with $s-$H\"older regularity.
\end{theorem}

Our applications of the above result are as follows.
\begin{coro}\label{coro_Amann}
Assume that the positive self-adjoint operator $B:D(B)\rightarrow Y$ has compact resolvents. Let $\Omega\subset \RR^d$ be a smooth bounded domain. Then,  
\begin{equation}\label{eq_coro_Amann1}
\forall \alpha\in[0,1)~,~~ H^2(\Omega,Y)\cap L^2(\Omega,D(B)) \doublehookrightarrow L^2(\Omega,D(B^\alpha))
\end{equation}
and, if $d=1$, $2$ or $3$,  
\begin{equation}\label{eq_coro_Amann2}
\forall \alpha\in\big[0,1-\frac d4\big)~,~~ H^2(\Omega,Y)\cap L^2(\Omega,D(B)) \doublehookrightarrow \Cc^0(\overline{\Omega},D(B^\alpha)).
\end{equation}
\end{coro}
\begin{demo}
We apply Theorem \ref{th_Amann} with:
\begin{itemize}
\item $E_1=D(B)$ and $E_0=Y$. Since $B$ has compact resolvents, we have $E_1\doublehookrightarrow E_0$.
\item $p_0=p_1=p=p_\theta=2$, $s_0=2$ and $s_1=0$,  
\item $\theta\in(0,1)$ with $\theta>\alpha$, so that $(E_0,E_1)_{\theta,p_\theta}=(Y,D(B))_{\theta,2}=D(B^{\theta})$ (see Proposition \ref{prop_interpolation_autoadjoint}), this interpolated space being compactly embedded in $D(B^\alpha)$.
\end{itemize}
Applying \eqref{eq_th_Amann1} with $s=0$, we obtain the compact embedding \eqref{eq_coro_Amann1}. To apply \eqref{eq_th_Amann2} to obtain \eqref{eq_coro_Amann2}, we need to have $s_\theta-d/p_\theta>0$, that is $2(1-\theta)-d/2>0$. Being able to choose a such $\theta\in (0,1)$ with $\theta>\alpha$ is equivalent to have $\alpha<1-d/4$.
\end{demo}

\begin{prop}\label{prop_compact}
Assume that the positive self-adjoint operator $B:D(B)\rightarrow Y$ has compact resolvents. Let $t\in [0,+\infty) \mapsto u(t)$ be a trajectory bounded in $H^2_\ulw (\RR^d,Y)\cap L^2_\ulw(\RR^d,D(B))$, uniformly with respect to $t\geq 0$. Then, for any sequence $(t_n)$ going to $+\infty$, there exists a subsequence $(t_{\varphi(n)})$ and a function $u_\infty\in H^2_\ulw (\RR^d,Y)\cap L^2_\ulw(\RR^d,D(B))$ such that, for all bounded set $\Omega\subset\RR^d$ and $\alpha\in [0,1)$,
\[u(t_{\varphi(n)})\xrightarrow[~~n\longrightarrow+\infty~~]{}u_\infty~~\text{ in }L^2(\Omega,D(B^\alpha)).\]
If moreover $\alpha<1-d/4$, then this convergence also holds in $\Cc^0(\overline{\Omega},D(B^\alpha))$.
\end{prop}
\begin{demo}
Up to an argument of diagonal extraction, it is sufficient to fix a bounded smooth set $\Omega\subset\RR^d$ and to study the compactness of a bounded sequence $(v_n)$ of $H^2(\Omega,Y)\cap L^2(\Omega,D(B))$. Since $(v_n)$ is bounded in the Hilbert space $H^2(\Omega,Y)\cap L^2(\Omega,D(B))$, up to extracting a subsequence, we can assume that $(v_n)$ weakly converges to a function $v_\infty\in H^2(\Omega,Y)\cap L^2(\Omega,D(B))$. The weak convergence ensures that the norm of $v_\infty$ is bounded by the same bound as the norms of $(v_n)$. Then, we use the compactness embeddings of Corollary \ref{coro_Amann} to obtain that, again up to extracting a subsequence, we can assume that the convergence to $v_\infty$ is strong in $L^2(\Omega,D(B^\alpha))$ and even in $\Cc^0(\overline{\Omega},D(B^\alpha))$ if $0\leq \alpha<1-d/4$.

To obtain the proposition from the above argument, we use a diagonal extraction of the sequence $u(t_n)$ to get the same type of convergence in each bounded set $\Omega$, with a limit $u_\infty\in H^2_\loc (\RR^d,Y)\cap L^2_\loc(\RR^d,D(B))$ independent of $\Omega$. It remains to notice that the local bounds of $u_\infty$ in $H^2(\Omega,Y)\cap L^2(\Omega,D(B))$ obtained by the arguments above follows from the local bounds on $u(t_n)$. Thus, since $u(t_n)$ is bounded in $H^2_\ulw (\RR^d,Y)\cap L^2_\ulw(\RR^d,D(B))$, we have that $u_\infty$ also belongs to $H^2_\ulw (\RR^d,Y)\cap L^2_\ulw(\RR^d,D(B))$.
\end{demo}



\section{Examples of application}\label{section_appli}
In this section, we show how to recover some natural evolution equations from our abstract equation
\begin{equation}\label{eq_appli} 
\partial_t u =\div(a(x) \grad u)-Bu+L(x)(u,\grad u)+F(u).
\end{equation}
We also discuss some of our results in this perspective.

\subsection{Classical real-valued parabolic systems}
Let $k\geq 1$ and $d\geq 1$ be natural integers. We consider the framework of the present paper with $Y=\RR^k$, $B=\id$, $L(x)(u,v)=u$ and $F(x,u)=f(x,u(x))$ with $f\in\Cc^1(\RR^{d}\times \RR,\RR^k)$. Then \eqref{eq_appli} becomes a classical parabolic system
\begin{equation}\label{eq_appli1}
\partial_t u(x,t) =\div(a(x) \grad u(x,t))  +f(x,u(x,t)). 
\end{equation}
where $u=(u_i)_{i=1\ldots k}$ is valued in $\RR^k$. This paper shows that we can apply the theory of abstract parabolic operators with non-dense domain to set a Cauchy problem in the weak version $L^2_\ulw(\RR^d,\RR^k)$ of the uniformly local spaces. 
We can apply Proposition \ref{prop_nonlin_loclip2} to a polynomial nonlinearity $f=u^{\gamma+1}$ and we notice that the condition $\alpha+\frac d4\frac{\gamma}{\gamma+1}<1$ allows any power $\gamma$ for $d=1$, $2$ or $3$ since $\alpha=0$. This naturally corresponds to the embedding $H^2\hookrightarrow L^p$ for all $p\geq 2$ in these dimensions. To apply the gradient structure as in Theorem \ref{th_gradient}, we require $\alpha+\left(\frac d4+\frac 12\right)\frac{\gamma}{\gamma+1}<1$, which is an empty condition for $d=1$ or $d=2$ but requires $\gamma<4$ for $d=3$. This last case corresponds to a potential of the type $V(u)=u^{\gamma'}$ with $\gamma'<6$, corresponding to the embedding $H^1\hookrightarrow L^6$, which is important and natural to be able to use the energy.

Finally, notice that, due to Corollary \ref{coro_Aul_sectoriel_s}, we could also recover a classical sectorial operator by working in the strong version $L^2_\uls(\RR^d,\RR^k)$ of the uniformly local spaces, leading to a semigroup $e^{-A_\uls t}$ being continuous at $t=0$, as already shown in \cite{Arrieta1}. Then, we can apply the usual theory to obtain a well-posed nonlinear Cauchy problem and to study the behaviour of solutions, as done in \cite{Arrieta2}.

\subsection{Anisotropic topology in $\RR^{d+d'}$.}

We consider a real-valued parabolic system as \eqref{eq_appli1} but we split the spatial domain into two subspaces: for $x\in\RR^d$ we consider the uniformly local topology and for $x'\in\RR^{d'}$ we use the classical $L^2$ topology. This allows to consider solutions with infinite energy along the $x$ variable but finite energy along the $x'$ variable, as the propagating pattern of Figure \ref{fig-corridor}.
This is done as follows: choose $Y=L^2(\RR^{d'})$, $B=-\Delta+\id$ the positive Laplacian operator in $\RR^{d'}$, $a\equiv \id$, $L\equiv 0$ and $F(u)=u^{\gamma+1}+u$.
Then our abstract parabolic equation \eqref{eq_appli} becomes the classical 
$$\partial_t u=\Delta u +u^{\gamma+1}~~~~~\text{ in }\RR^{d+d'}$$
but considering solutions in the space $L^2_\ulw(\RR^d,L^2(\RR^{d'}))$ instead of the classical $L^2(\RR^{d+d'})$. To be able to apply Theorem \ref{th_Cauchy_ul_loc_w} to define solutions, we have to choose $\alpha\in[0,1)$ such that $$D(B^\alpha)=H^{2\alpha}(\RR^{d'})\hookrightarrow L^{2(\gamma+1)}(\RR^{d'})~~\text{ and }~~\alpha+\frac d4\frac \gamma{1+\gamma}<1.$$
The Sobolev embedding holds if $\frac 12 - \frac{2\alpha}{d'}\leq \frac 1{2(\gamma+1)}$ that is if $\frac{d'}4\frac{\gamma}{\gamma+1}\leq \alpha$. Thus, we can properly set a Cauchy problem in our context if $\frac{d+d'}4\frac{\gamma}{\gamma+1}<1$, which depends only on the total spatial dimension $d+d'$, as naturally expected. Also notice that this is the usual condition to be able of defining simply the solutions of the semilinear parabolic equation in the base space $L^2(\RR^{d+d'})$. This underline that the presence of a partial uniformly local topology has no impact on the powers that are admissible for the nonlinear term.

\subsection{Advective reaction-diffusion PDE in cylinders}
Fix $d=1$ and consider the cylinder $\Omega=\RR\times\omega$, where $\omega\subset\RR^2$ is a bounded smooth domain. A point of $\Omega$ is denoted by $(x,y)$ with $x\in\RR$ and $y\in\omega$. Set $Y=L^2(\omega)$ and let $B=-\Delta_\omega$ be the positive Laplacian operator with Dirichlet boundary condition in $\omega$. 
Introduce a bounded vector field $w \in \Cc^0_b(\overline\Omega,\RR^3)$, set 
$L(x):(u,v)\in Y^4 \mapsto w\cdot v \in Y$ and $F(x,u)\equiv g(x,y,u)-c(x,y)u$ where $g\in\Cc^1_b(\overline\Omega\times \RR,\RR)$ is a smooth bounded function with bounded derivatives and $c\in \Cc^0_b(\overline \Omega,\RR)$ a bounded coefficient. Then, our abstract equation \eqref{eq_appli} becomes
\begin{equation}\label{eq_appli2}
\partial_t u =\Delta u + w\cdot \grad u-c(x,y)u + g(x,y,u).
\end{equation}
This type of PDE is studied for example in \cite{Berestycki_Nirenberg,Giletti}.
Notice that the function $F$ is a well-defined globally Lipschitz function in the sense of Corollary \ref{coro_globLip} with $\alpha=0$ and any $\beta>0$. 
Thus, Theorem \ref{th_Cauchy_ul_w} applies and we obtain that the advective reaction-diffusion PDE \eqref{eq_appli2} is well posed for any initial data in a space $D_{A_\ulw}(\beta,\infty)$ with $\beta\in(0,1)$, in particular for $u_0\in L^2_\ulw(\RR,H^1_0(\omega))\cap H^1_\ulw(\RR,L^2(\omega))$. We enhance that initial data as $u_0(x,y)=\chi_{x\in ]a,b[}$ or Heaviside functions are included in $D_{A_\ulw}(\beta,\infty)$ for small enough $\beta>0$ since these functions belong to $H^s_\loc(\Omega)$ for $s<1/2$. We also recall that, even if $f$ is very simple, we need to consider initial data slightly more regular than $L^2_\ulw(\RR,L^2(\omega))$ (in both directions) because otherwise, even for $f=0$, the Cauchy problem is ill-posed (see Figures \ref{fig-heat1} and \ref{fig-heat2}).

Let us apply Theorem \ref{th_gradient}. We can keep part of the convection with the following change of setting. Let $Y=L^2(\omega)$ endowed with the scalar product $\la u|\tilde u\ra_Y=\int_\omega e^{W(y)}u(y)\tilde u(y)\d y$ for a $W\in\Cc^1(\overline\omega,\RR)$. We set $B=-\Delta_\omega - \grad_y W(y)\cdot\nabla$. We can check that $B$ is positive self-adjoint on $Y$. 
To fit the context of Theorem \ref{th_gradient}, let us also assume that $g$ and $c$ do not depend on $x$. This leads to the equation 
\begin{equation}\label{eq_appli3}
\partial_t u =\Delta u + w\cdot \grad u - c(y)u + g(y,u)~~~\text{ with }w(y)=\grad_y W(y)
\end{equation}
We have $f(u)=g(u)-cu=-\grad_u V(u)$ in $Y$, in the sense of Theorem \ref{th_gradient}, with 
$$V(u)=-\int_\omega e^{W(y)} \left(\int_0^{u(y)}f(y,\nu)\d \nu\right)\d y$$ 
and the coercivity of $V$ means that $\inf_\omega c(y)>0$. Then Theorem \ref{th_gradient} shows that all the solutions of \eqref{eq_appli3} are global and remains bounded. Moreover, Proposition \ref{prop_compact} shows that the trajectories are asymptotically compact in all finite section of the cylinder.

To our knowledge, this paper is the first one to propose a suitable setting for this Cauchy problem of a parabolic equation in a cylinder in the context of uniformly local spaces. We are not claiming that it cannot be shown by a simpler non-abstract proof. However, we hope that the present paper makes the following facts clearer:
\begin{itemize}
\item Setting a Cauchy problem in spaces as $L^2(\Omega)$ misses most of the interesting patterns with infinite energy, as travelling fronts. If we consider a functional spaces with too much regular functions, we cannot consider the classical Heaviside-like initial data. Moreover, the local $H^1$-topology is the natural one to be able to use the formal energy, so working with a $L^\infty$-topology is not convenient. This leads to  considering the uniformly local Sobolev spaces.
\item The linear part of Equation \eqref{eq_appli2} does not generate a classical semigroup in the weak version $L^2_\ulw(\RR,L^2(\omega))$ of the uniformly local Lebesgue space due to a lack of continuity when $t\rightarrow 0^+$. This is due to the lack of density of the domain of the corresponding linear operator. We would meet the same problem when working in $L^\infty$ or $(\Cc^0,\|\cdot\|_\infty)$.
\item Contrarily to the real-valued case, working in the strong version $L^2_\uls(\RR,L^2(\omega))$ of the uniformly local spaces does not help: the same problems remain. 
\item Fortunately, there exists a suitable theory for this type of ``abstract parabolic operators'' with non-dense domain. It enables to set a relevant Cauchy theory in all the settings, even in the weak versions of the space.
\end{itemize}

\appendix

\section{The Bochner integral and spaces of vector-valued functions}\label{section_spaces}

Let $\Omega$ be a smooth open subset of $\RR^d$ and $Y$ a Banach space. The Bochner integral is the classical way to define an integral $\int_\Omega u(x)\d x$ where $u:\Omega\rightarrow Y$ is a vector-valued function. We recall here the basic facts, the interested reader can find more details in the original article \cite{Bochner} or in many textbooks as \cite{Arendt-Batty-et-al,Diestel-Uhl,Dunford-Schwartz}. A {\it simple} function is a function of the form $u(x)=\sum_{k=1}^n y_k \Un_{x\in A_k}$, where $(A_k)$ is a finite decomposition of the domain $\Omega$ in measurable sets. The notion of integral of a simple function is natural: $\int_\Omega u(x)\d x=\sum_{k=1}^n y_k|A_k|$. More generally, a {\it measurable} function is the almost everywhere limit of a sequence of simple functions. 
A measurable function $u$ is {\it Bochner-integrable} if there exists a sequence $(u_n)$ of simple functions such that $\int_\Omega |u(x)-u_n(x)|_Y \d x\rightarrow 0$ in $\RR$ and the integral of $u$ is then defined as the limit of the ones of $u_n$. 

Let $p\in [1,+\infty]$ and let $\rho\in\Cc^0(\Omega,\RR_+)$ be a non-negative smooth weight. We introduce the Banach space $L^p_\rho(\Omega,Y)$ of measurable functions $u:\Omega\rightarrow Y$ with finite $\|u\|_{L^p_\rho(\Omega,Y)}$ norm, as in Section \ref{section_basic_spaces}.
 Several natural properties of the classical integral extend to the Bochner's integral. We refer to \cite{Arendt-Batty-et-al,Diestel-Uhl,Kreuter} and state the classical results without proof.
\begin{theorem}[Bochner's theorem]\label{th_Bochner}
Let $\Omega\subset \RR^d$ and $\rho\in\Cc^0(\Omega,\RR_+)$ and let $Y$ be a Banach space. A measurable function $u:\Omega\rightarrow Y$ is Bochner-integrable if and only if the real function $|u|_Y:\Omega\rightarrow \RR$ is integrable. In the affirmative case, we have
\begin{equation}\label{estimation_bochner}
\left| \int_\Omega \rho(x) u(x)\d x\right|_Y \leq \int_\Omega \rho(x) |u(x)|_Y \d x. 
\end{equation}
In particular, if $p\in [1,+\infty]$, a measurable function $u:\Omega\rightarrow Y$ belongs to $L^p_\rho(\Omega,Y)$ if and only if $|u|_Y:\Omega\rightarrow \RR$ belongs to $L^p_\rho(\Omega,\RR)$.
\end{theorem}

\begin{prop}[Fubini's theorem]\label{prop_fubini}
Let $\Omega_1\subset \RR^{d_1}$ and $\Omega_2\subset \RR^{d_2}$ be two domains and let $f:\Omega_1\times\Omega_2\rightarrow Y$ be measurable, the classical integral $\int_{\Omega_1\times\Omega_2} |f(x,x')|_Y \d x \d x'$ being finite. Then $f$ is Bochner integrable and 
$$\int_{\Omega_1\times\Omega_2} f(x,x')\d x \d x' =\int_{\Omega_1}\left(\int_{\Omega_2} f(x,x')\d x' \right)\d x =\int_{\Omega_2}\left(\int_{\Omega_1} f(x,x')\d x \right)\d x'.$$ 
\end{prop}

\begin{prop}[Dominated convergence]\label{prop_cv_dom}
Let $f_n : \Omega \rightarrow Y$ be Bochner integrable functions. Assume that $f(x):= \lim_{n\rightarrow +\infty} f_n (x)$ exists a.e. and there exists an integrable function $g :\Omega\rightarrow \RR$ such
that $|f_n (x)|_Y \leq g(x)$ a.e. for all $n\in\NN$. Then $f$ is Bochner integrable and 
$$\lim_{n\rightarrow +\infty} \int_\Omega f_n(x) \d x = \int_\Omega f(x) \d x$$
Actually, we have the stronger convergence $\int_\Omega |f_n(x)-f(x)|_Y \d x\longrightarrow 0$.
\end{prop}

\begin{prop}[Commutations]\label{prop_commute}
If $Y$ and $Z$ are Banach spaces and if $B:Y\rightarrow Z$ is a bounded linear operator, and if $f:\Omega\rightarrow Y$ is Bochner integrable, then $B\circ f$ is Bochner integrable and $\int B\circ f = B (\int f)$.

For all $\phi\in Y^*$ and $f:\Omega\rightarrow Y$ Bochner integrable, $\int \langle f(x)|\phi\rangle_{Y,Y^*}$ is integrable (in the classical sense) and $\int \langle f(x)|\phi\rangle_{Y,Y^*}=\langle \int f(x)|\phi\rangle_{Y,Y^*}$

If $A:D(A)\subset Y\rightarrow Z$ is a closed linear operator, and if $f:\Omega\rightarrow Y$ is Bochner integrable with $f(x)\in D(A)$ for all $x\in\Omega$ and if $A\circ f$ is Bochner integrable, then $\int f$ belongs to $D(A)$ and $\int A\circ f = A(\int f)$.
\end{prop}

\begin{prop}[H\"older's inequalities]\label{prop_Holder}
Let $1\leq p\leq \infty$, set $q=1/(1-1/p)$ ($q=\infty$ if $p=1$). For all $u\in L^p_\rho(\Omega,Y)$ and $\lambda \in L^q_\rho(\Omega,\RR)$, then $\lambda u\in L^1_\rho(\Omega,Y)$ and 
$$\|\lambda u\|_{L^1_\rho(\Omega,Y)} \leq \|\lambda\|_{L^q_\rho(\Omega,\RR)} \|u\|_{L^p_\rho(\Omega,Y)}.$$ 
In addition, for all $u\in L^p_\rho(\Omega,Y)$ and $v \in L^q_\rho(\Omega,Y^*)$, then $\langle v|u\rangle_{Y^*,Y} \in L^1_\rho(\Omega,\RR)$ and 
$$\left| \int_\Omega \rho(x)\langle v|u\rangle_{Y^*,Y} \d x\right| \leq \int_\Omega \rho(x) |\langle v|u\rangle_{Y^*,Y}| \d x  ~\leq~ \|u\|_{L^p_\rho(\Omega,Y)} \|v\|_{L^q_\rho(\Omega,Y^*)}~.$$
\end{prop}
The last result shows that $L^q(\Omega,Y^*)$ can be embedded in $L^p(\Omega,Y)^*$. When dealing with general Banach spaces $Y$, the dual of $L^p(\Omega,Y)$ is not necessarily $L^q(\Omega,Y^*)$, even for finite $p$ and $q$. Radon-Nikodym property (see Definition \ref{defi_RNP}) is required to obtain this natural duality as shown in \cite{Bochner-Taylor}.

We can define the notion of weak derivative and the vector-valued Sobolev spaces $W^{k,p}_\rho(\Omega,Y)$ as in Section \ref{section_basic_spaces}.
We can recover the usual properties of Sobolev spaces. 
\begin{prop}\label{prop_dense_4}
Let $Y$ be a Banach space and let $k\in\NN$ and $p\in [1,+\infty)$. The space $\Cc^\infty_c(\RR^d,Y)$ of smooth and compactly supported functions is dense in $W^{k,p}(\RR^d,Y)$.
\end{prop}
\begin{demo}
To prove this result, it is sufficient to use the mollification and the fact that, if $u$ belongs to $W^{1,p}(\RR^d,Y)$, then 
$$\partial_{x_i}\left(\eta_\varepsilon\ast u\right)=\eta_\varepsilon\ast (\partial_{x_i}u)=(\partial_{x_i}\eta_\varepsilon)\ast u.$$
Indeed, using the definition of the weak derivative and Fubini theorem, we have that, for any test function $\varphi\in\Cc^\infty_c(\RR^d,\RR)$,
\begin{align*}
-\int_{\RR^d} (\eta_\varepsilon\ast u)(x) \partial_{x_i}\varphi(x) \d x &= -\int_{\RR^d}\int_{\RR^d} \eta_\varepsilon(x-x')u(x')\partial_{x_i}\varphi(x) \d x \d x' \\
& = \int_{\RR^d}\int_{\RR^d} \partial_{x_i}\eta_\varepsilon(x-x')u(x')\varphi(x) \d x \d x'.\\
& = -\int_{\RR^d}\int_{\RR^d} \partial_{x'_i}\eta_\varepsilon(x-x')u(x')\varphi(x) \d x \d x'.\\
&=\int_{\RR^d}\int_{\RR^d} \eta_\varepsilon(x-x')\partial_{x'_i}u(x')\varphi(x) \d x \d x'.
\end{align*}
This shows that $\eta_\varepsilon\ast u$ is a smooth function converging to $u$ in $W^{1,p}(\RR^d,Y)$ when $\varepsilon\rightarrow 0$. To obtain a compactly supported function, we can first truncate $u$ inside a sufficiently large ball with a sufficiently smooth truncation with small derivatives. Finally we notice that we can iterate the argument to higher derivatives.
\end{demo}

The integration by parts follows from the definition of the weak derivative and can be seen as its extension to vector-valued test functions. 
\begin{prop}[Integration by parts]\label{prop_inte_parts}
Let $Y$ be a Banach space and let $p$ and $q$ in $(1,+\infty)$ with $1/p+1/q=1$. Assume that $u\in W^{1,p}(\RR^d,Y)$ and $v\in W^{1,q}(\RR^d,Y^*)$. Then $\langle u|v\rangle_{Y,Y^*}$ belongs to $W^{1,1}(\RR^d,\RR)$ and, for any direction $x_i$ in $\RR^d$, we have that 
\begin{equation}\label{eq_inte_part_2}
\int_{\RR^d}  \langle \partial_{x_i} u(x)|v(x)\rangle_{Y,Y^*} \d x~=~ -\int_{\RR^d} \langle u(x)| \partial_{x_i} v(x)\rangle_{Y,Y^*} \d x.
\end{equation}
\end{prop}
\begin{demo}
First notice that the considered integral makes sense by Hölder inequality and the commutation between the integral and the dual product, see Propositions \ref{prop_commute} and \ref{prop_Holder}. Let $\eta_\varepsilon$ be a mollifier. The weak derivative is only defined via real-valued test functions. To extend it to our desired calculus, we write, using Propositions \ref{prop_fubini}, \ref{prop_commute} and \ref{prop_Holder}
\begin{align*}
 \int_{\RR^d}  \langle \partial_{x_i} u(x)|(\eta_\varepsilon\ast v)(x)\rangle_{Y,Y^*} \d x& =
  \int_{\RR^d}\int_{\RR^d}  \langle \partial_{x_i} u(x)|\eta_\varepsilon(x-x')v(x')\rangle_{Y,Y^*} \d x \d x'\\
&  =-\int_{\RR^d}\int_{\RR^d}  \partial_{x_i}\eta_\varepsilon(x-x')\langle u(x)|v(x')\rangle_{Y,Y^*} \d x \d x'\\
&  =\int_{\RR^d}\int_{\RR^d}  \partial_{x'_i}\eta_\varepsilon(x-x')\langle u(x)|v(x')\rangle_{Y,Y^*} \d x \d x'\\
&  =-\int_{\RR^d}\int_{\RR^d}  \eta_\varepsilon(x-x')\langle u(x)|\partial_{x'_i}v(x')\rangle_{Y,Y^*} \d x \d x'\\
&  =-\int_{\RR^d}\langle  u(x)|(\eta_\varepsilon\ast \partial_{x_i}v)(x)\rangle_{Y,Y^*} \d x\\
&  =-\int_{\RR^d}\langle  u(x)|\partial_{x_i}(\eta_\varepsilon\ast v)(x)\rangle_{Y,Y^*} \d x
\end{align*}
Then, it is sufficient to take the limit $\varepsilon\rightarrow 0$.
\end{demo}
\begin{coro}\label{coro_H2}
Assume that $Y$ is a Hilbert space. Then, for all $u\in H^2(\RR^d,Y)$, we have 
$$\|u\|_{H^1(\RR^d,Y)}^2 \leq \sqrt{d+1} \|u\|_{L^2(\RR^d,Y)} \|u\|_{H^2(\RR^d,Y)}~.$$
\end{coro}
\begin{demo}
Due to Proposition \ref{prop_inte_parts}, we have 
\[\sum_{i=1}^d \int_{\RR^d} \langle \partial_{x_i} u| \partial_{x_i} u  \rangle_{Y}= - \int_{\RR^d} \langle \Delta u| u  \rangle_{Y}.\]
Using the Cauchy-Schwarz version of H\"older inequality of Proposition \ref{prop_Holder}, we obtain 
\[ \sum_{i=1}^d \|\partial_{x_i} u\|_{L^2(\RR^d,Y)}^2 \leq \|u\|_{L^2(\RR^d,Y)} \|\Delta u\|_{L^2(\RR^d,Y)}~.\]
Thus, 
\begin{align*}
\|u\|_{H^1(\RR^d,Y)}^2 &\leq \|u\|_{L^2(\RR^d,Y)}^2 + \|u\|_{L^2(\RR^d,Y)} \big( \sum_i \|\partial^2_{x_i,x_i} u\|_{L^2(\RR^d,Y)} \big)\\ & \leq \|u\|_{L^2(\RR^d,Y)}^2 + \sqrt{d} \|u\|_{L^2(\RR^d,Y)} \sqrt{\sum_i \|\partial^2_{x_i,x_i} u\|_{L^2(\RR^d,Y)}^2}\\
& \leq \sqrt {d+1}  \|u\|_{L^2(\RR^d,Y)} \sqrt{\sum_i \|\partial^2_{x_i,x_i} u\|_{L^2(\RR^d,Y)}^2 + \|u\|_{L^2(\RR^d,Y)}^2}~.
\end{align*}
\end{demo}
 
The Sobolev embeddings are the same as the classical ones. It is noteworthy that the Radon-Nikodym property is not needed here, as shown in Section 6 of \cite{Arendt-Kreuter}.  
\begin{theorem}[Sobolev embeddings]\label{th_embeddings}
Let $\Omega$ be an open smooth subset of $\RR^d$ and $Y$ be a Banach space. Then, 
\begin{itemize}
\item[(i)] if $1<p<d$, then $W^{1,p}(\Omega,Y)\hookrightarrow L^{p^*}(\Omega,Y)$ where $p^*$ satisfies $\frac 1{p^*}=\frac 1p - \frac 1d$.
\item[(ii)] if $p=d$, then $W^{1,p}(\Omega,Y)\hookrightarrow L^{q}(\Omega,Y)$ for all $q\in [p,+\infty)$.
\item[(iii)] if $p>d$ then $W^{1,p}(\Omega,Y)\hookrightarrow L^{\infty}(\Omega,Y)$.
\end{itemize}
\end{theorem}
Notice that the embedding $W^{1,p}(\Omega,Y)\rightarrow L^\infty(\Omega,Y)$ for $p>d$ and Corollary \ref{coro_not_dense} shows that we cannot expect in general $W^{1,p}(\Omega,Y_2)$ to be dense in $W^{1,p}(\Omega,Y_1)$ even if $Y_2$ is dense in $Y_1$.


\section{Sectorial operators}\label{section_sectoriel}
We recall here some basic facts and estimates related to sectorial operators and analytic semigroups. Doing so, we enhance some proofs or formula, in which we will need to track later the dependence or independence with respect to some parameters. 
The missing definitions and proofs can be found in many textbooks, see for example \cite{Henry,Magal-Ruan}.
\begin{defi}\label{defi_secto}
A linear operator $A$ in a Banach space $X$ is called a {\bf sectorial operator} if it is a closed densely defined operator such that there exists $z_0\in\RR$ and $\phi\in (0,\pi/2)$ such that the sector 
$$S_{z_0,\phi}=\{ z \in \CC\setminus\{z_0\}~,~~ \phi \leq |\arg(z-z_0)| \leq \pi\}$$
is in the resolvent set of $A$ and there exists $M>0$ such that 
$$\forall z\in S_{z_0,\phi}~,~~\|(z\id-A)^{-1}\|_{\Lc(X)} \leq \frac {M}{|z-z_0|}~. $$
\end{defi}
\begin{prop}\label{prop_sa_secto1}
If $A$ is a densely defined self-adjoint operator in a Hilbert space $H$ and if $A$ is bounded from below, then $A$ is sectorial.
\end{prop} 
\begin{prop}\label{prop_sa_secto2}
Let $A$ be sectorial operator on $X$. There exists $\varepsilon>0$ such that, if $K$ is a linear operator with $D(A)\subset D(K)$ and if there exists $C>0$ such that
\[\forall u\in D(A)~,~~\|Ku\|\leq \epsilon \|Au\| + C\|u\|,\]
then $A+K:D(A)\rightarrow X$ is sectorial.
\end{prop}
\begin{demo}
See \cite{Henry}, Theorem 1.3.2 and the examples and remarks of p.19. We simply emphasize the key estimates for the perturbative argument are
\[\|K(z\id-A)^{-1}\|_{\Lc(X)} \leq \epsilon \|A(z\id-A)^{-1}\|_{\Lc(X)} + C \|(z\id-A)^{-1}\|_{\Lc(X)} \]
and 
\[\|(z\id-(A+K))^{-1}\|_{\Lc(X)}=\|(z\id-A)^{-1} (\id -K(z\id-A)^{-1})^{-1}\|_{\Lc(X)}.\] 
From them, we can see that the constants associated to the sectorial property of $A+K$ (as in Definition \ref{defi_secto}) only depends on $A$, $\epsilon$ and $C$. 
\end{demo}

\begin{prop}\label{prop_analytic_semigroup}
A sectorial operator $A$ in a Banach space $X$ generates an analytic semigroup $e^{-A t}$ on $X$ given by 
\begin{equation}\label{defi_semigroup2}
e^{-At}=\frac 1{2i\pi}\int_\Gamma (\lambda+A)^{-1}e^{\lambda t}\d \lambda
\end{equation}
where $\Gamma$ is a suitable contour.
\end{prop}

If $A$ is a sectorial operator for which the sector $S_{a,\phi}$ can be chosen with $a>0$, then we can define the fractional powers $A^\alpha$ of $A$ as the inverse of the operator
\begin{equation}\label{defi_fracto}
A^{-\alpha}=\frac 1{\Gamma(\alpha)} \int_0^\infty t^{\alpha-1}e^{-At}\d t.
\end{equation}

\section{Interpolation}\label{section_interpolation}

Let $E_0$ and $E_1$ (resp. $F_0$ and $F_1$) be two Banach spaces included in a larger one $\Ec$ (resp. $\Fc$). An interpolation theory is a way to define a scale of intermediate spaces $(E_\theta)_{\theta\in(0,1)}$ between $E_0\cap E_1$ and $E_0+E_1$, such that the following property holds: if $L:\Ec\rightarrow \Fc$ is a linear operator, whose restrictions $L:E_0\rightarrow F_0$ and $L:E_1\rightarrow F_1$ are continuous, then $L:E_\theta \rightarrow F_\theta$ is also continuous. The most usual ways to define such scale of intermediate spaces are more or less equivalent to one of these two methods:
\begin{itemize}
\item the {\it real interpolation}, introduced by Lions and Peetre in \cite{Lions,Lions-Peetre} as ``espaces de traces et de moyennes''. Consider $p\in [1,+\infty]$ and $\theta\in (0,1)$, the intermediate space $E_\theta$ obtained by the real interpolation method is denoted in this paper $$(E_0,E_1)_{\theta,p}$$
\item the {\it complex interpolation}, introduced by Calder\'on, Lions and Kre\u{\i}n in \cite{Calderon,Lions2,Krein} as ``analytic scales of Banach spaces''. Consider $\theta\in (0,1)$, the intermediate space $E_\theta$ obtained by the complex interpolation method is denoted in this paper $$[E_0,E_1]_{\theta}$$ 
\end{itemize}
The definitions of these methods are technical and will not be recalled in this paper, we refer to \cite{Bergh-Lofstrom}, \cite{Triebel} or the above historical references. Both methods of interpolation are related as follows (see \cite{Bergh-Lofstrom} or Section I.2 of \cite{Amann-livre}). 
\begin{prop}\label{prop_interpolations}
Let $E_1\hookrightarrow E_0$ be two Banach spaces included in a larger one. We have the chain of inclusions
\begin{align*}
E_0 \cap E_1 &\hookrightarrow (E_0,E_1)_{\zeta,q} \hookrightarrow (E_0,E_1)_{\eta,1} \hookrightarrow [E_0,E_1]_\eta \\  
& \hookrightarrow (E_0,E_1)_{\eta,\infty} \hookrightarrow (E_0,E_1)_{\xi,q} \hookrightarrow E_0+E_1
\end{align*}
for $1\leq q<\infty$ and $0<\zeta < \eta < \xi < 1$.
\end{prop}

The following classical result is proved in \cite{Lions-Peetre} and \cite{Cwikel}, see also many textbooks on interpolation as \cite{Bergh-Lofstrom}. 
\begin{prop}\label{prop_interpolation_Lp}
Let $\Omega$ be an open subset of $\RR^d$ and let $E_0$ and $E_1$ be two Banach spaces. Let $\theta\in (0,1)$, $p_0$ and $p_1\in [1,+\infty]$ and set $p_\theta$ such that $\frac 1{p_\theta}=\frac {1-\theta}{p_0} + \frac \theta {p_1}$. Then we have
\[\big(L^{p_0}(\Omega,E_0) , L^{p_1}(\Omega,E_1)\big)_{\theta,p_\theta} = L^{p_\theta}\big(\Omega,(E_0,E_1)_{\theta,p_\theta}\big). \]
\end{prop}

The case of Sobolev spaces is more involved. Even for real valued functions, the interpolated space between two Sobolev spaces is a Besov space and not necessarily a Sobolev space, see \cite{Bergh-Lofstrom}. In this paper, we will use the following embedding result.
\begin{coro}\label{coro_interpolation_Sob}
Let $\Omega$ be an open subset of $\RR^d$ and let $E_0$ and $E_1$ be two Banach spaces. Let $\theta\in (0,1)$, $p_0$ and $p_1\in [1,+\infty]$ and set $p_\theta$ such that $\frac 1{p_\theta}=\frac {1-\theta}{p_0} + \frac \theta {p_1}$. Then we have
\[\big(W^{1,p_0}(\Omega,E_0) , W^{1,p_1}(\Omega,E_1)\big)_{\theta,p_\theta} \hookrightarrow W^{1,p_\theta}\big(\Omega,(E_0,E_1)_{\theta,p_\theta}\big).\]
\end{coro}
\begin{demo}
If $u$ belongs to $\big(W^{1,p_0}(\Omega,E_0) , W^{1,p_1}(\Omega,E_1)\big)_{\theta,p_\theta}$, then Proposition \ref{prop_interpolation_Lp} above shows that $u$ belongs to $L^{p_\theta}\big(\Omega,(E_0,E_1)_{\theta,p_\theta}\big)$. Then, we notice that, for any $i=1,\ldots,d$, the linear operator $\partial_{x_i}$ defines a continuous operator from $W^{1,p_j}(\Omega,E_j)$ into $L^{p_j}(\Omega,E_j)$ ($j=0,1$). The interpolation property for continuous linear operators and Proposition \ref{prop_interpolation_Lp} shows that $\partial_{x_i}$ extends to a continuous operator from $\big(W^{1,p_0}(\Omega,E_0) , W^{1,p_1}(\Omega,E_1)\big)_{\theta,p_\theta}$ into  $L^{p_\theta}\big(\Omega,(E_0,E_1)_{\theta,p_\theta}\big)$. 
\end{demo}

The problem of mixing the Sobolev regularities with interpolating spaces is one of the problems of ``maximal regularity''. We will use here the following result
\begin{prop}\label{prop_interpolation_Amann}
Let $E_0$ and $E_1$ be two Banach spaces such that $E_1\hookrightarrow E_0$ densely. Then, for any $p\in (1,\infty)$, 
\[L^p(\RR,E_1)\cap W^{1,p}(\RR,E_0)\,\hookrightarrow\,L^\infty(\RR,(E_0,E_1)_{1-1/p,p})\]
\end{prop}
\begin{demo}
We refer to Theorem III.4.10.2 of \cite{Amann-livre}. Since we consider the $L^\infty-$norm, up to using a partition of the unity, we can extend the result from a subinterval of $\RR$ to the whole space $\RR$.
\end{demo}

The interpolations between domains of fractional powers of sectorial operators are proved in Triebel \cite[Theorems 1.15.2 and 1.18.10]{Triebel} (see also \cite{Amann-livre}). 
\begin{prop}\label{prop_interpolation_autoadjoint}
Let $X$ be a Banach space and let $A$ be a densely defined sectorial operator from $D(A)$ into $X$. Then, for $\alpha\in (0,1)$,
\[(X,D(A))_{\alpha,1} \hookrightarrow  (D(A^\alpha)) \hookrightarrow (X,D(A))_{\alpha,\infty}~.\]
If $Y$ is a Hilbert space and $B$ a positive self-adjoint operator on $Y$, then
\[D(B^\alpha) = [Y,D(B)]_\alpha = (Y,D(B))_{\alpha,2}.\]
\end{prop}


\end{document}